\documentclass[preprint,12pt]{elsarticle}

\usepackage{comment}
\usepackage{tikz}

\usetikzlibrary{calc}
\usepackage{graphicx}%
\usepackage{multirow}%
\usepackage{amsmath,amssymb,amsfonts}%
\usepackage{amsthm}%
\usepackage{mathrsfs}%
\usepackage[title]{appendix}%
\usepackage{xcolor}%
\usepackage{textcomp}%
\usepackage{manyfoot}%
\usepackage{booktabs}%
\usepackage{algorithm}%
\usepackage{algorithmicx}%
\usepackage{algpseudocode}%
\usepackage{listings}%
\usepackage{url}
\usepackage{hyperref}
\usepackage[acronym,nonumberlist]{glossaries}

\newcommand{\R}[1]{}
\newcommand{\revone}[1]{\textcolor{black}{#1}}
\newcommand{\revtwo}[1]{\textcolor{black}{#1}}

\DeclareMathOperator*{\argmin}{arg\,min}

\makenoidxglossaries   

\newacronym{LSTM}{LSTM}{long short-term memory}
\newacronym{BC}{BC}{boundary condition}
\newacronym{ROM}{ROM}{reduced-order model}
\newacronym{DOF}{DOF}{degrees of freedom}
\newacronym{FOM}{FOM}{full-order model}
\newacronym{POD}{POD}{proper orthogonal decomposition}
\newacronym{PDE}{PDE}{partial differential equation}
\newacronym{RK4}{RK4}{Runge-Kutta 4}
\newacronym{G-POD}{G-POD}{space-global proper orthogonal decomposition}
\newacronym{L-POD}{L-POD}{space-local proper orthogonal decomposition}
\newacronym{LO-POD}{LO-POD}{space-local proper orthogonal decomposition with overlapping subdomains}
\newacronym{SVD}{SVD}{singular value decomposition}

\newcommand{\bPhi}{\boldsymbol{\Phi}}
\newcommand{\bGamma}{\boldsymbol{\Gamma}}

\journal{Computers \& Fluids}

\begin{document}

\begin{frontmatter}



\title{Modeling Advection-Dominated Flows with Space-Local Reduced-Order Models}


\author[1,eindhoven]{T. van Gastelen}
\author[1]{W. Edeling}
\author[1,eindhoven]{B. Sanderse}


\address[1]{Centrum Wiskunde \& Informatica, 
            Science Park 123, 
            Amsterdam,
            The Netherlands}
\address[eindhoven]{Centre for Analysis, Scientific Computing and Applications, Eindhoven University of Technology,
    PO Box 513, Eindhoven,
    5600 MB, The Netherlands}

\begin{abstract}
\Glspl{ROM} are often used to accelerate the simulation of large physical systems. However, traditional \gls{ROM} techniques, such as \gls{POD}-based methods, often struggle with advection-dominated flows due to the slow decay of singular values. This results in high computational costs and potential instabilities.

This paper proposes a novel approach using space-local \gls{POD} to address the challenges arising from the slow singular value decay. Instead of global basis functions, our method employs local basis functions that are applied across the domain, analogous to the finite element method, but with a data-driven basis. By dividing the domain into subdomains and applying the space-local \gls{POD}, we obtain a sparse representation that generalizes better outside the training regime. This allows the use of a larger number of basis functions compared to standard \gls{POD}, without prohibitive computational costs. To ensure smoothness across subdomain boundaries, we introduce overlapping subdomains inspired by the partition of unity method.

Our approach is validated through simulations of the 1D and 2D advection equation. We demonstrate that using our space-local approach, we obtain a \gls{ROM} that generalizes better to flow conditions not included in the training data. In addition, we show that the constructed \gls{ROM} inherits the energy conservation and non-linear stability properties from the full-order model. Finally, we find that using a space-local \gls{ROM} allows for larger time steps.
\end{abstract}


\begin{highlights}
\item Novel space-local POD approach for advection-dominated flows.
\item Local basis improves generalizability and sparsity of ROMs. 
\item Overlapping subdomains ensure smooth solutions.
\item ROMs conserve energy and allow larger simulation time steps
\end{highlights}


\end{frontmatter}

\section{Introduction}\label{sec:introduction}

Simulating large physical systems is an ongoing challenge in the field of computational sciences. This especially becomes challenging when dealing with multiscale systems. Such systems exhibit interesting behavior at various spatial and temporal scales. A prominent example and our main incentive for this work are turbulent flows described by the incompressible Navier-Stokes equations. Systems described by these equations feature the formation of turbulent eddies of a range of different sizes. The significant difference in size between the largest and the smallest eddies gives rise to the multiscale nature of the problem. The problem with simulating such systems is that they require high-resolution computational meshes to obtain accurate simulations and small time steps. This places a significant burden on the available computational resources \cite{sasaki2002navier, agdestein2024discretizefirstfilternext}.


To make simulation of turbulent flows feasible, one typically resorts to Reynolds averaged Navier-Stokes \cite{RANS}, large eddy simulation \cite{sagaut2006large}, and \acrfullpl{ROM} \cite{brunton2019data,projection_ROMs_survey}. Here we focus on the latter approach. In reduced-order modelling, data is used to speed up simulations. The data can be obtained from simulations or experiments. The collected data is then used to identify the most critical features of the flow. This is typically done by a \acrfull{POD} of the collected flow data. The resulting features are then used to construct a reduced basis. By projecting the fluid flow equations on this basis one obtains the widely used \gls{POD}-Galerkin \gls{ROM} \cite{POD1, POD2, projection_ROMs_survey}. However, two common issues with \gls{POD}-Galerkin \glspl{ROM} are their stability and their accuracy for convection-dominated systems. In \cite{podbenjamin}, it is shown that the stability issue can be resolved by making sure that the energy-conserving property of the Navier-Stokes equations is still satisfied under Galerkin projection. For the incompressible Navier-Stokes equations, the only prerequisite for an energy-conserving \gls{ROM} is that the discretization is structure-preserving, i.e.\ it is such that it inherits the energy conservation property from the continuous equations. 

However, the stability property derived in \cite{podbenjamin} does not guarantee accuracy; it is possible to have energy-stable simulations that are highly inaccurate. This issue can already be observed on academic test cases such as the linear advection equation. For example, in \cite{brunton2019data}, the authors discuss the problem of a single traveling wave through a periodic domain. Even though it is clear that a one-dimensional representation of the system exists, it is not recovered by the \gls{POD} algorithm. The result is a slow decay in singular values of the snapshot matrix. This means that many \gls{POD} modes are required to accurately describe the flow. The slow decay in singular values is often related to the Kolmogorov $N$-width. This is a measure of how well the solution space of a \gls{PDE} can be represented by a linear combination of $N$ basis functions \cite{advection_autoencoder_and_n_width, n_width}. Advection-dominated flows are notorious for displaying a slow Kolmogorov $N$-width decay, requiring a large $N$ for accurate simulation. The resulting \gls{ROM} is then expensive to evaluate \cite{expensive_1, expensive_2}. When not including a sufficient number of basis functions, inaccurate results are obtained, e.g., they contain oscillations \cite{oscillations_1, oscillations_2}. 

Different approaches to deal with this problem have been suggested. 
One way of dealing with this is by taking a local in time approach \cite{time_local_farhat}. In this approach, the \gls{ROM} switches basis for different time intervals. This means that the \gls{ROM} can employ a smaller basis, since each basis is specialized to deal with a specific time interval. Switching between bases during simulation can also be done in a structure-preserving (energy-conserving) manner \cite{KLEIN2024112697}. A related approach is to update the basis during the \gls{ROM} simulation to make it more robust to changes in simulation conditions \cite{adaptive_basis}.
In \cite{Grimberg}, a Petrov-Galerkin approach is suggested, which leads to more stable \glspl{ROM} when using a small \gls{POD} basis than the standard Galerkin approach.
In \cite{ROM_closure}, it is suggested to add a closure model to the \gls{ROM} to account for a small \gls{POD} basis. The closure model is then tasked with modeling the interaction between the part of the flow that is covered by the \gls{POD} basis and the part that is removed. 
Another suggestion is to construct \glspl{ROM} on non-linear reduced subspaces instead of linear ones. In \cite{klein2024entropystablemodelreductiononedimensional}, a more accurate representation of the solution space is obtained using rational quadratic manifolds. Here, the \gls{ROM} construction was also performed in a structure-preserving (entropy-stable) manner, yielding stability of the \gls{ROM}.

Machine learning approaches have also gained traction for the construction of \glspl{ROM}. In \cite{POD_LSTM}, a \gls{LSTM} neural network is used instead of Galerkin projection, as the computational complexity of \gls{LSTM} is more favorable than Galerkin projection-based \glspl{ROM}. They also allow one to take larger time steps \cite{MUCKE2021101408}. In \cite{MUCKE2021101408}, this approach is combined with an autoencoder to reduce the dimensionality of the system. It is demonstrated that a significantly greater reduction in \gls{DOF} can be achieved using this approach compared to the standard \gls{POD}-Galerkin approach. 
In \cite{advection_autoencoder_and_n_width}, an advection-aware autoencoder is suggested to reduce the dimensionality of this system. This is achieved by introducing an additional decoder which is trained to reconstruct a ``shifted" version of the encoded snapshots. This shift can, for example, be a snapshot taken at a later time. This forces the latent space of the autoencoder to become aware of the dominant advective features of the flow. 

In this work, we focus on an alternative approach to address the issue of \gls{ROM} accuracy in advection-dominated flows. Namely, we propose to use the idea of a space-local \gls{POD} to tackle the slow singular value decay of advection-dominated flows. The idea is as follows: instead of obtaining a set of global basis functions that span the entire domain, we obtain a set of space-local basis functions. These basis functions are repeated across the entire domain, similarly to how a finite element basis covers the domain \cite{FEM_book}. Furthermore, they are only nonzero on their designated subdomain. The advantage of this is that the resulting Galerkin projected operators are sparse. This makes them much cheaper to evaluate when simulating the \gls{ROM}. For example, in \cite{farhat_local_POD} a speedup of up to 1.5 orders of magnitude is reported with respect to a standard \gls{POD}-Galerkin \gls{ROM}. Our approach differs from the one presented in \cite{farhat_local_POD}. In our approach, the solution data in each subdomain is treated with the same local \gls{POD} basis, rather than obtaining a different local \gls{POD} basis for each subdomain. 
This approach ensures that a feature observed in part of the domain can also be represented in a different part of the domain using the same basis. In this way, less data is required to obtain a \gls{POD} basis that generalizes well. 
Note that this approach is designed specifically for problems where the dynamics are similar throughout the domain, such as a traveling wave. For problems where the dynamics are more variable throughout the domain, the approach suggested by \cite{farhat_local_POD} is likely still superior, as it builds a specialized basis for each subdomain. \R{comment_13_a}\revone{For diffusion-dominated problems the standard space-global \gls{POD} approach is likely still superior due to the rapid Kolmogorov $N$-width decay \cite{Quarteroni2015ROM}.} Note that while our approach does not directly solve the slow Kolmogorov $N$-width decay (we still use linear approximations), the sparsity of the basis allows us to use a much larger number of basis functions. Sparsity and generalizability are therefore key features of our approach.
\R{comment_1a}\revone{Furthermore, the authors note this approach is similar to the methodology described in \cite{CHUNG2024DDROM} where different types of subdomains were specified, each with a shared \gls{POD} basis. These subdomains were then used as building blocks to extrapolate the \gls{POD} basis obtained from small spatial domains to larger domains. Our work differs in the fact that it focuses on temporal extrapolation for time-dependent \glspl{PDE}, as opposed to steady-state problems. In addition, we introduce another novel idea in this work: the use of overlapping subdomains (to avoid discontinuities at the subdomain boundaries). Lastly, we show that our space-local \glspl{ROM} satisfy energy conservation if the \gls{FOM} does. In this way, stability of the \gls{ROM} is guaranteed.}



This paper is structured in the following way. 
In Section \ref{sec:FOM} we introduce the \gls{FOM}, namely a central difference discretization of the 1D advection equation. The central difference discretization ensures the scheme is energy conserving. In Section \ref{sec:POD} we introduce the \gls{POD} approaches. We begin with the standard \gls{POD} approach and then introduce two space-local approaches, one with and one without overlapping subdomains. In Section \ref{sec:ROM} we use Galerkin projection to project the \gls{FOM} onto the \gls{POD} basis and show that the resulting \glspl{ROM} satisfy energy conservation. Finally, in Section \ref{sec:results}, we evaluate the different \gls{POD}-based approaches in numerical experiments using a set of different metrics. In addition, we assess the performance of the \glspl{ROM} on a case where we extrapolate beyond the data used to construct the \gls{POD} basis. Furthermore, we investigate the energy-conserving properties of the \glspl{ROM}. To conclude this section, we apply the introduced methodologies to a 2D advection equation test case and evaluate the computational efficiency, among other metrics. Finally, in Section \ref{sec:conclusions}, we present our main findings and suggest future research topics.

\section{Full-order model}\label{sec:FOM}

\subsection{Advection equation}

In this work, we focus on developing a reduced-order model for the linear advection equation in both one-dimensional (1D) and two-dimensional (2D) settings. This equation is chosen as it exhibits similar difficulties as the incompressible Navier-Stokes equations when applying model reduction. For the sake of simplicity, we discuss the properties of the system for 1D, but the ideas easily carry over to 2D. To start, we consider a scalar solution $u(x,t)$ to the linear advection equation
\begin{equation}\label{eq:advection_equation}
    \frac{\partial u}{\partial t} = - c \frac{\partial u}{\partial x},
\end{equation}
with initial condition $u(x,0) = u_0(x)$ and constant $c$. \R{comment_2}\revone{In this work we stick to periodic \glspl{BC}}. An important property of this equation is that the total energy of the system
\begin{equation}\label{eq:energy}
    E := \frac{1}{2}(u,u),
\end{equation}
is conserved, where 
\begin{equation}
    (a(x),b(x)) := \int_\Omega a(x) b(x)  \text{d}\Omega
\end{equation}
on the spatial domain $\Omega$.
This can easily be shown using the product rule of differentiation:
\begin{equation}\label{eq:energy_conservation}
    \frac{\text{d}E}{\text{d}t} =  \frac{1}{2} \frac{\text{d}}{\text{d}t}(u,u) = \left(u, \frac{\partial u }{\partial t}\right)  = - c\left( u, \frac{\partial u}{\partial x}\right) =  c\left(\frac{\partial u}{\partial x} ,u \right)  = 0.
\end{equation}
In the final step, we carried out integration by parts and used the fact that the boundary term cancels on periodic domains. The energy-conserving property of this equation will be mimicked by the discretization and by the \glspl{ROM} developed in this work. This leads to unconditionally stable methods, as in \cite{podbenjamin}.

\subsection{Finite difference discretization}\label{sec:discretization}
Although equation \eqref{eq:advection_equation} can be easily solved exactly, for more complex equations, this is not the case. In general, we approximate the solution by representing $u(x,t)$ on a grid. In this case, we employ a uniform grid with grid spacing $h = \frac{| \Omega|}{N}$ such that $\text{u}_i(t) \approx u(x_i,t)$. The approximated solution is contained within the state vector $\mathbf{u}(t) \in \mathbb{R}^{N}$. To approximate the spatial derivative, we use a central difference approximation:
\begin{equation}\label{eq:stencil}
    \frac{\partial u}{\partial x}|_{x_i} \approx \frac{\text{u}_{i+1}-\text{u}_{i-1}}{2h}.
\end{equation}
This leads to the following semi-discrete system of equations
\begin{equation}\label{eq:semi-discrete}
    \frac{\text{d}\mathbf{u}}{\text{d}t} = -c \mathbf{D}\mathbf{u},
\end{equation}
where the linear operator $\mathbf{D} \in \mathbb{R}^{N \times N}$ is skew-symmetric and encodes the stencil in \eqref{eq:stencil}. For the time integration, we use a classic \gls{RK4} scheme \cite{RK4_Butcher:2007}. This time integration scheme introduces a small energy conservation error, which is negligible in our test cases. Alternatively, one can use the implicit midpoint method for exact energy conservation in time \cite{SANDERSE2013}. Equation \eqref{eq:semi-discrete} will be regarded as the full-order model (\gls{FOM}). We note that other discretization techniques, such as the finite element method, can also be employed to derive a \gls{FOM}; our framework in Section \ref{sec:POD} is also applicable in this context.

\subsection{Energy conservation of the FOM}

It is well known that this \gls{FOM} mimics the energy-conservation property \eqref{eq:energy_conservation} in a discrete setting. 
This can be shown by defining the discretized energy as
\begin{equation}
    E_h := \frac{h}{2}\mathbf{u}^T\mathbf{u},
\end{equation}
\R{comment_3}\revone{where the inclusion of $h$ is such that this definition discretely represents the inner product in \eqref{eq:energy}.}

Using the product rule, we obtain the following evolution equation for the energy
\begin{equation}\label{eq:discr_E_conservation}
    \frac{\text{d}E_h}{\text{d}t} =\frac{h}{2}  \frac{\text{d}\mathbf{u}^T\mathbf{u}}{\text{d}t}  = h \mathbf{u}^T \frac{\text{d}\mathbf{u}}{\text{d}t} =  - c h \mathbf{u}^T \mathbf{D}\mathbf{u}= c h \mathbf{u}^T \mathbf{D}^T\mathbf{u} = 0.
\end{equation}
In the final step, we used the fact that $\mathbf{D}$ is skew-symmetric, i.e., $\mathbf{D} = -\mathbf{D}^T$, to see that the energy is conserved using this stencil. This yields both stability and consistency with the continuous equation.

\section{Local and global POD}\label{sec:POD}

To reduce the computational cost of solving the \gls{FOM}, we construct a \gls{ROM}. For this purpose, we use simulation data to build a data-driven basis. This basis is obtained through a \gls{POD} of the simulation data \cite{brunton2019data,projection_ROMs_survey}. The most common approach is to employ a global \gls{POD} basis, defined over the entire simulation domain, similar to a Fourier basis. An alternative is a local basis, as proposed in \cite{farhat_local_POD,CHUNG2024DDROM}. In this section, we discuss both approaches. We will present our version of the space-local \gls{POD} framework for finite difference discretizations, see \eqref{eq:stencil}. \R{comment_4a}\revone{However, the ideas can also be applied to different discretization techniques, as shown in \cite{CHUNG2024DDROM}.} Our space-local approach has some key differences from the one presented in \cite{farhat_local_POD}, which will be highlighted. 

\subsection{Global POD}

In the global approach we first express the discrete solution $u_h(x,t)$ in terms of an orthogonal basis $\{
\psi_i\}$:
\begin{equation}\label{eq:basis_representation}
        u(x,t) \approx u_h(x,t) = \sum_{i=1}^N \text{c}_{i}(t)\psi_i(x).
\end{equation}
Here $\psi_i$ can represent, for example, a Fourier basis or a localized box function in the case of the presented finite difference discretization. In the finite difference case we have $\mathbf{c}(t)=\mathbf{u}(t)$ and 
\begin{equation}\label{eq:FD_basis}
    \psi_i(x) = \begin{cases}
        1 \quad &\text{if } x_i - \frac{h}{2} \leq x < x_i + \frac{h}{2}, \\
        0 \quad &\text{elsewhere}.
    \end{cases}
\end{equation}
For finite element discretizations, one could use, e.g., L\"{o}wdin orthogonalization to express the solution in terms of an orthogonal basis \cite{aiken1980lowdin}. 

The snapshot matrix $\mathbf{X}$ is constructed from the coefficient vector at different points in time
\begin{equation}\label{eq:global_POD_snapshot}
    \mathbf{X} = \begin{bmatrix}
        \mathbf{c}(t_1) & \mathbf{c}(t_2) & \ldots & \mathbf{c}(t_s)
    \end{bmatrix} \quad \in \mathbb{R}^{N \times s},
\end{equation}
where $s$ is the number of snapshots.
We decompose this snapshot matrix using a \gls{SVD} \cite{SVD_complexity}
\begin{equation}
    \mathbf{X} = \tilde{\boldsymbol{\Phi}}\boldsymbol{\Sigma}\mathbf{V}^T.
\end{equation}
We use the first $r$ left-singular vectors in $\tilde{\boldsymbol{\Phi}} \in \mathbb{R}^{N\times N}$ to obtain a reduced set of orthonormal basis vectors $\boldsymbol{\Phi} \in \mathbb{R}^{N\times r}$.
This basis minimizes the projection error in the Frobenius norm:
\begin{equation}\label{eq:minimal_projection_error}
    \bPhi = \argmin_{\mathbf{M} \in \mathbb{R}^{N\times r}} ||\mathbf{X} - \mathbf{M}\mathbf{M}^T\mathbf{X}||^2_F,
\end{equation}
under the orthonormality constraint $\bPhi^T\bPhi = \mathbf{I}$ \cite{POD1,POD2,brunton2019data,projection_ROMs_survey}.
We can project the coefficient vector onto this basis as follows:
\begin{equation}\label{eq:transform}
    \mathbf{a}(t) = \boldsymbol{\Phi}^T \mathbf{c}(t) \quad \in \mathbb{R}^r,
\end{equation}
such that projecting back onto the \gls{FOM} space yields the approximation 
\begin{equation}
    \mathbf{c}(t) \approx \mathbf{c}_r(t) := \bPhi \mathbf{a}(t).
\end{equation}
By substituting this into \eqref{eq:basis_representation} we obtain the approximated solution $u_r(x,t)$:
\begin{equation}
    u_h(x,t) \approx u_r(x,t)  =  \sum_{i=1}^N \text{c}_{r,i}(t)\psi_i(x).
\end{equation}
We can also write this approximation in terms of the \gls{POD} basis expansion:
\begin{equation}\label{eq:approximation}
\begin{split}
    u_r(x,t)  &=  \sum_{i=1}^N \text{c}_{r,i}(t)\psi_i(x) = \sum_{i=1}^N (\bPhi \mathbf{a}(t))_i \psi_i(x) = \sum_{i=1}^N (\sum_{j=1}^r \Phi_{ij}\text{a}_j(t)) \psi_i(x)  \\ &=  \sum_{i=1}^N \sum_{j=1}^r \Phi_{ij}\text{a}_j(t) \psi_i(x) = \sum_{j=1}^r   \text{a}_j(t) \sum_{i=1}^N \Phi_{ij}\psi_i(x) = \sum_{j=1}^r   \text{a}_j(t) \phi_j(x),
\end{split}
\end{equation}
where the reduced \gls{POD} basis $\{\phi_j\}$ is obtained by applying the following transformation:
\begin{equation}\label{eq:POD_basis}
    \phi_j(x) = \sum^N_{i=1}\Phi_{ij}\psi_i(x).
\end{equation}
As the obtained basis functions span the entire domain, we refer to this approach as \gls{G-POD}.

\subsection{Space-local POD}\label{sec:L-POD}

\R{comment_4b}\revone{In \gls{L-POD} we take a different approach from \gls{G-POD}. We start by assuming a uniform finite difference discretization for the discrete 
\gls{FOM} solution $u_h(x,t)$\footnote{If $u_h(x,t)$ would stem from a simulation on an unstructured grid, one could potentially first project the discrete solution on a uniform grid, after which the presented methodology could still be applied. We consider this outside the scope of this research.}.}
Once the solution is represented on this grid, we subdivide the domain $\Omega$ into $I$ non-overlapping subdomains $\Omega_i = [\alpha_i,\beta_i)$, $i=1\ldots I$, such that $\alpha_i < \beta_i = \alpha_{i+1} < \beta_{i+1}$.
Later in this section, these subdomains will be used to construct the \gls{L-POD} basis. 
Furthermore, we assume each of these subdomains contains exactly $J$ grid points. This means the total number of grid points $N$ has to equal the product $N=I \cdot J$. This is essential for our methodology to work, as for the \gls{L-POD} procedure, each subdomain is treated equally in the snapshot matrix, which is only possible if they contain the same number of grid points. This helps us obtain a basis that generalizes better outside the training data for advection-dominated problems, where the dynamics are similar throughout the domain.
This approach is fundamentally different from what is presented in \cite{farhat_local_POD}, where each subdomain has its own \gls{POD} basis. The latter is more suitable for problems where the dynamics differ throughout the domain, such that the local bases can be specialized to these dynamics. If $N/I$ is not an integer, one could possibly resolve this by projecting the \gls{FOM} solution on a compatible grid of $M$ grid points for which $M/I$ is an integer. 

\begin{figure}[ht]
    \centering
    \includegraphics[width = 0.6\textwidth]{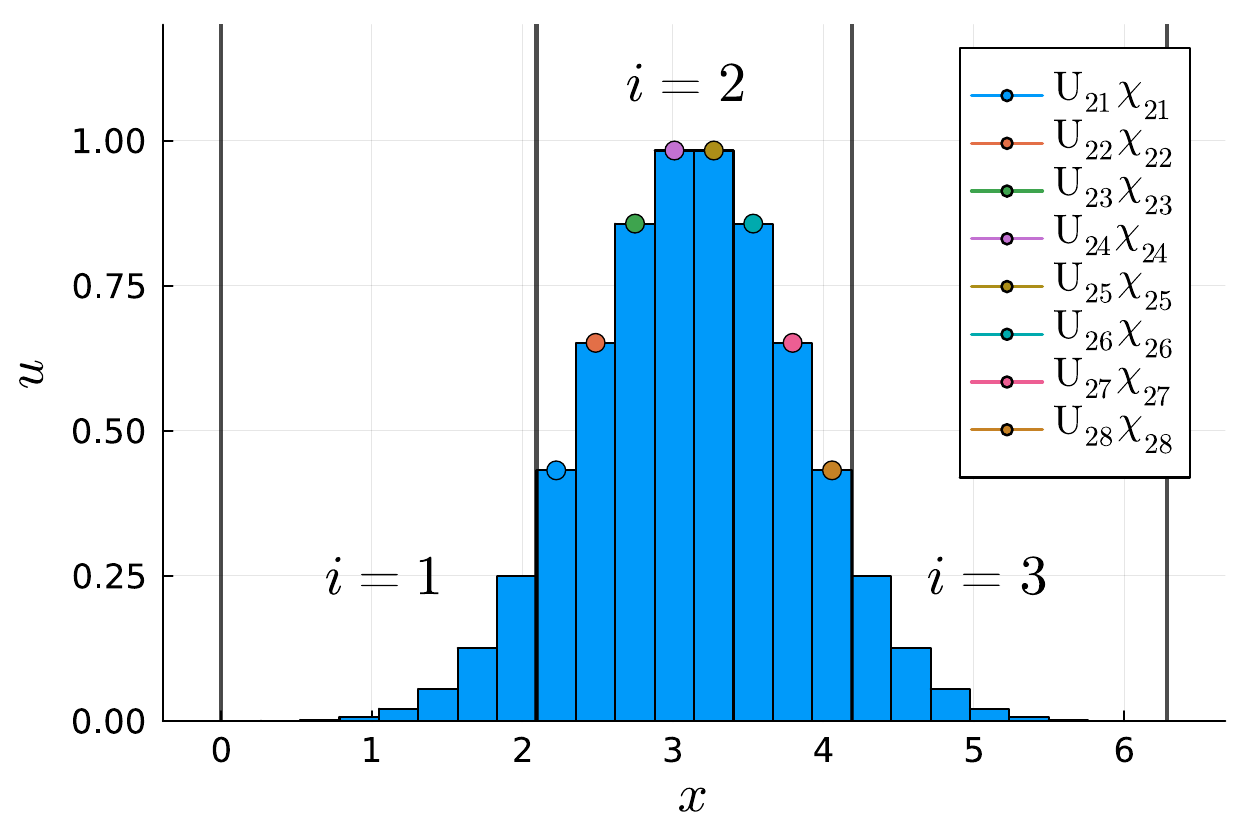}
    \caption{A Gaussian wave discretized by a finite difference scheme represented by a local basis of box functions for $I=3$ subdomains and $J = 8$ points per subdomain. The edges of the subdomains are indicated by the vertical grey lines.}
    \label{fig:local_basis}
\end{figure}

As stated earlier, we use a common finite difference basis $\{\chi_{ij}\}$ of size $J$ within each subdomain $\Omega_i$  to describe the solution:
\begin{equation}\label{eq:local_expansion}
        u_h(x,t) = \sum_{i=1}^{I} \sum_{j=1}^{J}\text{U}_{ij}(t)\chi_{ij}(x),
\end{equation}
where $\chi_{ij}$ is only nonzero within $\Omega_i$ and $\mathbf{U}(t) \in \mathbb{R}^{I \times J}$ contains the coefficient values. Note that in this case $\mathbf{U}(t)$ is simply a reshaped version of the solution vector $\mathbf{u}(t)$. 
\R{comment_4c}\revone{The finite difference basis functions are defined as}
\begin{equation}
    \chi_{ij}(x) = \begin{cases} 1 \quad &\text{if } \alpha_i + (j-\frac{1}{2})h \leq x < \alpha_i + (j+\frac{1}{2})h, \\ 
    0 \quad &\text{elsewhere},
    \end{cases}
\end{equation} 
similarly to \eqref{eq:FD_basis}.
 An example for a Gaussian solution profile for $I=3$ and $J = 8$ is given in Figure \ref{fig:local_basis}.

After obtaining $\mathbf{U}(t)$ at different points in time, the snapshot matrix is constructed as follows
\begin{equation}\label{eq:local_POD_snapshot}
    \mathbf{X}_\ell = \begin{bmatrix}
        \mathbf{U}^T(t_1) & \mathbf{U}^T(t_2) & \ldots & \mathbf{U}^T(t_s)
    \end{bmatrix} \quad \in \mathbb{R}^{J \times Is }.
\end{equation}
In this way, each subdomain is treated equally in the snapshot matrix.
This is what differentiates our work from \cite{farhat_local_POD} \R{comment_1b}\revone{and similar to what is done \cite{CHUNG2024DDROM}}. \R{comment_6}\revone{A schematic representation of this is shown in Figure \ref{fig:snapshot_matrix}, where the \gls{G-POD} snapshot matrix $\mathbf{X}$ is reshaped into the \gls{L-POD} snapshot matrix $\mathbf{X}_\ell$.}
\begin{figure}[ht]
    \centering
\tikzset{every picture/.style={line width=0.75pt}}        

\begin{tikzpicture}[x=1pt,y=1pt,scale=0.8, every node/.style={scale=0.8}]

\draw[fill={rgb,255:red,208;green,2;blue,27}](60,30) rectangle (130,70);
\draw[fill={rgb,255:red,245;green,166;blue,35}] (60,70) rectangle (130,110);
\draw[fill={rgb,255:red,74;green,144;blue,226}] (60,110) rectangle (130,150);

\node at (40,50) {$\Omega_3$};
\node at (40,90) {$\Omega_2$};
\node at (40,130) {$\Omega_1$};

\node at (95,165) {$\mathbf{X} \in \mathbb{R}^{N \times s}$};

\draw[->] (140,90) -- (200,90);
\node at (170,100) {Reshape};

\draw[fill={rgb,255:red,74;green,144;blue,226}] (220,70) rectangle (290,110);
\draw[fill={rgb,255:red,245;green,166;blue,35}] (290,70) rectangle (360,110);
\draw[fill={rgb,255:red,208;green,2;blue,27}] (360,70) rectangle (430,110);

\node at (255,55) {$\Omega_1$};
\node at (325,55) {$\Omega_2$};
\node at (395,55) {$\Omega_3$};

\node at (325,125) {$\mathbf{X}_{\ell} \in \mathbb{R}^{J \times Is}$};

\end{tikzpicture}
    \caption{Schematic representation of the \gls{G-POD} snapshot matrix $\mathbf{X}$ being reshaped into the \gls{L-POD} snapshot matrix $\mathbf{X}_\ell$ for $I = 3$.}
    \label{fig:snapshot_matrix}
\end{figure}
Using a single common snapshot matrix allows us to obtain a space-local \gls{POD} basis that generalizes well outside the training data.
The way we divide the domain into subdomains is decided by what results in an \gls{L-POD} basis that generalizes best on a validation data set, see section \ref{sec:proj_error}.

Note that this matrix has fewer rows but more columns than the \gls{G-POD} snapshot matrix, see \eqref{eq:global_POD_snapshot}. This makes the \gls{SVD} cheaper to compute for a large number of snapshots $s$ \cite{SVD_complexity}.
As in the global case, we use a \gls{SVD} of the snapshot matrix to obtain a truncated basis $\hat{\bGamma} \in  \mathbb{R}^{J \times q}$ with $q < J$ from the left-singular vector. In terms of the local basis $\{ \chi_{ij}\}$ the \gls{POD} basis is written as
\begin{equation}
    \gamma_{ij}(x) = \sum^J_{k=1} \chi_{ik}(x)\hat{\Gamma}_{kj},
\end{equation}
analogous to \eqref{eq:POD_basis}.
In this basis, the approximated solution $u_\ell$ is obtained as 
\begin{equation}\label{eq:local_approximation}
    u_h(x,t) \approx u_\ell(x,t) = \sum^{I}_{i=1}  \sum^{q}_{j=1}\text{A}_{ij}(t) \gamma_{ij}(x),
\end{equation}
similarly to the global case, see \eqref{eq:approximation}. 
The coefficients $\text{A}_{ij}$ are obtained as
\begin{equation}\label{eq:a}
    \mathbf{a}_\ell = \bGamma^T \mathbf{u},
\end{equation}
where
\begin{equation}\label{eq:gamma}
    \bGamma = \begin{bmatrix}
        \hat{\bGamma} & & \\
        & \ddots & \\
        & & \hat{\bGamma} 
    \end{bmatrix}\quad \in \mathbb{R}^{N \times r}.
\end{equation}
This matrix, containing non-overlapping blocks, is constructed such that multiplying by its transpose projects $\mathbf{u}$ onto the \gls{L-POD} basis. 
For the mapping from $\mathbf{a}_\ell \in \mathbb{R}^{Iq}$ to $\mathbf{A}\in \mathbb{R}^{I \times q}$ we follow the same convention as for $\mathbf{u}$ and $\mathbf{U}$. The effective number of basis functions for \gls{L-POD} is $Iq = r$ with $q < J$. This basis is orthonormal, i.e. $\bGamma^T \bGamma = \mathbf{I}$. The sparsity of $\bGamma$, as opposed to $\bPhi$, is what yields a \gls{ROM} which is cheaper to evaluate than a \gls{G-POD}-based \gls{ROM} \cite{farhat_local_POD}.

\begin{figure}[ht]
    \centering
    \includegraphics[width = 0.6\textwidth]{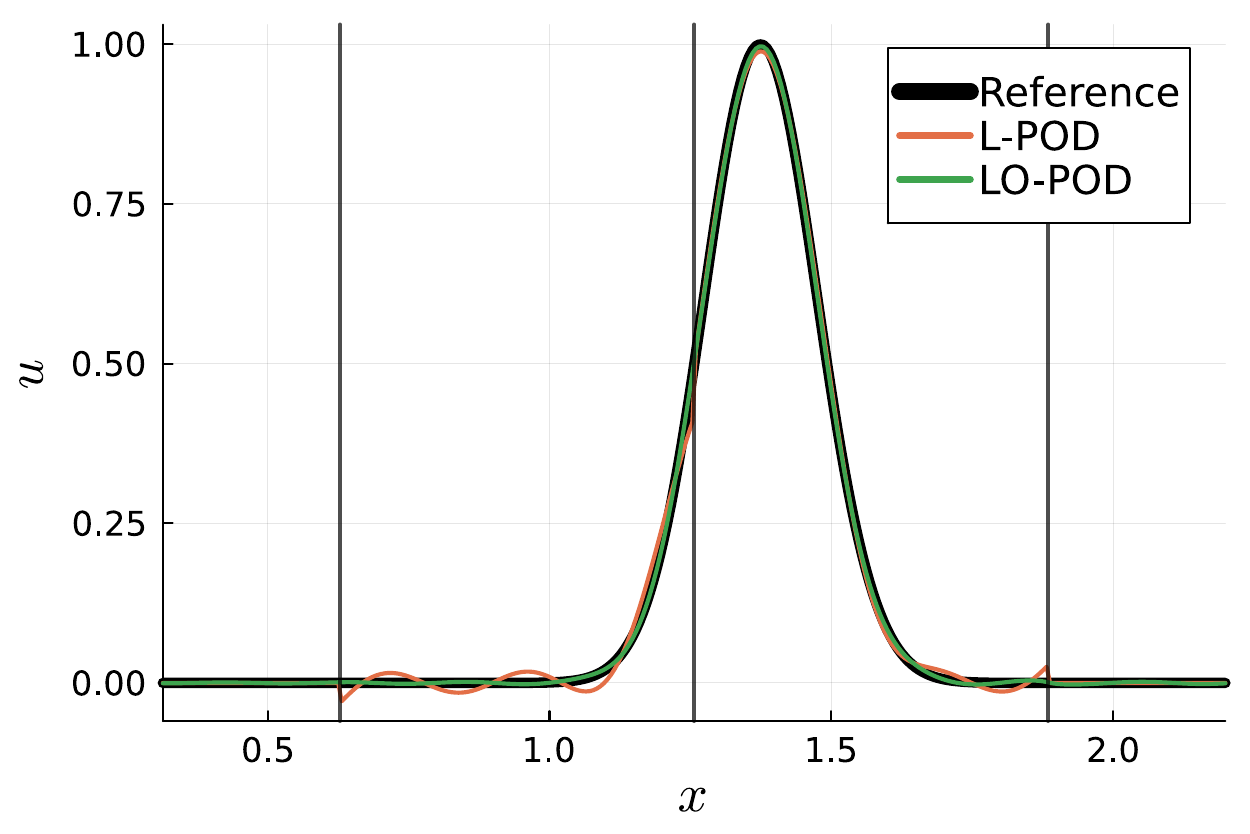}
    \caption{\gls{L-POD} and \acrshort{LO-POD} representation of a Gaussian wave. The data used to obtain the local \gls{POD} basis is explained in Section \ref{sec:results}. The depicted snapshot is part of the snapshot matrix used to obtain the \gls{POD} basis. The edges of the subdomains are indicated by the vertical grey lines. Only part of the domain $\Omega = [0,2\pi)$ is shown.}
    \label{fig:discontinuity}
\end{figure}

An example of an \gls{L-POD} approximation $u_\ell$ to a Gaussian wave is shown in Figure \ref{fig:discontinuity}. For the \gls{L-POD} we used $I=10$ subdomains with $q=6$ basis functions per subdomain. The training data used to obtain the basis are discussed in the results section (Section \ref{sec:results}), as shown in Figure \ref{fig:reference}.
Although the \gls{L-POD} approximation is accurate in most of the domain, some oscillations appear near the boundaries. In addition, the approximation is not guaranteed to be smooth on the edges of the subdomains, see Figure \ref{fig:discontinuity}. Discontinuities appear when transitioning from one subdomain to the next. In Section \ref{sec:results}, we will show that these discontinuities tend to grow increasingly severe when using the \gls{L-POD} basis to construct a \gls{ROM}, see Figure \ref{fig:trajectories}. To remedy this issue, we introduce space-local \gls{POD} with overlapping subdomains in the next section. 

\subsection{Local POD with overlapping subdomains}

In finite element discretizations, similar discontinuities are resolved by imposing continuity of the approximation, while using a local basis. In this work, we aim to achieve the same by introducing a novel space-local \gls{POD} formulation with \textit{overlapping} subdomains. We will refer to this approach as \acrshort{LO-POD}. To start off, we subdivide the domain $\Omega$ into $I$ overlapping subdomains $\Omega_i = [\alpha_i,\beta_i)$. This means $\beta_{i-2} = \alpha_i < \beta_{i-1} = \alpha_{i+1} < \beta_i = \alpha_{i+2}$.
Once again, each subdomain contains $J$ grid points. Note that due to the overlap, each grid point is now located in two subdomains. This means $\frac{I \cdot J }{2} = N$, as opposed to $I \cdot J = N$ for \gls{L-POD}.

\begin{figure}[ht]
    \centering
    \includegraphics[width = 0.6\textwidth]{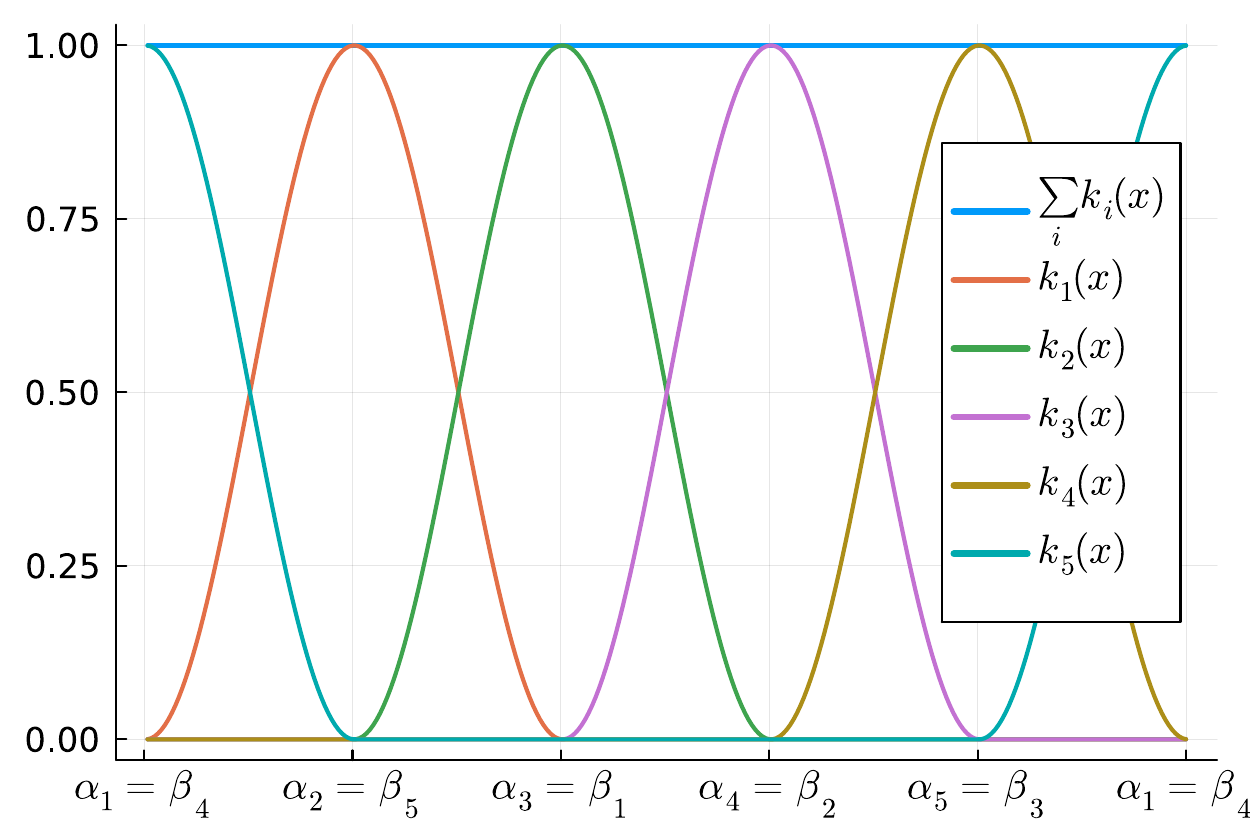}
    \caption{Kernels $k_i$ for a subdivision of the periodic domain into five overlapping subdomains of equal size.}
    \label{fig:partition_unity}
\end{figure}

To obtain a smooth approximation $u_\ell$, we require the \acrshort{LO-POD} basis to smoothly decay to zero on the edge of the subdomain. This is enforced through a post-processing step of the local snapshot matrix given by
\eqref{eq:local_POD_snapshot}.
For this purpose, we introduce a kernel $k_i(x)$ to divide the solution between the subdomains. This places the following constraint on the kernels:  
\begin{equation}\label{eq:kernel_constraint}
    \sum^I_{i=1} k_i(x) = 1,
\end{equation}
such that this set of functions forms a partition of unity. 
Here, we propose the following kernel
\begin{equation}
    k_i(x) = \begin{cases}
        \sin^2(\frac{1}{2}\frac{x-\alpha_i}{\alpha_i - \alpha_{i+1}}\pi)  \quad &\text{if }  \alpha_i \leq x < \alpha_{i+1},\\
        \sin^2(\frac{1}{2}\frac{x-\alpha_{i+1}}{\alpha_{i+2} - \alpha_{i+1}}\pi + \frac{1}{2}\pi) \quad &\text{if } \alpha_{i+1} \leq x < \alpha_{i+2}, \\
        0 \quad &\text{elsewhere},
    \end{cases}
\end{equation}
which is chosen as it smoothly decays to zero at the subdomain boundaries. This approach is inspired by the partition of unity method \cite{PUM}. Here, the right half of each subdomain is overlapped by the subdomain to its right, and the left half by the subdomain to its left. A visualization of the kernels is displayed in Figure \ref{fig:partition_unity}. 
Different kernels and different amounts of overlap between subdomains can also be considered. However, we consider this outside the scope of this paper. 

Using the introduced kernel, the coefficients $\mathbf{U}$ for the local expansion in \eqref{eq:local_expansion} are obtained as
\begin{equation}\label{eq:post_processing}
    \text{U}_{ij}(t) = \frac{(\chi_{ij}(x), k_i(x)u_h(x,t))}{(\chi_{ij}(x),\chi_{ij}(x))}.
\end{equation}
These integrals can simply be approximated using the midpoint rule for integration \cite{midpoint_dragomir1998some}.
The remaining procedure stays the same as for \gls{L-POD}, i.e., we build the snapshot matrix in \eqref{eq:local_POD_snapshot} and carry out a \gls{SVD} to obtain $\hat{\bGamma}$. However, what changed is that the blocks in $\bGamma$, see \eqref{eq:gamma}, are now overlapping. This means the basis is no longer orthonormal, i.e. $\bGamma^T\bGamma \neq \mathbf{I}$. Due to the non-orthonormality, the coefficients $\mathbf{a}_\ell$ for the expansion in \eqref{eq:local_approximation} are obtained by solving the following linear system 
\begin{equation}
(\bGamma^T\bGamma)\mathbf{a}_\ell = \bGamma^T\mathbf{u},
\end{equation}
as opposed to \eqref{eq:a}.
Obtaining the coefficients from the \gls{FOM} solution is therefore more expensive. 
In addition, evaluating the resulting \gls{ROM} requires solving a linear system, as will be discussed in Section \ref{sec:ROM}. This makes the \acrshort{LO-POD} \gls{ROM} more expensive to evaluate than its non-overlapping counterpart \gls{L-POD}, see Sections \ref{sec:Galerkin_projection} and \ref{sec:convergence}.
However, looking at Figure \ref{fig:discontinuity}, we find that using the \acrshort{LO-POD} results in a much smoother approximation of the Gaussian wave. Note that the parameters are kept the same, i.e., $I=10$ and $q=6$. In Section \ref{sec:proj_error}, the associated error will be quantified more precisely.

\section{Space-local, energy-conserving reduced-order model}\label{sec:ROM}

\subsection{Projection on reduced basis}\label{sec:projection}

Consider again $\mathbf{c}$, the state vector of the full-order model, and $\bPhi \in \mathbb{R}^{N \times r}$ a reduced basis where $r<N$, obtained from either \gls{G-POD}, \gls{L-POD}, or \acrshort{LO-POD}. For the space-local approaches, this operator was referred to as $\bGamma$, see Section \ref{sec:L-POD}. The corresponding `Gram' matrix $\mathbf{S} \in \mathbb{R}^{r\times r}$ of this basis is computed as 
\begin{equation}
    \mathbf{S} = \bPhi^T \bPhi.
\end{equation}
We can project $\mathbf{c}$ onto the subspace spanned by the \gls{POD} basis as
\begin{equation}\label{eq:projector}
    \mathbf{c}_r = \underbrace{\bPhi \mathbf{S}^{-1}\bPhi^T}_{=:\mathbf{P}} \mathbf{c}.
\end{equation}
Note that for both \gls{G-POD} and \gls{L-POD}, computing the inverse of $\mathbf{S}$ is trivial, as it is simply the identity. This is because the basis is orthonormal. However, for \acrshort{LO-POD} the basis is non-orthogonal, which makes computing its inverse less trivial. As stated earlier, this is a downside to \acrshort{LO-POD} and increases its computational cost. 
It is quite straightforward to see that $\mathbf{P}$ is idempotent as $\mathbf{P}^2 = \mathbf{P}$.
The \gls{POD} coefficient vector is obtained as
\begin{equation}\label{eq:POD_coefficient_vector}
    \mathbf{a} = \mathbf{S}^{-1}\bPhi^T\mathbf{c}.
\end{equation}

\subsection{Galerkin projection}\label{sec:Galerkin_projection}

Having obtained a reduced basis, we construct a \gls{ROM} for the advection equation as follows \cite{POD1,POD2,brunton2019data,projection_ROMs_survey}: based on \eqref{eq:semi-discrete} we define the residual $\mathbf{r} \in \mathbb{R}^N$ as 
\begin{equation}
    \mathbf{r}(\mathbf{u}_r) := \frac{\text{d}\mathbf{u}_r}{\text{d}t}+c\mathbf{D}\mathbf{u}_r \neq \mathbf{0},
\end{equation}
where we replaced $\mathbf{u}$ by the \gls{POD} approximation $\mathbf{u}_r = \bPhi \mathbf{a}$. As $\mathbf{u}$ comes from the finite difference discretization of the advection equation, we have $\mathbf{c}(t) = \mathbf{u}(t)$.  
Next, we carry out the Galerkin projection of the residual on the \gls{POD} basis, according to \eqref{eq:POD_coefficient_vector}, and set this to zero:
\begin{equation}
    \mathbf{S}^{-1}\bPhi^T\mathbf{r}(\mathbf{u}_r) = \mathbf{0}.
\end{equation}
This ensures the residual is orthogonal to the basis. 
This results in the following \gls{ROM}:
\begin{equation}\label{eq:ROM}
     \frac{\text{d}\mathbf{a}^\text{ROM}}{\text{d}t} := -c\mathbf{S}^{-1}\mathbf{A}\mathbf{a}^\text{ROM},
\end{equation}
where the \gls{ROM} operator $\mathbf{A} \in \mathbb{R}^{r\times r}$ is defined as 
\begin{equation}
   \mathbf{A} := \bPhi^T \mathbf{D}\bPhi. 
\end{equation}
Note that we introduced $\mathbf{a}^\text{ROM}$ here as the \gls{POD} state vector predicted by \gls{ROM}. This is typically not equal to the true $\mathbf{a}$, see \eqref{eq:POD_coefficient_vector}, past $t=0$.
Going from the \gls{FOM} to the \gls{ROM}, we reduced the \gls{DOF} in the system from $N$ to $r$. In the \gls{G-POD} case $\mathbf{A}$ is typically dense, whereas in the \gls{L-POD}/\acrshort{LO-POD} it is sparse. The sparsity decreases the cost of evaluating the \gls{ROM} \cite{farhat_local_POD}. The amount of nonzero entries in the \gls{G-POD} \gls{ROM} operator scales with $\mathcal{O}(r^2)$. For \gls{L-POD} it scales with $\mathcal{O}(((n+1)q)^2I)$, where $n$ is the number of neighbors per subdomain. This accounts for the interaction between the subdomains. For \acrshort{LO-POD}, the neighbors of the neighbors have to be included in the interactions due to the overlap. This results in a scaling of $\mathcal{O}(((\tilde{n}+1)q)^2I)$, where $\tilde{n}$ is the number of subdomains that can be reached by crossing at most one subdomain starting from one of the subdomains. For 1D problems $n =2$ and $\tilde{n} = 4$. For a large number of subdomains $I$ and a small number of \gls{POD} basis functions per subdomain $q$, with the total number of \gls{POD} modes being $I \cdot q = r$, this scales more favorably than \gls{G-POD}. This allows us to include a larger number of basis functions in the \gls{POD} basis without sacrificing computational efficiency. Hyper-reduction techniques, including energy-conserving ones, can also be employed to further sparsen the \gls{ROM} operators \cite{hyperreduction_1,hyperreduction_2,KLEIN2024112697}.

\subsection{Energy conservation of the ROM}

To ensure stability of the \glspl{ROM}, we aim to mimic the energy conservation property of the \gls{FOM}, see \eqref{eq:discr_E_conservation}. To investigate if the constructed \glspl{ROM} satisfy this property we define the \gls{ROM} energy $E_r^\text{ROM}$ as
\begin{equation}
    E_r^\text{ROM} = \frac{h}{2}(\mathbf{u}_r^\text{ROM})^T(\mathbf{u}_r^\text{ROM}),
\end{equation}
where $\mathbf{u}_r^\text{ROM} := \bPhi \mathbf{a}^\text{ROM}$.
Employing \eqref{eq:ROM} we obtain the \gls{ROM} evolution of $\mathbf{u}_r^\text{ROM}$ as 
\begin{equation}
    \frac{\text{d} \mathbf{u}_r^\text{ROM}}{\text{d}t} = \bPhi\frac{\text{d}\mathbf{a}^\text{ROM}}{\text{d}t} =-c\bPhi\mathbf{S}^{-1}\mathbf{A}\mathbf{a}^\text{ROM}.
\end{equation}
The evolution of the \gls{ROM} energy follows as 
\begin{equation}\label{eq:change_in_energy}
\begin{split}
    \frac{\text{d}E_r^\text{ROM}}{\text{d}t} &= h(\mathbf{u}_r^\text{ROM})^T \frac{\text{d} \mathbf{u}_r^\text{ROM}}{\text{d}t}  = -ch(\mathbf{u}_r^\text{ROM})^T\bPhi\mathbf{S}^{-1}\mathbf{A}\mathbf{a}^\text{ROM} \\ &= -ch(\mathbf{a}^\text{ROM})^T \underbrace{\bPhi^{T} \bPhi}_{=\mathbf{S}}\mathbf{S}^{-1}\bPhi^T \mathbf{D} \bPhi \mathbf{a}^\text{ROM} \\ &= -ch(\mathbf{a}^\text{ROM})^T \bPhi^T \mathbf{D} \bPhi \mathbf{a}^\text{ROM} = -ch(\mathbf{u}_r^\text{ROM})^T \mathbf{D}\mathbf{u}_r^\text{ROM}= 0,
\end{split}
\end{equation}
employing the product rule of differentiation. From the final expression, we conclude that a \gls{ROM} based on Galerkin projection inherits energy conservation and stability from the \gls{FOM}. This is true for both orthogonal and non-orthogonal projections. This means that not only the \gls{G-POD} \gls{ROM} inherits this property, as shown in \cite{podbenjamin}, but also the newly introduced space-local approaches \gls{L-POD} and \acrshort{LO-POD}.

\section{Results \& discussion}\label{sec:results}

\subsection{Test case setup}

\begin{figure}[ht]
    \centering
    \includegraphics[width = 0.48\textwidth]{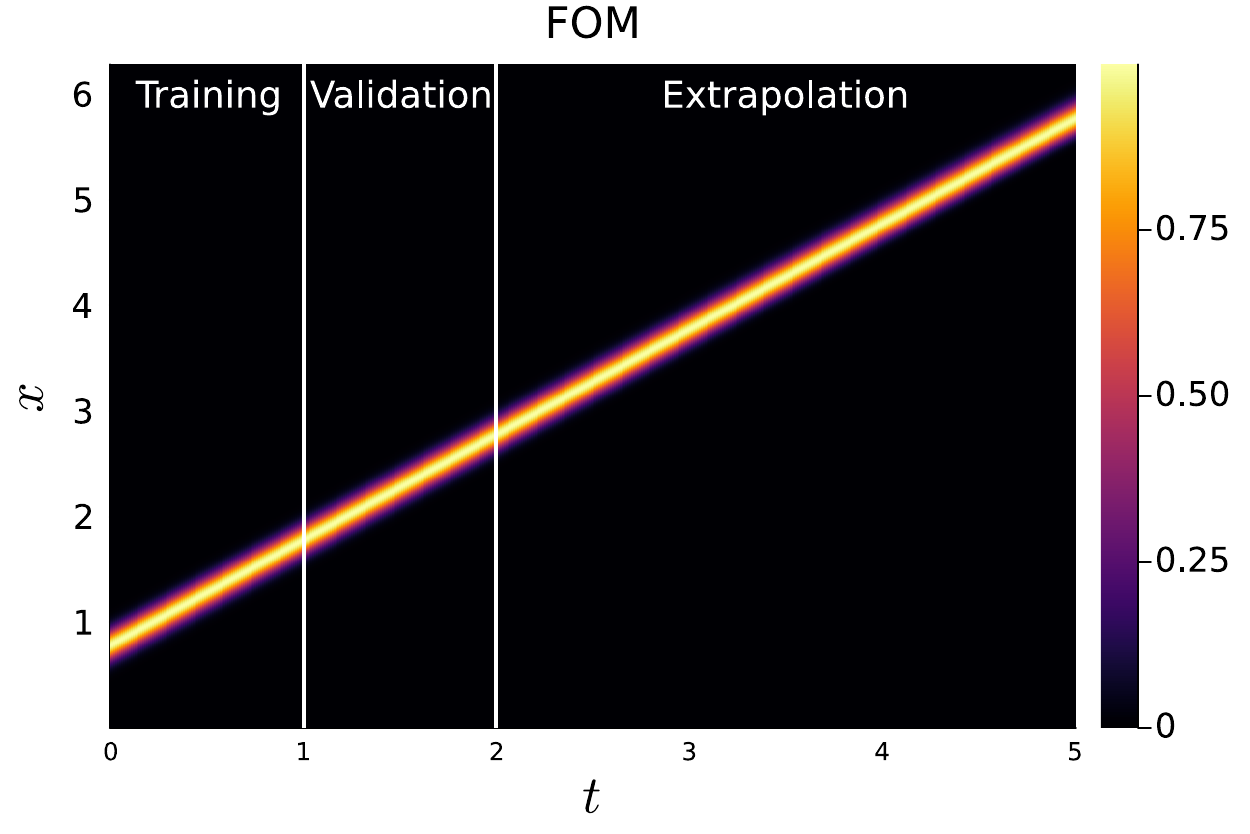}
    \caption{Reference simulation of a Gaussian wave being advected throughout the domain. White bars separate the training data from the validation data and the validation data from the extrapolation data, respectively.}
    \label{fig:reference}
\end{figure}

To construct a \gls{ROM}, we first require data from the \gls{FOM}. As stated earlier, we use the finite difference discretization of the advection equation, detailed in Section \ref{sec:discretization}, for this purpose. The system is simulated on a periodic domain $\Omega = [0,2 \pi)$ discretized with $N=1000$ grid points for $t = [0,5]$ and constant $c=1$. A time step size of $\Delta t = 0.01$ is used for the time integration, using a \gls{RK4} scheme. This is the largest time step size that still yielded stable simulations. Snapshot data for the \gls{ROM} construction is collected in the interval $t = [0,1]$. We refer to this as the training data. Data generated in the interval $t=(1,2]$ will be referred to as validation data. This data is used to evaluate the generalization of the \gls{POD} basis. The remaining part of the simulation data, $t=(2,5]$, will be used to evaluate the extrapolation capabilities of the \glspl{ROM}. As an initial condition, we use a Gaussian wave centered around $x = \frac{1}{4}\pi$, namely 
\begin{equation}\label{eq:init_cond}
    u(x,0) = \exp(-50(x - \frac{1}{4}\pi)^2).
\end{equation}
A visualization of this simulation is displayed in Figure \ref{fig:reference}.

For the \glspl{ROM} we consider the three bases discussed in this work: the global \gls{POD} basis \gls{G-POD}, the space-local \gls{POD} basis \gls{L-POD}, and the space-local \gls{POD} basis with overlapping subdomains \acrshort{LO-POD}. For the space-local approaches, the domain is subdivided into $I$ subdomains. For the local basis $\{\chi _{ij}\}$, we use $J$ box functions contained within each subdomain, such that $IJ = N$. In this way, the local basis aligns with the \gls{FOM} finite difference basis, see \eqref{eq:FD_basis}.  The integrals in the \acrshort{LO-POD} post-processing step, equation \eqref{eq:post_processing}, are approximated using the midpoint rule for integration \cite{midpoint_dragomir1998some}. For the \gls{ROM} time integration, we use the same \gls{RK4} scheme as for the \gls{FOM}.

To evaluate the \glspl{ROM} we evaluate the difference between the \gls{ROM} solution $\mathbf{u}^\text{ROM}_r := \bPhi \mathbf{a}^\text{ROM}$ and the \gls{FOM} solution $\mathbf{u}^\text{FOM}$: 
\begin{equation}\label{eq:error}
\underbrace{\frac{\mathbf{u}^\text{ROM}_r(t) -\mathbf{u}^\text{FOM}(t)}{||\mathbf{u}^\text{FOM}(t)||_2}}_{:=\text{solution error}} = \underbrace{\frac{\mathbf{u}_r^\text{ROM}(t) -\mathbf{P}\mathbf{u}^\text{FOM}(t)}{||\mathbf{u}^\text{FOM}(t)||_2}}_{:=\text{ROM error}} + \underbrace{\frac{\mathbf{P}\mathbf{u}^\text{FOM}(t)-\mathbf{u}^\text{FOM}(t)}{||\mathbf{u}^\text{FOM}(t)||_2}}_{:=\text{projection error}},
\end{equation}
where $\mathbf{P}$ projects the solution on the \gls{POD} basis, see \eqref{eq:projector}.
This difference will be referred to as the solution error. In \eqref{eq:error}, this error is decomposed as the sum of the \gls{ROM} error (the error made by the \gls{ROM} during the time integration) and the projection error (the error made by the reduced basis approximation). 
Our implementation, in Julia \cite{Julia-2017}, of the introduced methodologies and experiments is freely available on Github, see \url{https://github.com/tobyvg/local_POD_overlap.jl}. 

\subsection{Projection error}\label{sec:proj_error}

\begin{figure}[ht]
    \centering
    \includegraphics[width = 0.48\textwidth]{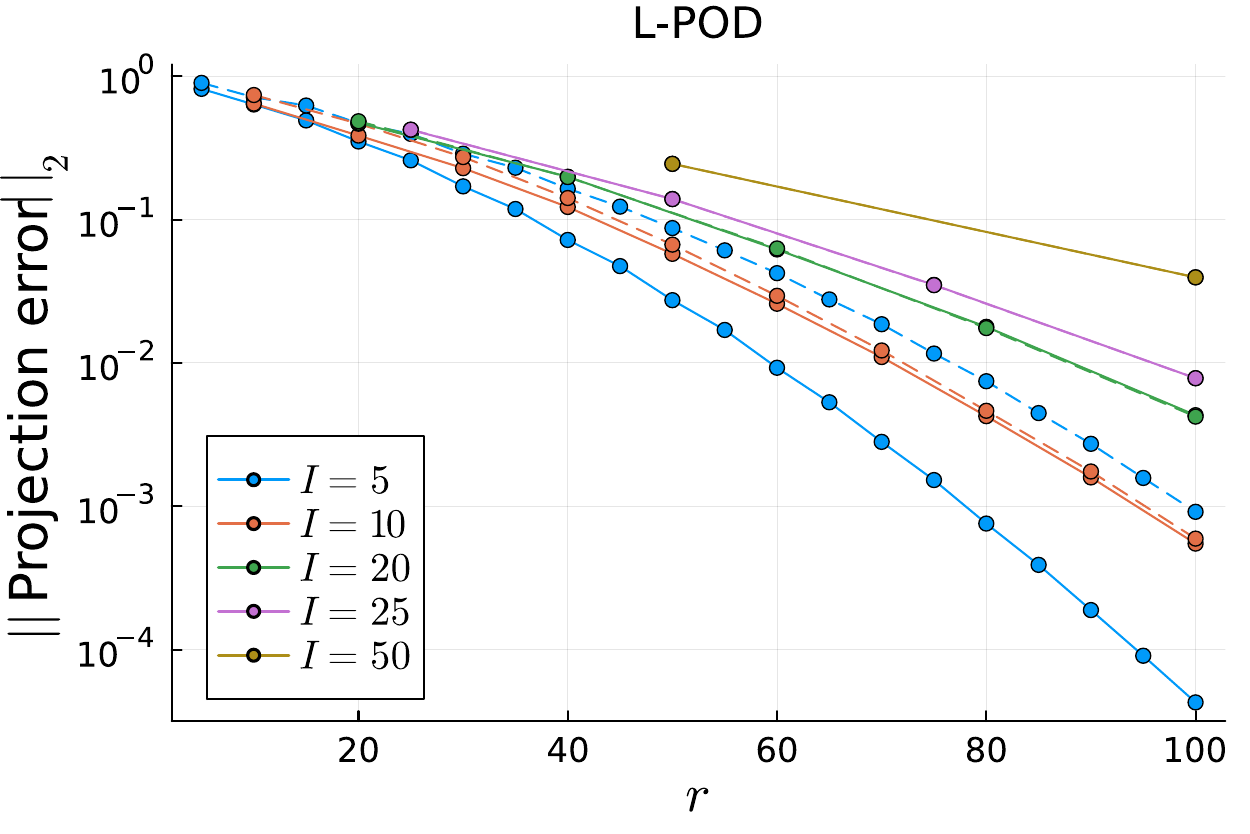}
    \includegraphics[width = 0.48\textwidth]{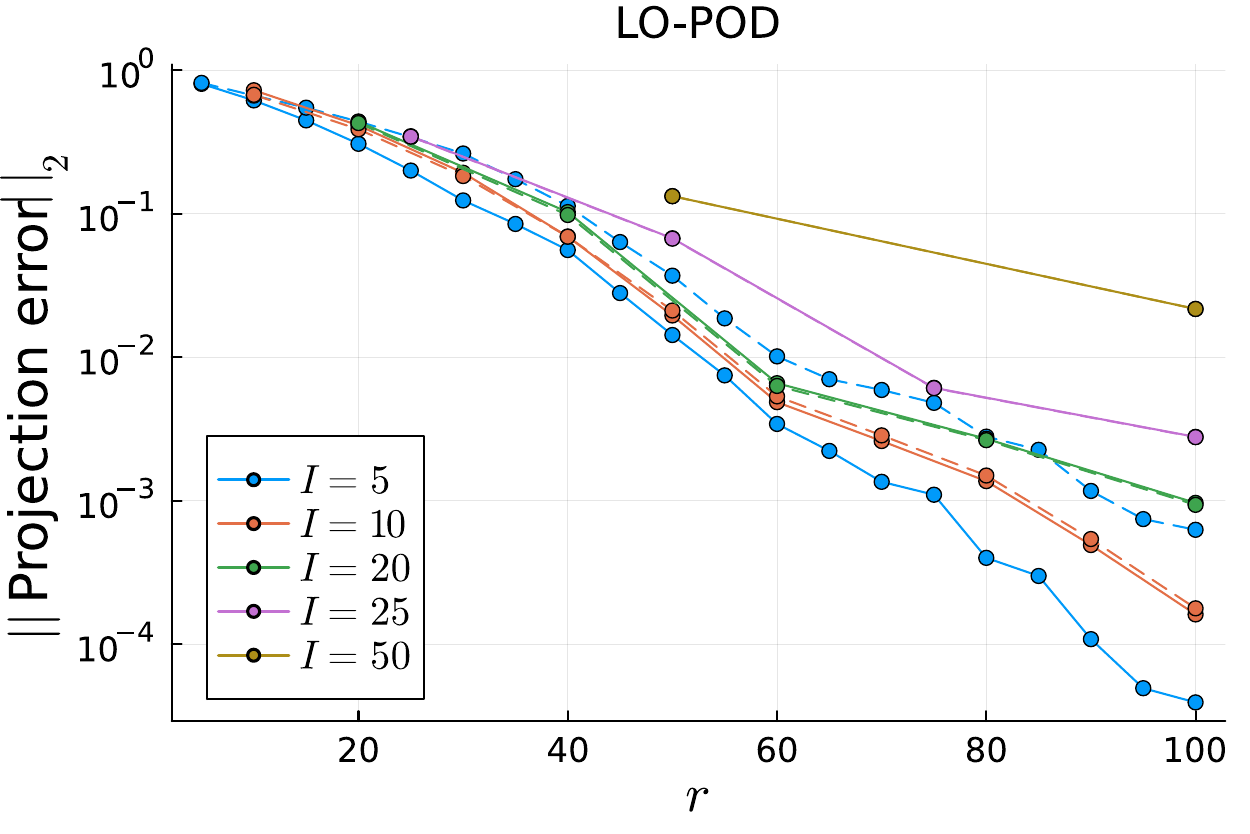}
    \includegraphics[width = 0.48\textwidth]{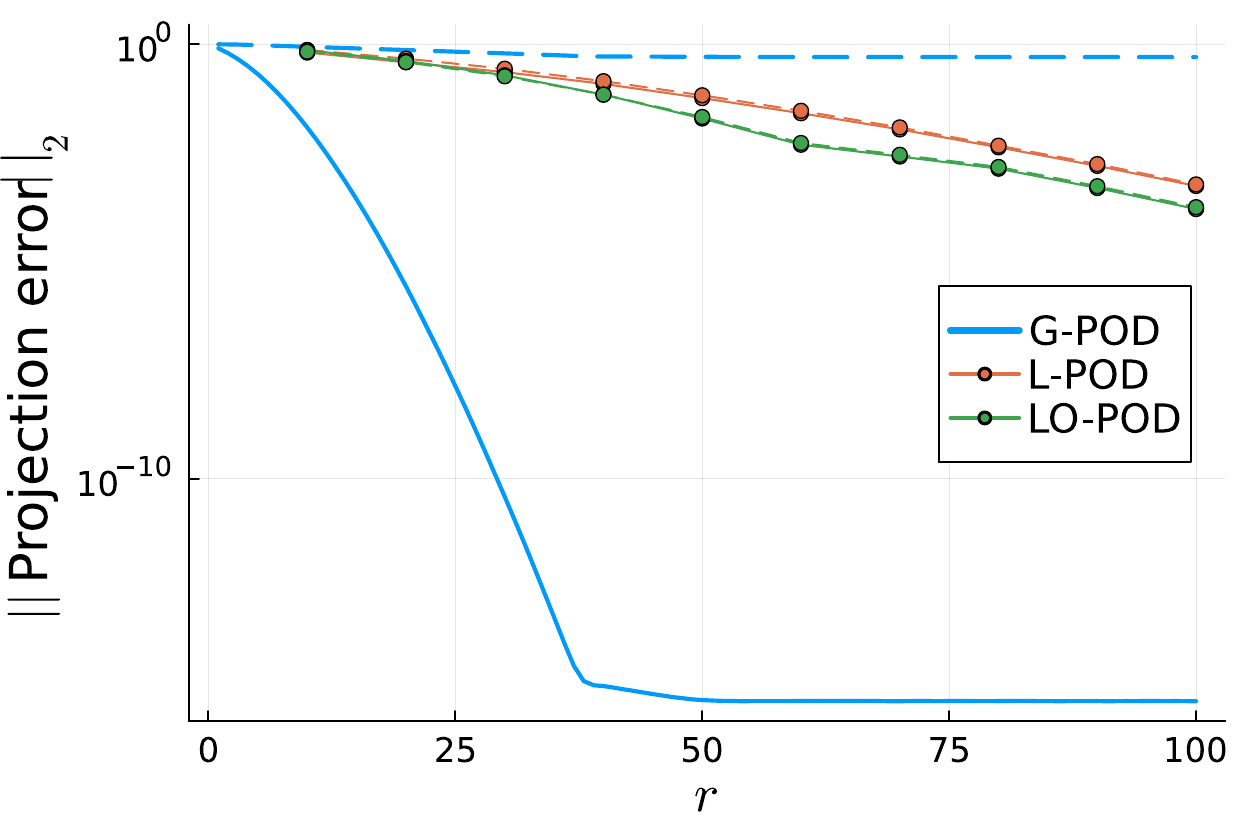}
    \caption{(Top) Projection error evaluated over the training (solid line) and validation (dashed line) data for different $I$ for the space-local approaches. (Bottom) Projection error for $I=10$ for \gls{L-POD} and \acrshort{LO-POD} evaluated over the training (solid line) and validation (dashed line) data. The projection error for \gls{G-POD} is also depicted.}
    \label{fig:optimization}
\end{figure}

For the construction of the \gls{L-POD} and \acrshort{LO-POD} \glspl{ROM}, we have to determine the number of subdomains $I$ and the number of modes per subdomain $q$, with $Iq=r$, which results in a basis that generalizes best outside the training data. To do this, we evaluate the projection error, see \eqref{eq:error},
on both the training and validation data. The results are depicted in Figure \ref{fig:optimization} for different $I$, averaged over both the training (solid line) and validation data (dashed line).
A basis that generalizes well should perform well on both the training and validation data. 
We observe that for $I=5$, the training error is lowest, but the validation error is rather large, indicating that the basis does not generalize well. For higher values of $I$, the training and validation errors are much closer, but tend to increase with increasing $I$. In this case, $I=10$ is considered optimal, as the training and validation errors are small and of similar size. This is true for both \gls{L-POD} and \acrshort{LO-POD}. For the remainder of this text, we will therefore stick to $I=10$ for the construction of the local \glspl{ROM}. In general, the choice of $I$ is likely problem-dependent, and an a priori study similar to Figure \ref{fig:optimization} needs to be performed to determine an optimal value.

Next, we compare the performance of the space-local approaches against each other, as well as to \gls{G-POD}. This is also displayed in Figure \ref{fig:optimization} (bottom plot). We observe that the projection error converges faster for \acrshort{LO-POD}, as compared to \gls{L-POD}. This can be explained by the fact that the \acrshort{LO-POD} basis functions are constructed from twice as many finite difference points as for \gls{L-POD}, due to the overlapping subdomains. This means the fit is carried out over twice as many \gls{DOF}, which yields a lower projection error. 
The difference between convergence on the training and validation data is slight for both space-local approaches. On the other hand, for \gls{G-POD}, the convergence of the projection error on the training data is very fast, see \eqref{eq:minimal_projection_error}. However, on the validation data, the error hardly converges. This is because the \gls{G-POD} is only suited for representing the wave on the left side of the domain, as detailed in the next section. We conclude that the space-local approaches yield a basis that generalizes better than the global approach for this particular problem. 

\subsection{Resulting basis}

\begin{figure}[ht]
    \centering
    \includegraphics[width = 0.48\textwidth]{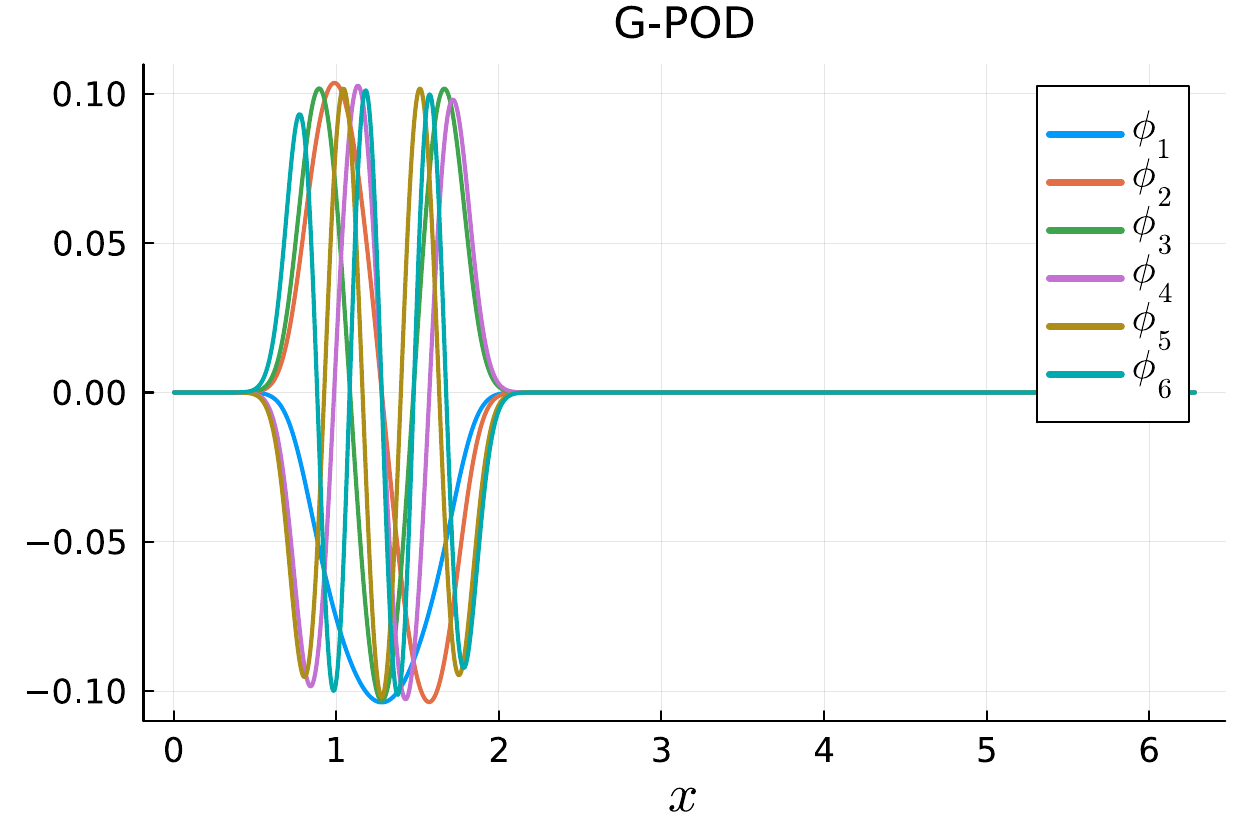}
    \includegraphics[width = 0.48\textwidth]{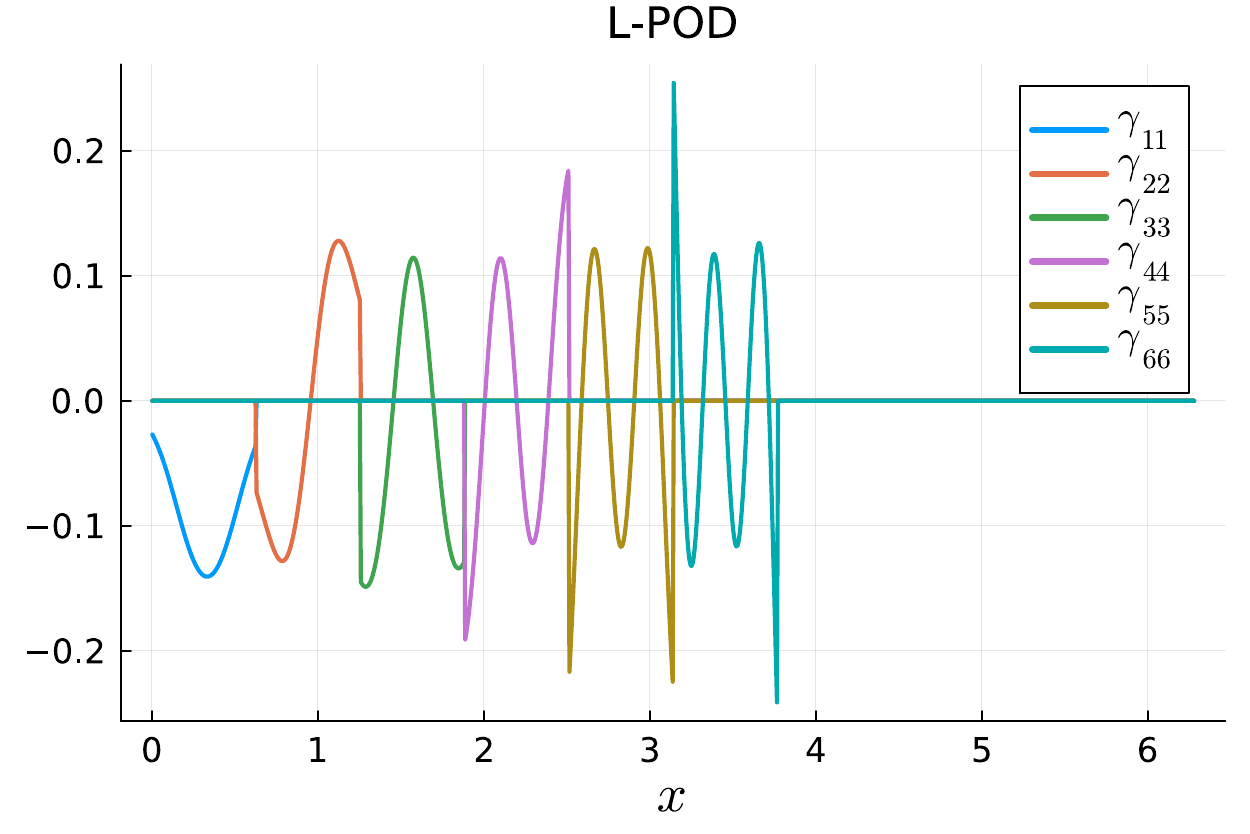}
    \includegraphics[width = 0.48\textwidth]{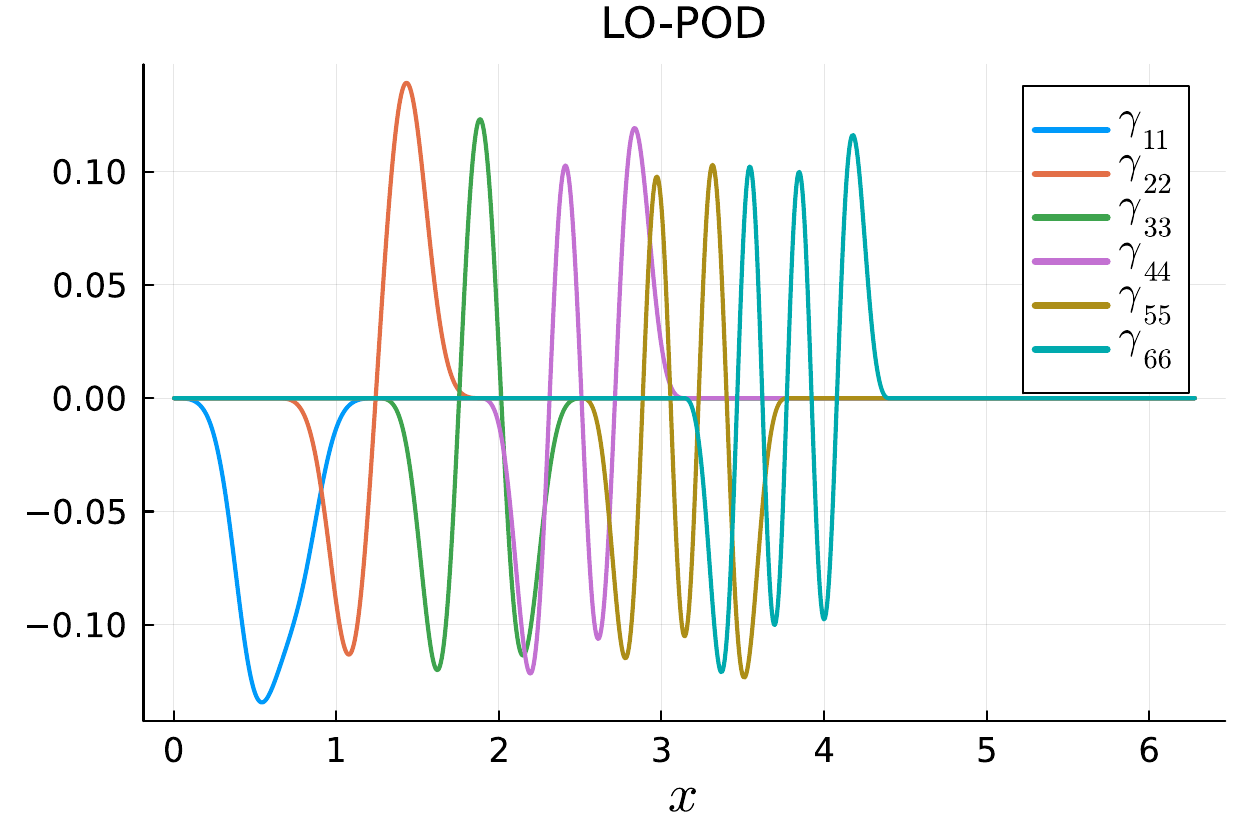}    
    \caption{The first six \gls{G-POD} modes and the first six space-local \gls{POD} modes for \gls{L-POD} and \acrshort{LO-POD}. For visualization purposes, the space-local \gls{POD} modes are displayed on adjacent subdomains.}
    \label{fig:basis}
\end{figure}

The resulting local basis functions for $q=6$ are displayed in Figure \ref{fig:basis} along with the first six \gls{G-POD} modes.
We find that the \gls{G-POD} modes are only nonzero where the traveling wave is represented in the training data. This explains why the error on the validation data hardly converged in Figure \ref{fig:optimization}. 
For the space-local approaches, where a common basis is ``copied/pasted" across the entire domain, this is not an issue. Regarding \gls{L-POD}, we observe that the obtained basis functions do not smoothly decay to zero at the edge of the subdomain, but instead end abruptly. Finally, we observe that for \acrshort{LO-POD} the post-processing procedure in \eqref{eq:post_processing} indeed results in a basis that smoothly decays on the edge of the subdomains.

\subsection{Accuracy and energy conservation of ROM}\label{sec:error_and_structure}

\begin{figure}[ht]
    \centering
    \includegraphics[width = 0.48\textwidth]{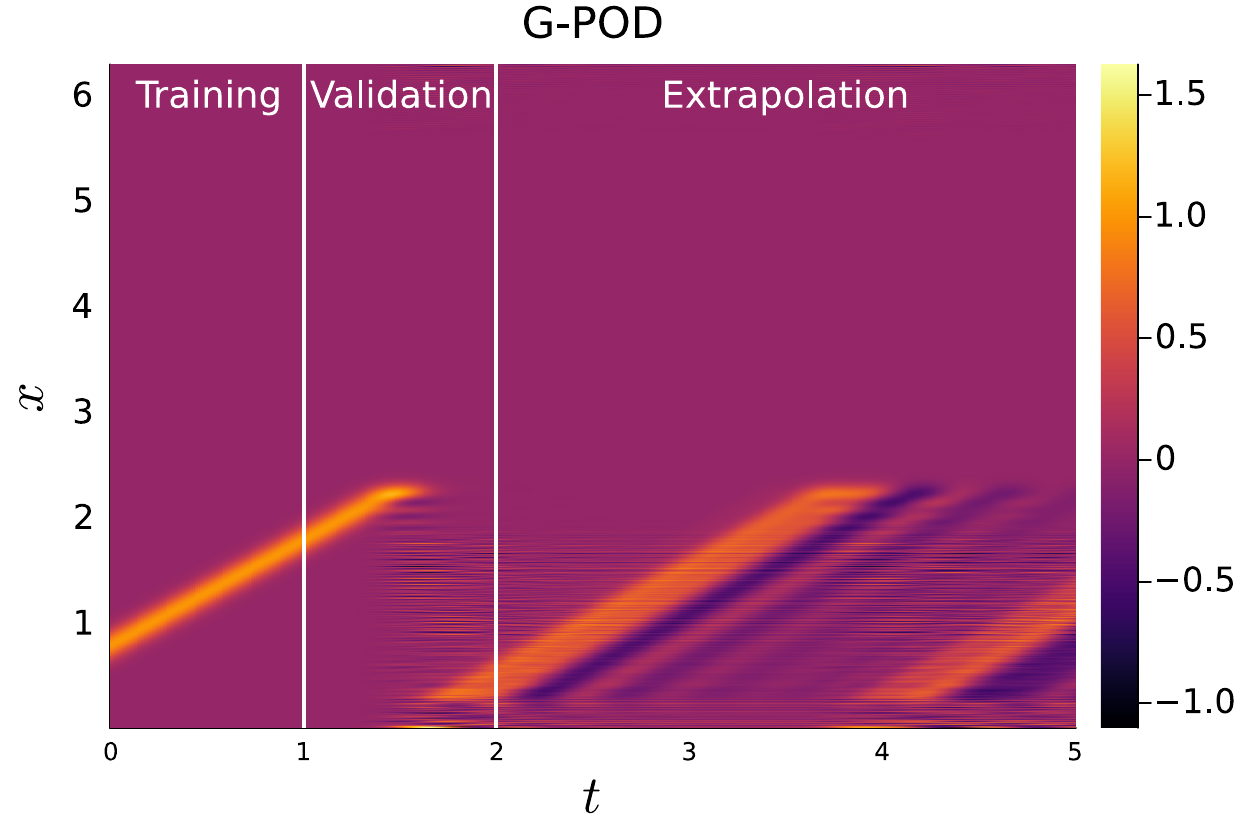}
    \includegraphics[width = 0.48\textwidth]{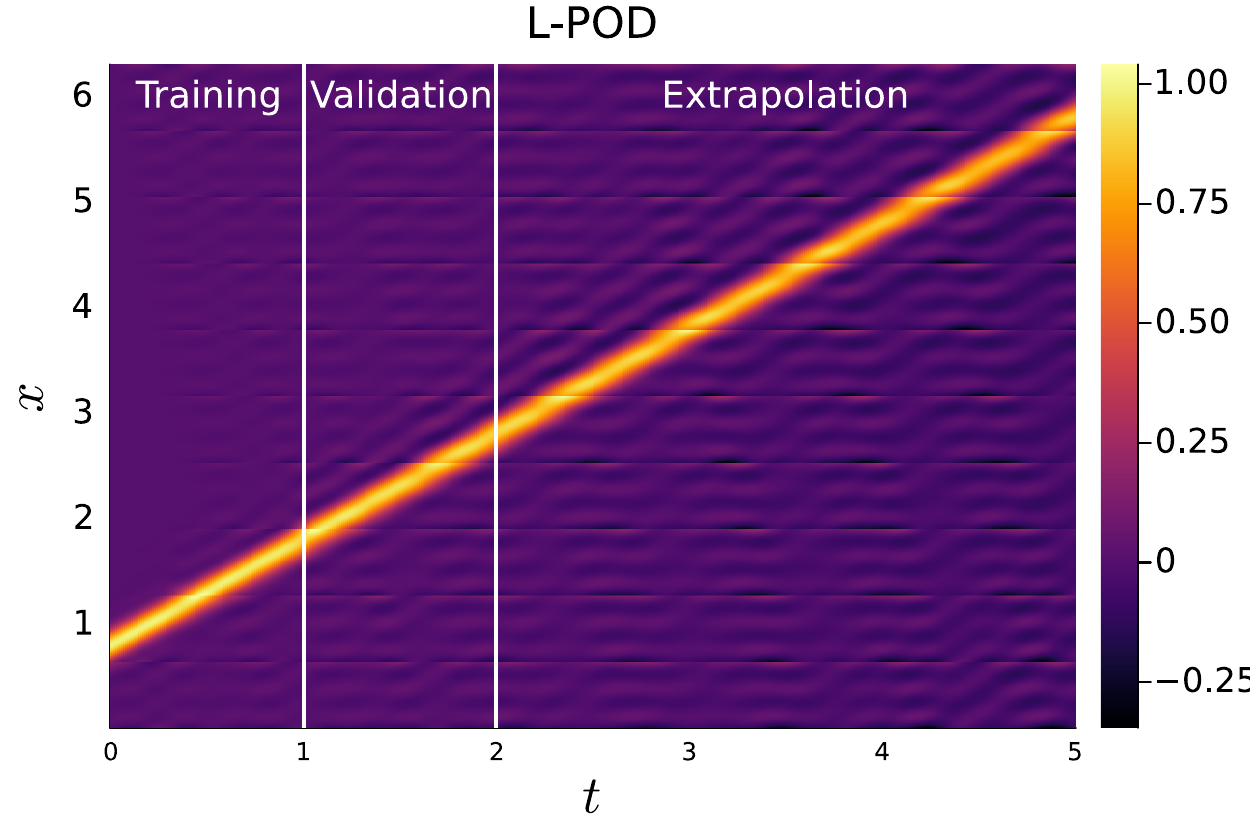}
    \includegraphics[width = 0.48\textwidth]{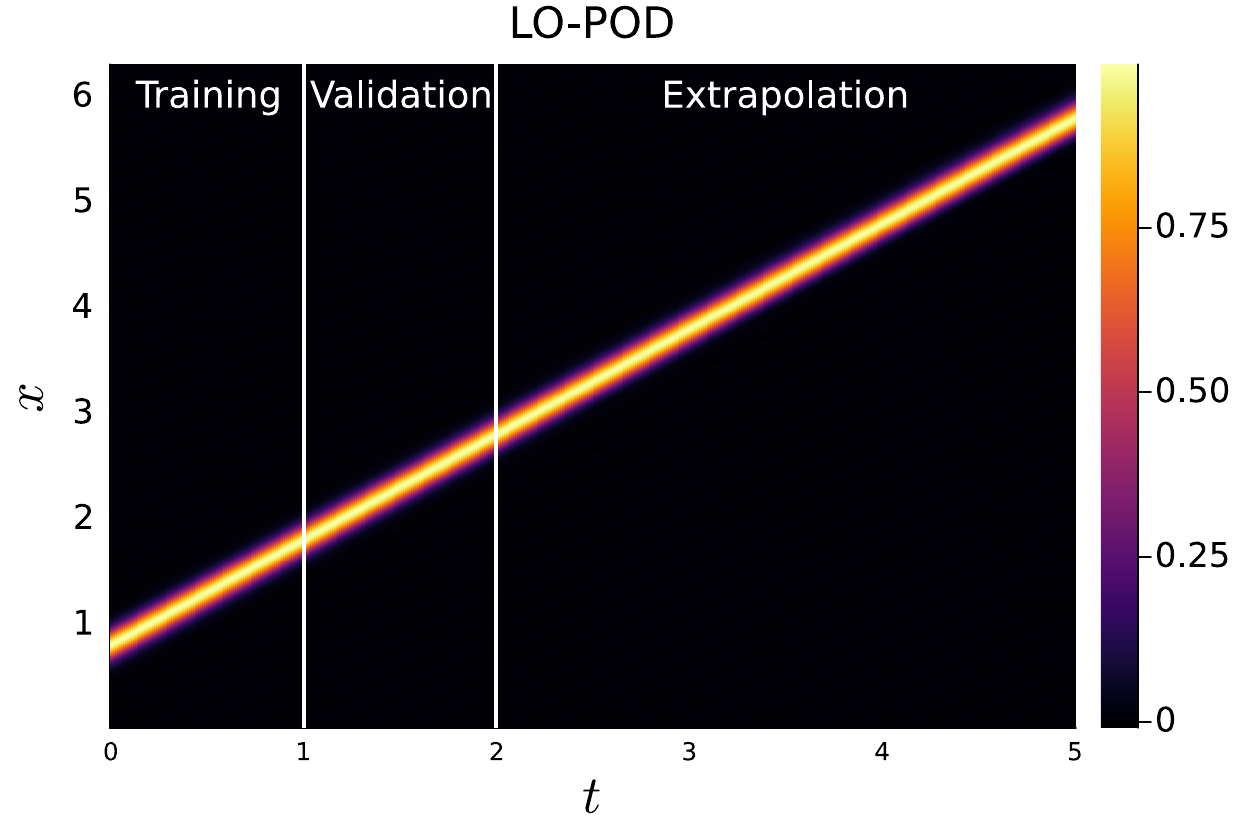}
    \includegraphics[width = 0.48\textwidth]{reference.pdf}
    \caption{Trajectories of $\mathbf{u}_r$ for the different \glspl{ROM}, along with the \gls{FOM} trajectory. From left to right, the white lines indicate the end of the training data and validation data, respectively.}
    \label{fig:trajectories}
\end{figure}

Having obtained a basis for each of the \gls{POD} approaches we construct a set of \glspl{ROM}. Based on the results in Section \ref{sec:proj_error}, we select $r=60$, $I=10$, and $q=6$ for the space-local approaches. To keep the size of the basis the same, we choose $r=60$ for \gls{G-POD}. The \gls{ROM} results are obtained by evaluating \eqref{eq:ROM} from the initial condition in \eqref{eq:init_cond}. The resulting simulations up to $t=5$ are presented in Figure \ref{fig:trajectories}.
For \gls{G-POD} we find that after the training data the performance degrades significantly. This exposes the limitations of \gls{G-POD}, as it is unable to adapt to wave traveling outside the training range. The space-local approaches perform better in this regard. Both \gls{L-POD} and \acrshort{LO-POD} are capable of extrapolating the traveling of the wave past the training region. However, \gls{L-POD} seems to suffer from discontinuities in $\mathbf{u}_r^\text{ROM}$, as the edges of the subdomains become increasingly visible as the simulation progresses. \acrshort{LO-POD} does not suffer from this issue and smoothly extrapolates the solution past the training region. This means the smooth reconstruction coming from \acrshort{LO-POD} indeed improves the quality of the simulation for the same number of \gls{DOF}.

\begin{figure}[ht]
    \centering
    \includegraphics[width = 0.48\textwidth]{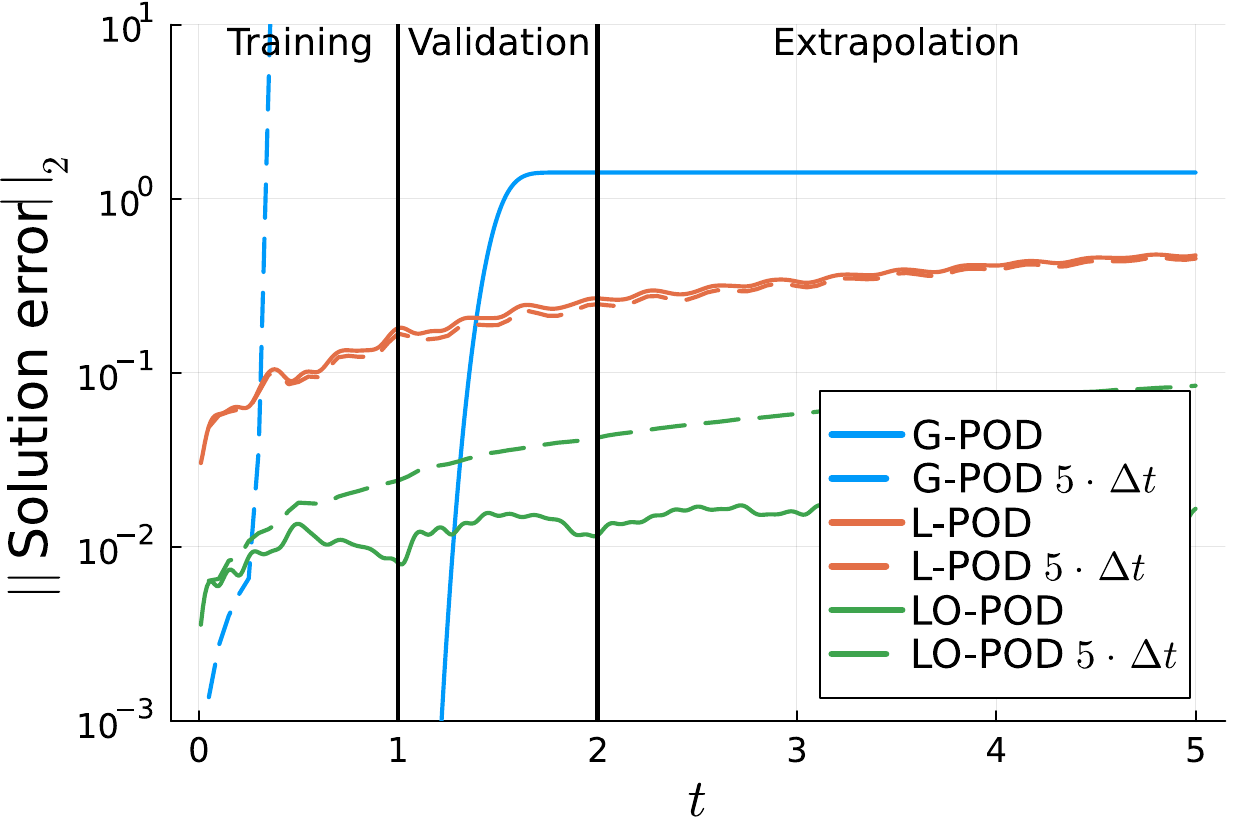}
    \includegraphics[width = 0.48\textwidth]{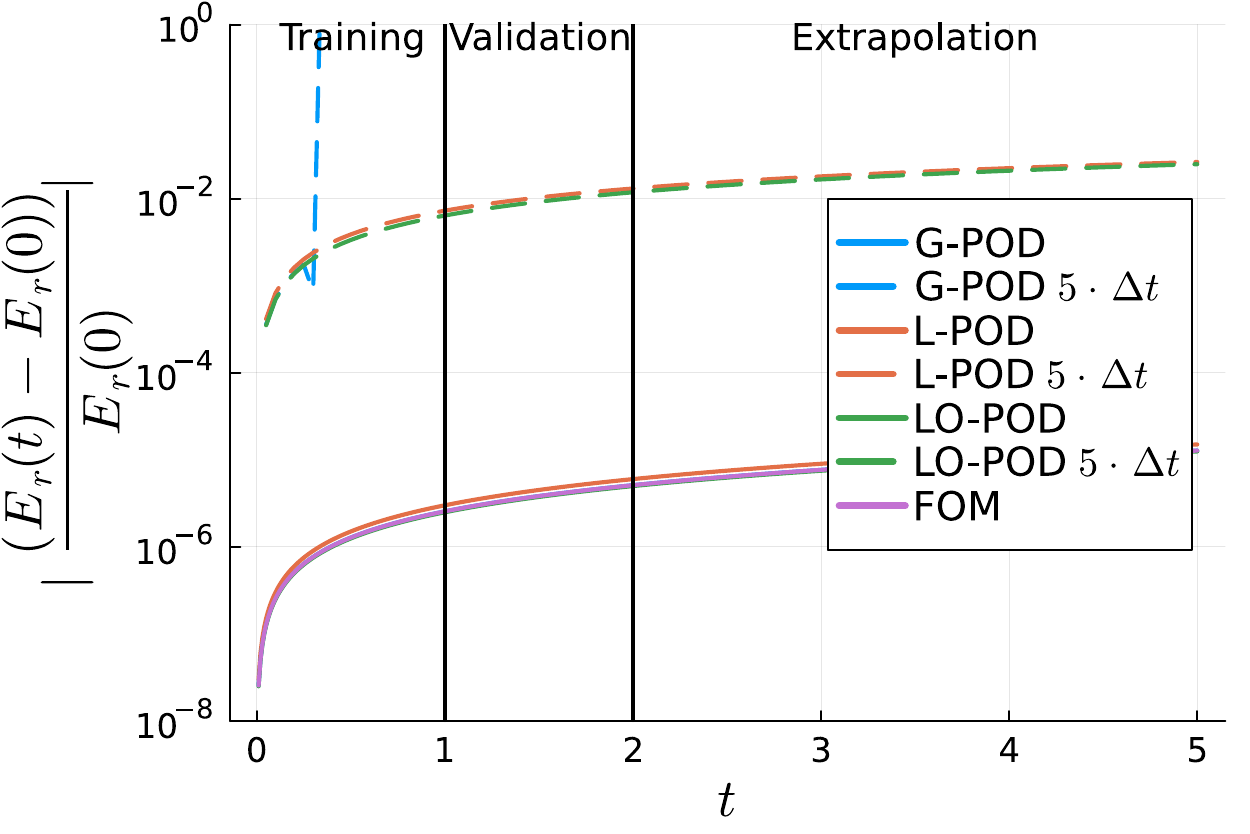}
    \includegraphics[width = 0.48\textwidth]{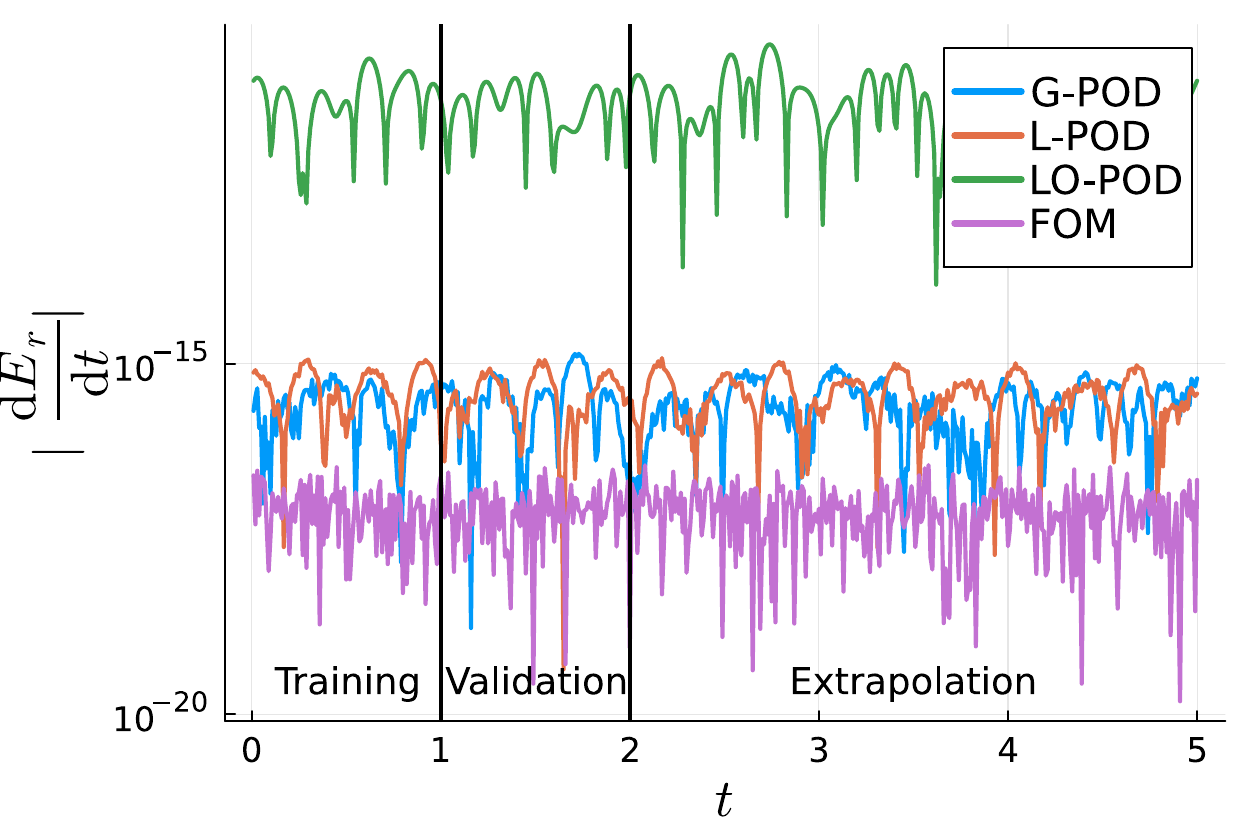}
    \includegraphics[width = 0.48\textwidth]{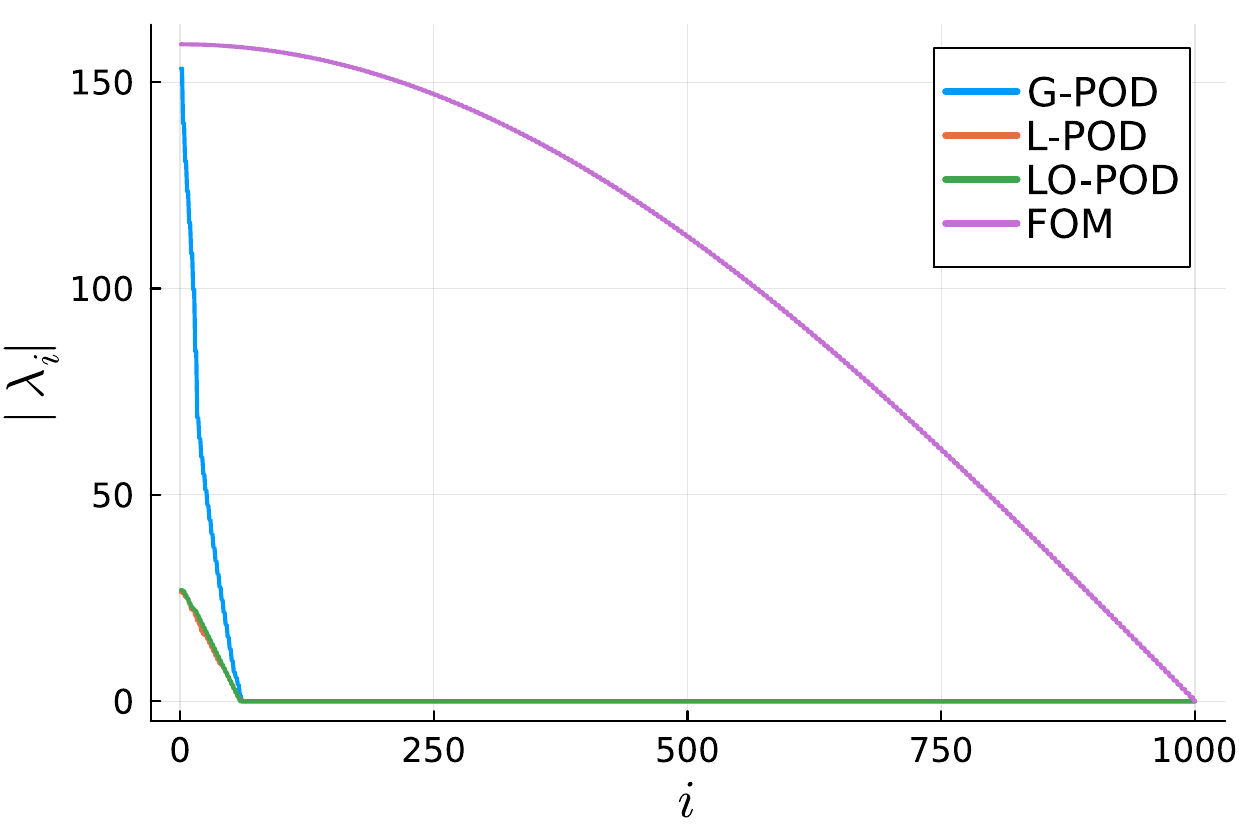}
    \caption{(Top-left) Solution error for each of \glspl{ROM} during the simulation, presented for both $\Delta t = 0.01$ and $\Delta t = 0.05$. From left to right, the black lines indicate the end of the training data and validation data, respectively. (Top-right) Relative change in energy with respect to the start of the simulation. (Bottom-left) Instantaneous change in energy, computed by evaluating \eqref{eq:change_in_energy}. (Bottom-right) Eigenvalues of the projected \gls{FOM} operator.}
    \label{fig:error}
\end{figure}

Next, we evaluate the \glspl{ROM} on a set of performance metrics. The first metric is the solution error, see \eqref{eq:error}. The other metrics (defined on the figure axes) focus on how well the energy is conserved during the simulation, namely the total change in energy and the instantaneous change in energy. The results are shown in Figure \ref{fig:error}.
The results are depicted for both the \gls{FOM} time step size $\Delta t$ (solid line) and a five times larger time step $5 \cdot \Delta t$ (dashed line). This larger time step is based on an eigenvalue analysis of the \gls{ROM} operator, also presented in Figure \ref{fig:error}. For this analysis, we determined the eigenvalues of the operator $\mathbf{D}$ projected on the \gls{ROM} basis. This operator is given by $\mathbf{P}\mathbf{D}\mathbf{P}$. It can be shown that replacing $\mathbf{D}$ in the \gls{FOM}, see \eqref{eq:semi-discrete}, by this projected operator is equivalent to integrating the \gls{ROM}, see \eqref{eq:ROM}. The eigenvalues $\lambda_i$ are ordered according to the magnitude of their absolute values. We observe that the largest eigenvalue for the space-local approaches is roughly five times smaller than that of \gls{G-POD} and the \gls{FOM}. This is likely due to the fact that the \gls{G-POD} basis functions contain higher frequencies than the space-local basis functions, see Figure \ref{fig:basis}. Using a smaller \gls{G-POD} basis could alleviate this.  Based on this analysis, we also evaluate the performance of the \glspl{ROM} for a five times larger time step \cite{eigenvalues_BALDAUF20086638}. 

Looking at the solution error, we observe that \gls{G-POD} performs best in the training region (the error is so small that it is outside the plotting range). However, after leaving the training region ($t>1$), the performance degrades rapidly. On the other hand, for the space-local approaches, we do not observe this jump of error outside the training region, but rather a steady increase. Importantly, outside the training region, the space-local approaches are much more accurate than \gls{G-POD}. In particular, \acrshort{LO-POD} improves by more than one order of magnitude upon both \gls{G-POD} and \gls{L-POD} outside the training region. 

When increasing the time step size by fivefold, the simulation quickly becomes unstable for \gls{G-POD}. For \gls{L-POD}, both time steps yield stable simulations, giving results that are very close to each other. For \acrshort{LO-POD} and a smaller time step size, the error seems to converge to an equilibrium. However, for a larger time step, it increases steadily. This means there is likely still a benefit to taking a smaller time step. An interesting continuation of this research would be to find a way to systematically determine the ``sweet spot" between increasing the time step and maintaining accuracy of the \gls{ROM}. \R{comment_7}\revone{This could for example, be done by considering several different time step sizes and both implicit and explicit integration methods, or adaptive time-stepping based on an error estimator \cite{Johnson1987}.}

In terms of energy conservation, we observe that all \glspl{ROM} conserve the energy as predicted by the theoretical analysis, except for a time discretization error incurred by the use of \gls{RK4}. This error increases when the time step size is increased. Only for \gls{G-POD} with an increased time step size, the simulation becomes unstable. For \acrshort{LO-POD}, the numerical error in the instantaneous change in energy, see \eqref{eq:change_in_energy}, is the largest. This is likely due to the linear system that needs to be solved to evaluate the \gls{ROM} for a non-orthogonal basis, see \eqref{eq:ROM}. However, the time discretization error is the primary source of error, as the change in energy during the simulation is roughly the same for \gls{L-POD} and \acrshort{LO-POD}.

\subsection{Convergence with increasing ROM dimension}\label{sec:convergence}

Next, we look at the convergence of the solution error as we increase the size of \gls{ROM} basis $r$. In addition, we compare the \glspl{ROM} to a finite difference discretization with the same number of \gls{DOF}. Note that for the \glspl{ROM} $\text{DOF} = r$. For the finite difference simulations, we first project the solution on the \gls{FOM} grid using linear interpolation. We then compute the difference, such as in \eqref{eq:error}, and compute the $L_2$-norm to quantify the error. For the local \glspl{ROM} we stick to $I=10$ subdomains, while increasing the number of basis functions per subdomain $q$ to increase $r$. We consider a maximum of $r=100$ for the \glspl{ROM}. 

\begin{figure}[ht]
    \centering
    \includegraphics[width = 0.48\textwidth]{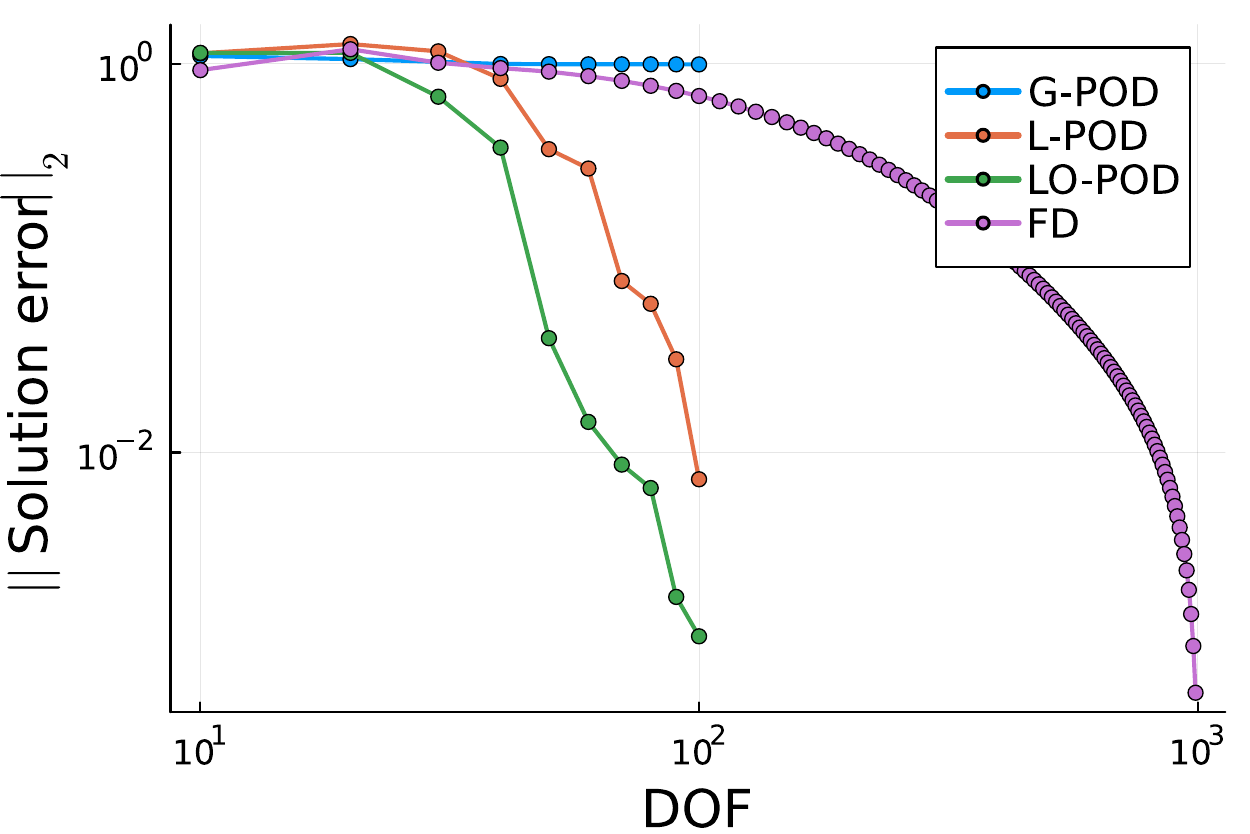}
    \includegraphics[width = 0.48\textwidth]{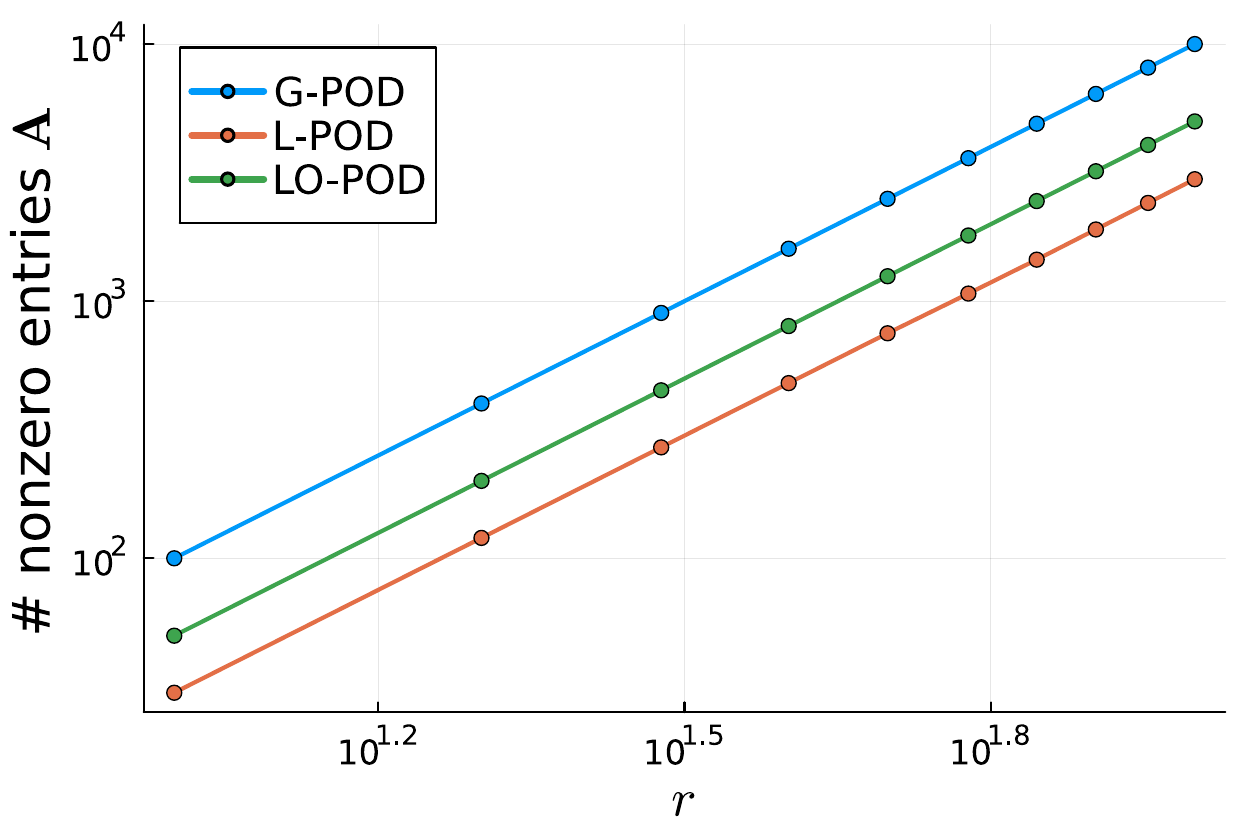}
    \caption{(Top-left) Solution error averaged over a simulation up to $t=5$. Results are presented for the different \glspl{ROM}, as well as a finite difference discretization (FD) with the same number of \gls{DOF}. (Bottom-right) Number of nonzero entries in the \gls{ROM} operator $\mathbf{A}$.}
    \label{fig:convergence}
\end{figure}

The results are depicted in Figure \ref{fig:convergence}.
Regarding the local \glspl{ROM}, we observe rapid convergence of the solution error as we increase the number of \gls{DOF}, with \acrshort{LO-POD} consistently outperforming \gls{L-POD}. This is in line with the projection error convergence discussed in Section \ref{sec:proj_error}. Regarding \gls{G-POD}, we observe no convergence of the solution error. Regarding the finite difference discretization, we find it converges at a much slower rate than the space-local \glspl{ROM}.

To obtain a conclusive answer about the scaling of the \glspl{ROM}, we considered the number of nonzero entries in the \gls{ROM} operator $\mathbf{A}$. This is also depicted in Figure \ref{fig:convergence}. Here we find that the number of nonzero entries scales according to a power law, with \gls{L-POD} scaling the most favorably. This is in line with the discussion presented in Section \ref{sec:Galerkin_projection}. Hyper-reduction can be carried out to further reduce the computational cost of the \glspl{ROM} \cite{hyperreduction_1,hyperreduction_2,KLEIN2024112697}. Finally, one could also increase the time step size, as discussed in Section \ref{sec:error_and_structure}, to further decrease the computational cost of the space-local \glspl{ROM}.

\subsection{2D advection-diffusion equation}

To further evaluate the viability of the space-local \glspl{ROM}, we consider the 2D advection-diffusion equation:
\begin{equation}
     \frac{\partial u(\mathbf{x},t)}{\partial t}  = \underbrace{-\nabla \cdot (\mathbf{V}(\mathbf{x}) u(\mathbf{x},t))}_{\text{advection}} + \underbrace{\nu \nabla^2   u(\mathbf{x},t)}_{\text{diffusion}},
\end{equation}
where some quantity $u(\mathbf{x},t)\in \mathbb{R}$ is advected through space $\mathbf{x} = \begin{bmatrix}
    x & y
\end{bmatrix}^T \in \mathbb{R}^2$ by a stationary velocity field $\mathbf{V}\in \mathbb{R}^2$ which is divergence free, i.e. $\nabla \cdot \mathbf{V} = 0$. The fact that the velocity field is now space-dependent, as opposed to the 1D case where it was uniform, adds to the difficulty of the test case. \R{comment_8}\revone{As opposed to the 1D case this system also contains diffusion. This results in a passive spread of $u$ throughout the domain.}  
The diffusion rate is determined by the scalar viscosity $\nu \geq 0$. In our case, we employ a finite volume discretization with a staggered grid for the velocity field \cite{benjamin_thesis}. This is done to preserve the skew-symmetry of the operator. The diffusive term is discretized using a simple second-order scheme. 

For the \gls{FOM} we consider the following setting: A periodic domain $\Omega = [-\pi,\pi] \times [-\pi,\pi]$, a velocity field given by $\mathbf{V}= \begin{bmatrix}
    \cos(x-y) & \cos(x-y)
\end{bmatrix}^T$ discretized on a $256 \times 256$ uniform grid, and a viscosity of $\nu = 10^{-3}$. The initial condition is given by the sum of three Gaussians, each centered on one of the streamlines of the flow where the velocity is highest: 
\begin{equation}
\begin{split}
    u(\mathbf{x},0) =& -\exp(-2x^2 -2 y^2) + \exp(-2\left(x-\frac{\pi}{2}\right)^2 -2\left(y+\frac{\pi}{2}\right)^2 )  \\ & + \exp(-2\left(x+\frac{\pi}{2}\right)^2 -2\left(y-\frac{\pi}{2}\right)^2 ).
\end{split}    
\end{equation}
The \gls{FOM} is integrated in time using a \gls{RK4} scheme with $\Delta t = 0.025$. Every 16 time steps, we save a snapshot of $u(x,t)$ for the construction of the \gls{ROM} basis. \R{comment_9}\revone{Note that the performance of the \glspl{ROM} converge to the \gls{FOM}, when doing the Galerkin projection. To resolve this closure, models can be included in the \gls{ROM} \cite{snyder2022reducedordermodelclosures}.}
The snapshots collected in the interval $t = [0,12]$ are used as training data, the interval $t = (16,20]$ is used as validation data, and the interval $t= (20,40]$ is used to test the extrapolation capabilities of the \glspl{ROM}. The \glspl{ROM} are also integrated in time explicitly using a \gls{RK4} scheme. \R{comment_11b}However, for the \acrshort{LO-POD} \gls{ROM} we resort to the implicit Crank-Nicolson method, see \cite{FEM_book}, to limit the computational cost. This scheme enables us to achieve second-order accuracy in time, while solving only a single linear system at each time step, rather than at each evaluation of \eqref{eq:ROM} in the \gls{RK4} scheme. To solve the resulting system efficiently, we employ an LU-decomposition \cite{MITTAL2002131_LU}.  This is required, as the \acrshort{LO-POD} basis is non-orthogonal, i.e. $\mathbf{S} \neq \mathbf{I}$.

Regarding the subdomain decomposition, we use the fact that the \gls{FOM} is discretized on a uniform grid. In this way, we can simply subdivide the domain uniformly into subdomains. For \gls{L-POD}, this results in a straightforward decomposition. For \acrshort{LO-POD} the decomposition is chosen such that each subdomain overlaps with its four neighboring subdomains which share a vertex in the middle of the considered domain. This is depicted in Figure \ref{fig:LO-decomp}.
\begin{figure}[ht]
    \centering
\begin{tikzpicture}

\def\op{0.4}
\def\offset{0.1}

\def\s{2}

\draw[blue, ultra thick] (0,0) -- (\s,0) -- (\s,\s) -- (0,\s) -- cycle;

\draw[green, ultra thick] (0,0) -- (-\s,0) -- (-\s,\s) -- (0,\s) -- cycle;

\draw[orange, ultra thick] (0,0) -- (-\s,0) -- (-\s,-\s) -- (0,-\s) -- cycle;

\draw[purple, ultra thick] (0,0) -- (\s,0) -- (\s,-\s) -- (0,-\s) -- cycle;

\fill[red, opacity=\op] (-1,-1) rectangle (1,1);
\draw[red, ultra thick] (-1,-1) rectangle (1,1);

\node at (-1,-1) [anchor=south west] {$\Omega_i$};

\end{tikzpicture}
    \caption{Subdomain $\Omega_i$ and its overlapping subdomains for \acrshort{LO-POD} in 2D.}
    \label{fig:LO-decomp}
\end{figure}

\R{comment_10b}The number of subdomains $I$ for the space-local approaches and the number of modes $q$ is chosen according to the performance on the validation set. 
This is presented in \ref{app:2D_optimization}. Based on this, we settle on $I = 8 \times 8$ and $q = 20$ for \gls{L-POD} and $q = 15$ for \acrshort{LO-POD}. We settle on these values for $q$ as these are the lowest values which surpass a projection error threshold of $10^{-2}$ on the validation set for this value of $I$. The value of $I$ is chosen as it generalizes well from the training data set to the validation data set. For \gls{G-POD}, we choose $r = q = 31$, as the snapshot matrix contains $31$ snapshots for \gls{G-POD} such that the \gls{SVD} only produces $31$ basis functions.
For the size of the time step, we choose the largest value that does not degrade the \gls{ROM} performance. This is presented in \ref{app:2D_time_step}. The time step sizes used are presented in Table \ref{tab:ROM_results_overview}.

\subsubsection{2D basis}

\begin{figure}[ht]
    \centering
    \includegraphics[width = 1\textwidth]{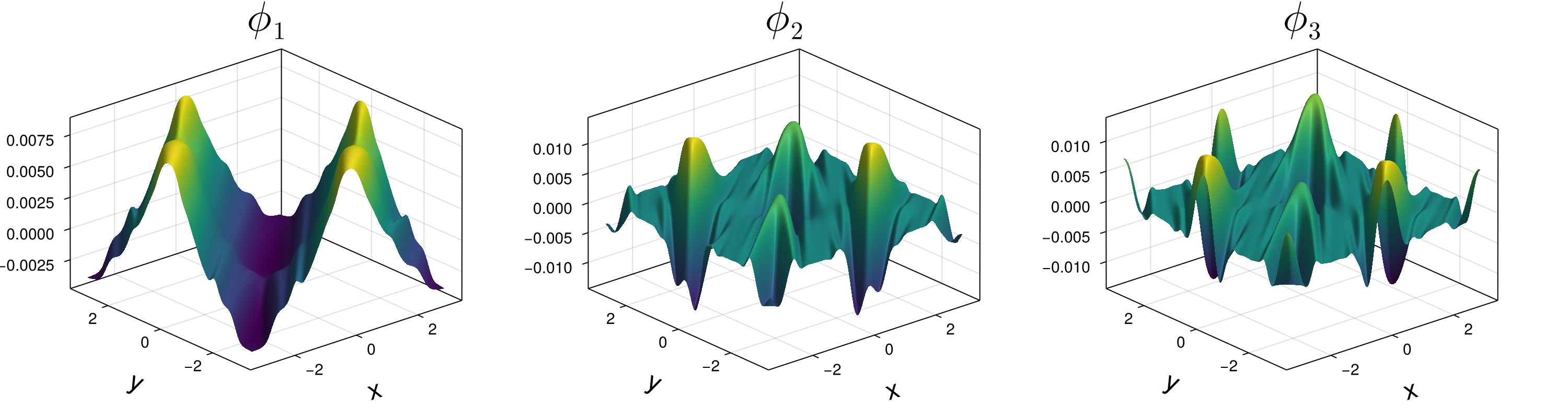}
    \includegraphics[width = 1\textwidth]{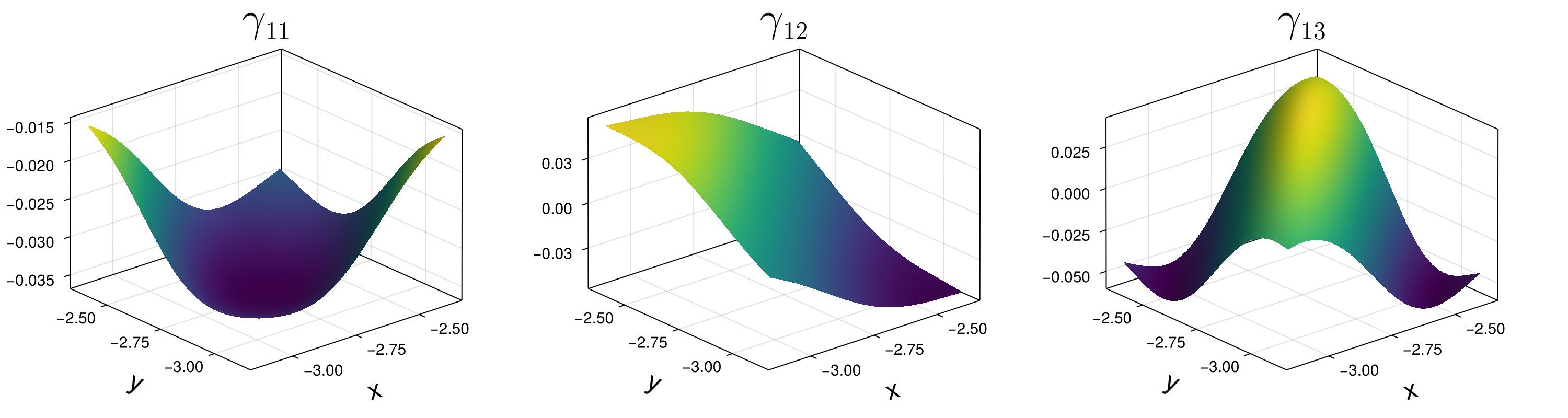}
    \includegraphics[width = 1\textwidth]{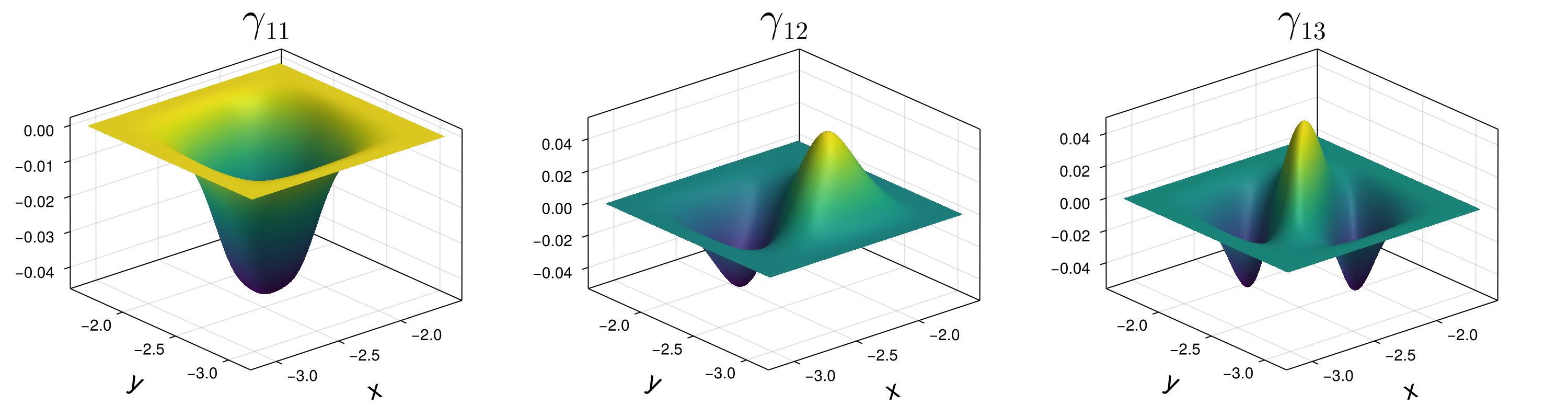}
    \caption{The first three \gls{G-POD} (top) modes and the first three space-local \gls{POD} modes for \gls{L-POD} (middle) and \acrshort{LO-POD} (bottom) for the 2D advection test case.}
    \label{fig:2D_basis}
\end{figure}

The resulting basis computed with the \gls{SVD} for each of the \gls{POD} approaches is depicted in Figure \ref{fig:2D_basis}.
From the basis functions, we can see that the flow runs diagonally across the domain. This is expected from the prescribed velocity field for which $\text{V}_1 = \text{V}_2$. The \gls{G-POD} modes contain information on the most dominant features of the flow. From these modes, we can see that some interesting dynamics occur in the region where the flow changes direction. The second and the third mode also seem to be quite similar, but slightly shifted. This is common for advection-dominated flows and results in a slow singular value decay \cite{brunton2019data}. For the space-local approaches, we find that the basis functions contain less information about the flow, but seem more similar to a generic basis. This likely allows them to generalize better on the validation set. We provide more insight on this in \ref{app:2D_optimization}. For \acrshort{LO-POD} we find that the basis functions smoothly decay to zero at the edge of the subdomain, as enforced by the kernel introduced in \ref{app:2D_kernel}. This allows \acrshort{LO-POD} to produce smooth approximations. The basis resulting from \acrshort{L-POD} does not guarantee this. 

\subsubsection{ROM performance}

In this section, we discuss the performance of different \glspl{ROM} and compare to both the \gls{FOM} and a finite volume discretization, which is twice as coarse as the \gls{FOM} in each direction ($128 \times 128$ as compared to $256 \times 256$ for the \gls{FOM}). This will be referred to as coarse \gls{FOM}. The results are summarized in Table \ref{tab:ROM_results_overview} for a full simulation up to $t=40$ (including the training, validation, and extrapolation region).

\begin{table}[]
\caption{Computation times (\textbf{comp time}) in seconds on a standard laptop CPU and average solution errors (\textbf{sol err}) computed for the entire simulation $t \in [0,40]$ for each of the \glspl{ROM} for the 2D advection test case. A coarse \gls{FOM} with a resolution of $128 \times 128$ is also included for comparison. For the space-local approaches, the number of \gls{DOF} is reported as the multiplication $I  \times q = r$. Reported computation times are an average of $5$ replica simulations. The average error for the gradient of the solution (\textbf{grad err}) is also depicted. }
\begin{footnotesize}
\begin{tabular}{ccccccc}
\toprule
 & \textbf{DOF} ($I \times q$) & \textbf{comp time} & \textbf{$||\text{sol err} ||_2$} & \textbf{$|| \text{grad err} ||_{2}$}  & $\Delta t$ \\ 
 \midrule
\gls{G-POD} & \textbf{31} & $\boldsymbol{2.5 \times 10^{-3}}$ & 0.267\phantom{0} &  0.636 & \textbf{0.2\phantom{00}} \\ 
\gls{L-POD} & $8^2 \times 20 = 1280$ & 0.207\phantom{00} & \textbf{0.0353} & 0.213 & \textbf{0.2\phantom{00}} \\ 
\acrshort{LO-POD} & $8^2 \times 15 = 960$ & 1.02\phantom{000} & 0.0354 &  \textbf{0.164} & 0.05\phantom{0} \\ 
coarse \gls{FOM} & $128^2 = 16384$ & 4.57\phantom{000} & 0.226\phantom{0} & 0.514 & 0.1\phantom{00} \\ 
\gls{FOM} & $256^2 = 65536$ & 25.9\phantom{000} & - & - & 0.025 \\ 
\bottomrule
\end{tabular}\label{tab:ROM_results_overview}
\end{footnotesize}
\end{table}

For the \gls{G-POD} \gls{ROM} we find a short computation time, due to the limited number of modes, but also a large solution and gradient error. The latter is defined as
\begin{equation}
    \text{gradient error} := \frac{\mathbf{G} \mathbf{u}^\text{ROM}_r(t) -\mathbf{G} \mathbf{u}^\text{FOM}(t)}{||\mathbf{G} \mathbf{u}^\text{FOM}(t)||_2},
\end{equation}
where $\mathbf{G} \in \mathbb{R}^{2N \times N}$ is a forward difference approximation of the gradient operator, such that the elements of $\mathbf{G}\mathbf{u}_r^\text{ROM}$ approximate $\nabla u$ in different parts of the domain. These large errors are likely the result of the \gls{G-POD} basis's limited generalization capabilities. For the space-local approaches, we observe that the computation time is roughly two orders of magnitude larger than for \gls{G-POD}. \R{comment_10a}\revone{As we have access to more space-local \gls{POD} modes than global ones, we choose to include enough to accurately represent the solution, see \ref{app:2D_optimization} for details. This increases the computation time with respect to \gls{G-POD}, see Table \ref{tab:ROM_results_overview} for the exact number of \gls{DOF} for each \gls{ROM}. } However, using the space-local approaches, we still obtain a speedup between one and two orders of magnitude with respect to the \gls{FOM}. The \gls{L-POD} \gls{ROM} is roughly 5 times faster than the \acrshort{LO-POD} \gls{ROM}. This is mostly because the explicit \gls{RK4} scheme (fourth order accurate) used for \gls{L-POD} allows for larger time steps without loss of accuracy, as compared to the Crank-Nicolson scheme (only second order accurate) used for \acrshort{LO-POD}, see \ref{app:2D_time_step} for an analysis on the time step size. \R{comment_11a}\revone{As mentioned earlier, the Crank-Nicolson scheme is used for \acrshort{LO-POD} as it limits the number of linear systems needed to be solved per time step to one}. \gls{L-POD} and \acrshort{LO-POD} both produce a similar solution error. However, looking at the gradient error \acrshort{LO-POD} obtains a lower error. This is likely because \acrshort{LO-POD} does not suffer from discontinuities at the subdomain boundaries, while \gls{L-POD} does, see Figure \ref{fig:discontinuity}. Both space-local approaches outperform the coarse grid finite volume discretization, both in computation time and error metrics.

\begin{figure}[ht]
    \centering
    \includegraphics[width = 0.32\textwidth]{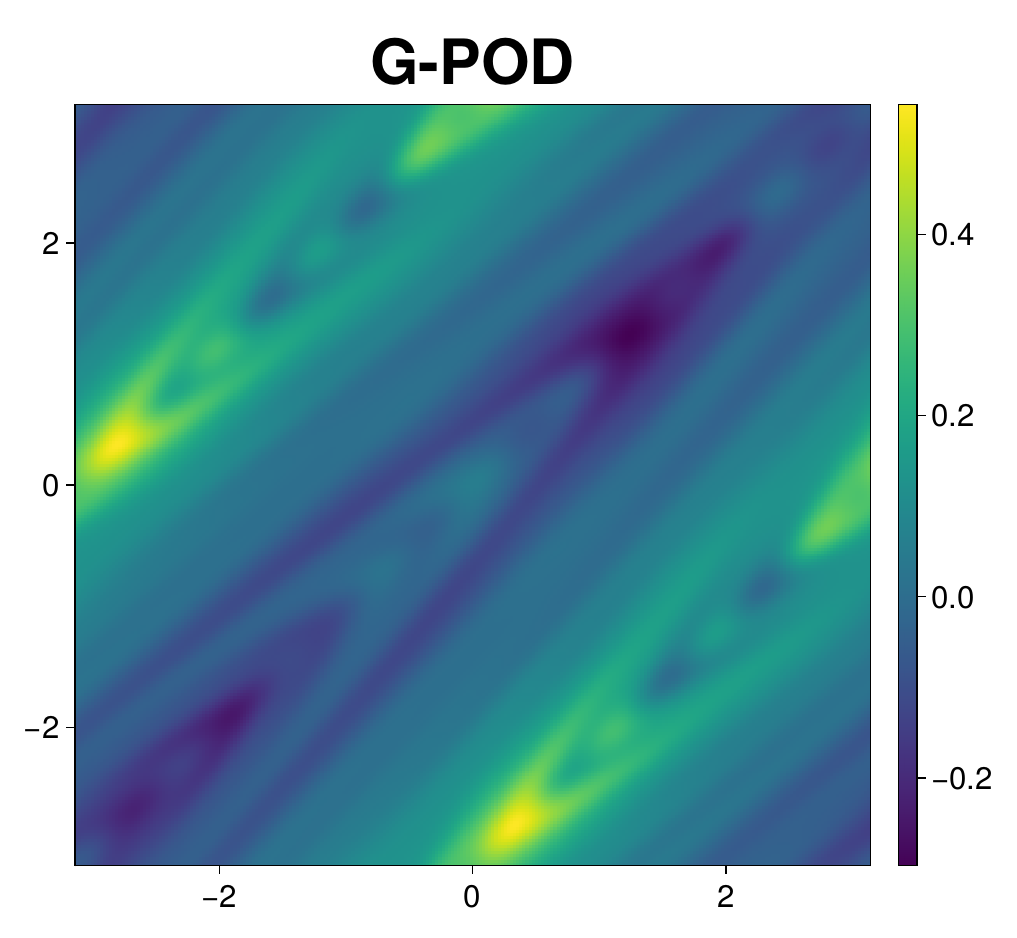}
    \includegraphics[width = 0.32\textwidth]{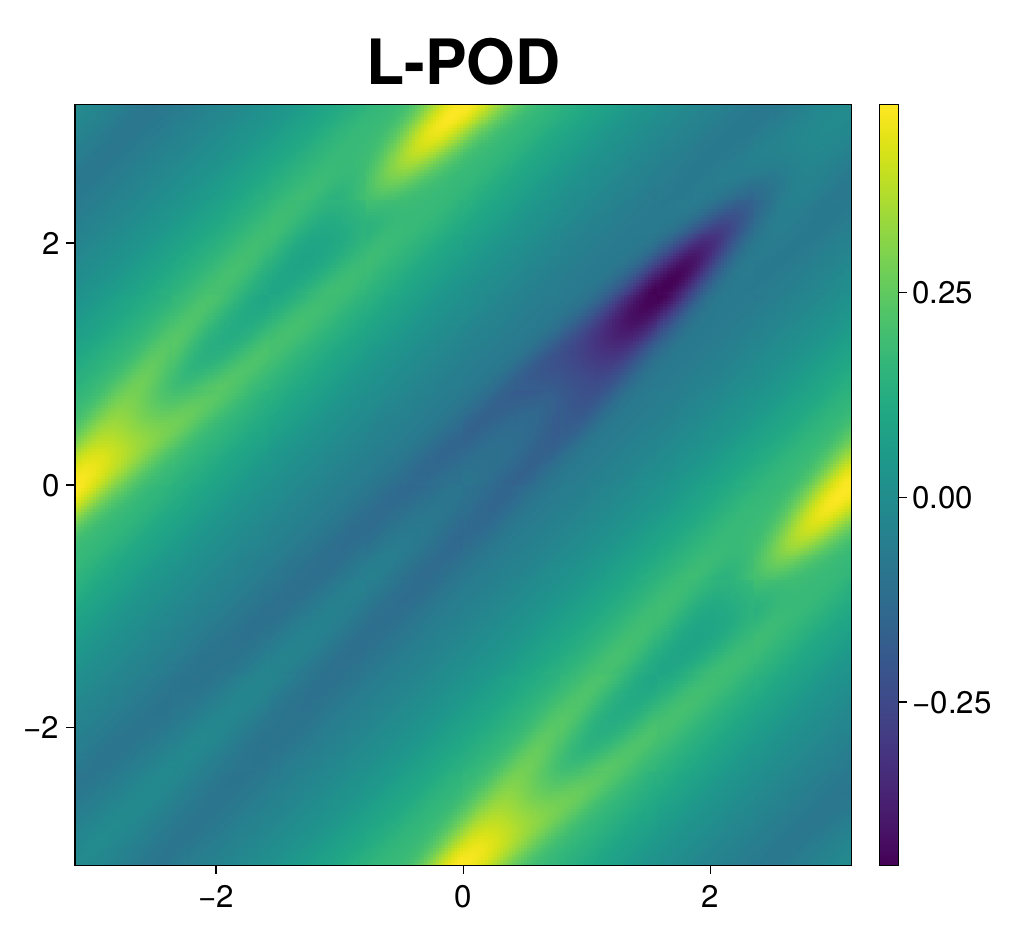}
    \includegraphics[width = 0.32\textwidth]{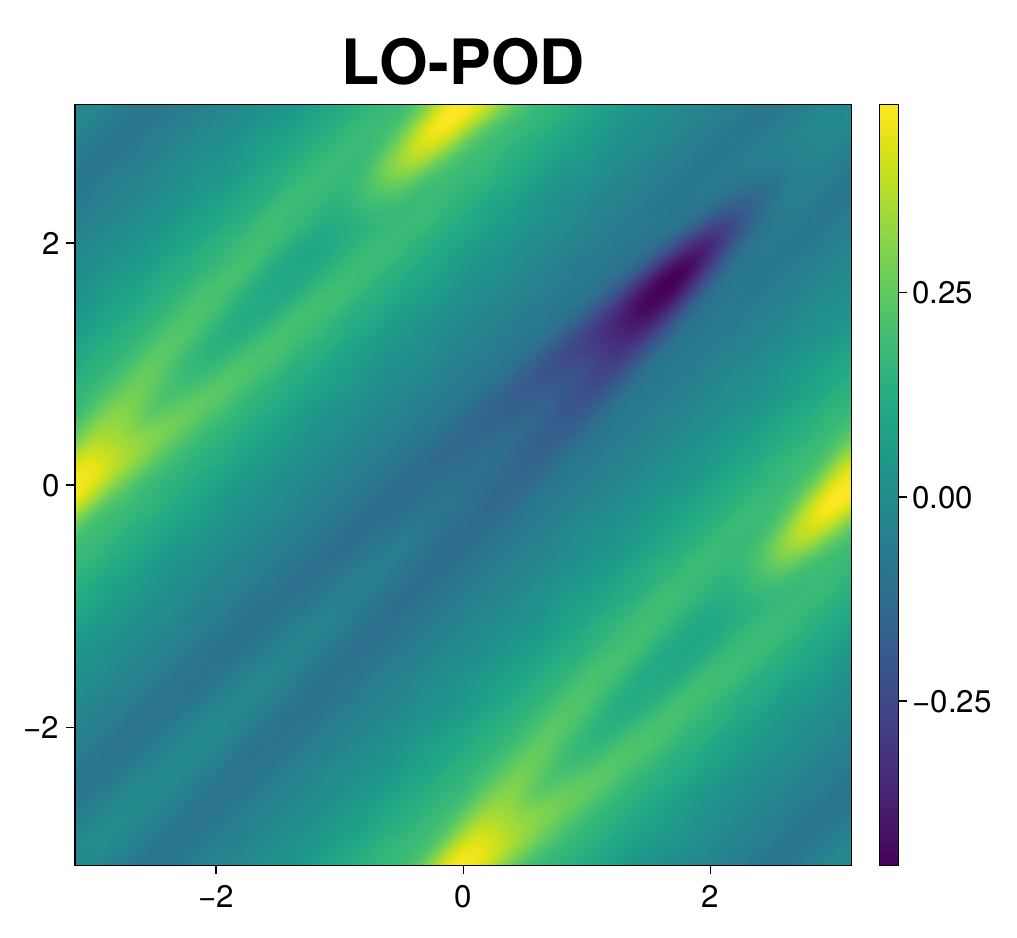}
    \includegraphics[width = 0.32\textwidth]{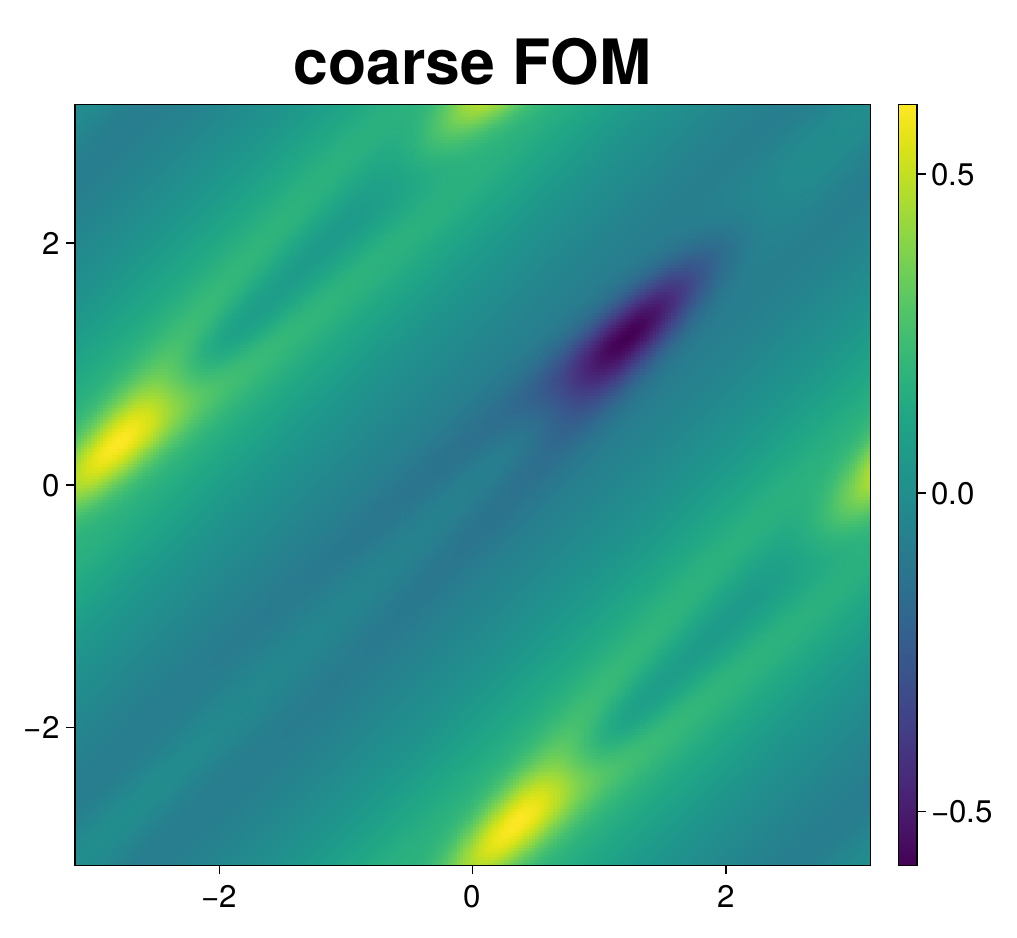}
    \includegraphics[width = 0.32\textwidth]{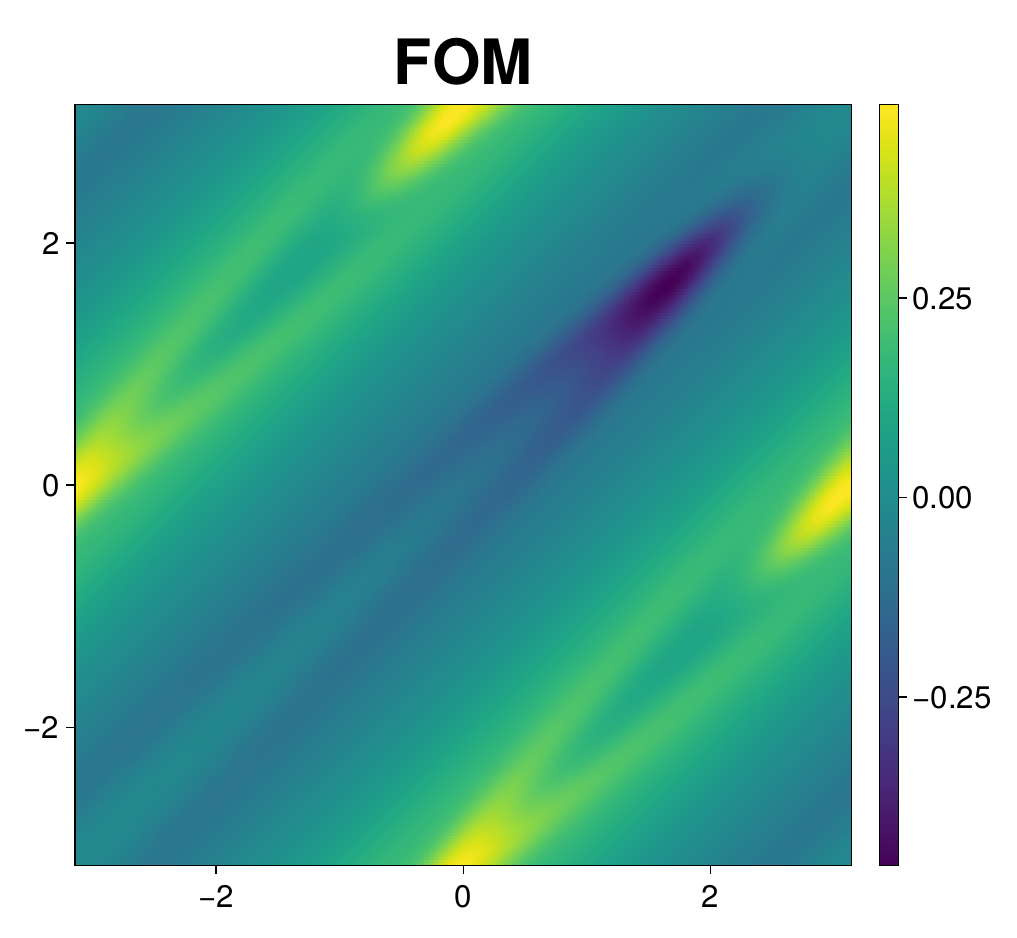}
    \caption{Solution at the end of the simulation at $t=40$ for the 2D advection test case produced by the different \glspl{ROM} as well as a coarse \gls{FOM}.}
    \label{fig:2D_heatmaps}
\end{figure}

Next, we consider the solution at the end of the simulation ($t=40$), see Figure \ref{fig:2D_heatmaps}. 
Here we find that \gls{G-POD} does not capture the solution well and contains unphysical artifacts, whereas the space-local approaches do reproduce the solution well. The coarse-grid finite volume discretization also captures the solution quite well, but is slightly smoothed and shifted (due to numerical diffusion and dispersion).

\begin{figure}[ht]
    \centering
    \includegraphics[width = 0.48\textwidth]{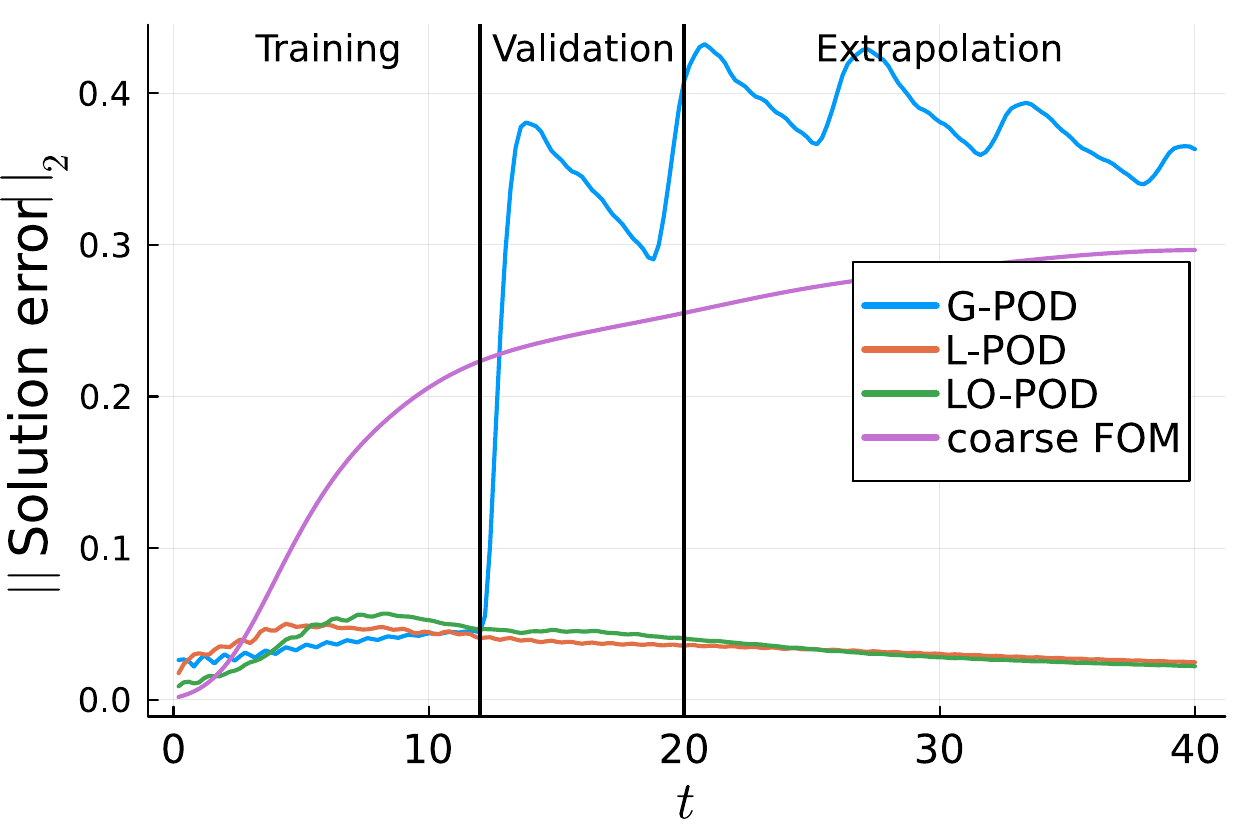}
    \includegraphics[width = 0.48\textwidth]{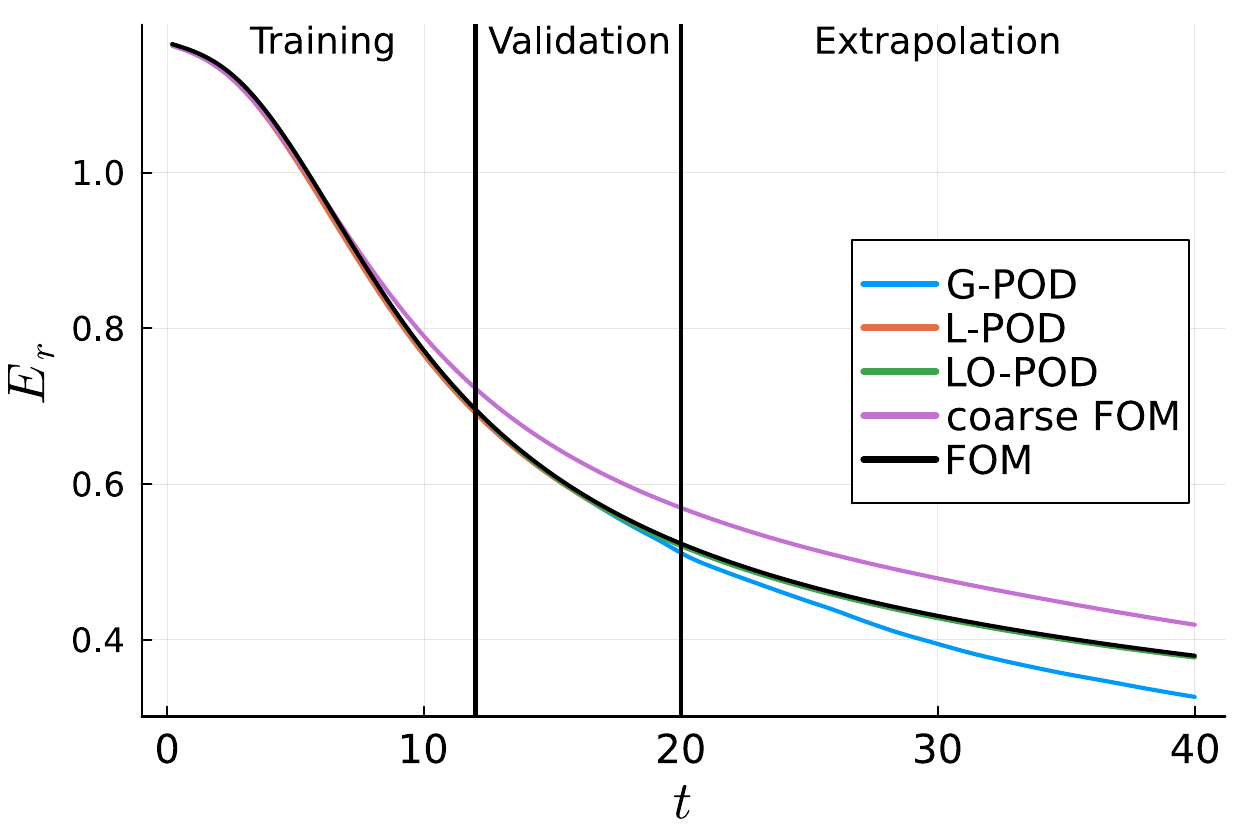}
    \caption{(Left) Solution error for each of the \glspl{ROM} during the simulation of the 2D advection test case. (Right) Energy trajectories during the simulation. Energy trajectories from \gls{L-POD} and \acrshort{LO-POD} overlap with the \gls{FOM} trajectory. From left to right, the black lines indicate the end of the training data and validation data, respectively.}
    \label{fig:2D_trajectories}
\end{figure}

Finally, we consider the solution error and energy trajectory for each of these models, see Figure \ref{fig:2D_trajectories}.
 We find that the \gls{G-POD} \gls{ROM} performs well within the training region, but it does not extrapolate well. For the space-local approaches, this is not an issue as both perform comparably, both inside and outside the training region, and for both the solution error and the energy. These approaches also outperform the coarse \gls{FOM}, which slowly accumulates error during the simulation. 

\subsection{\R{editor_a}\revtwo{Applicability to 2D Navier-Stokes}}\label{sec:navier-stokes}

As a final test case, we evaluate the applicability to the 2D incompressible Navier-Stokes equations by assessing how well the \gls{POD} bases generalize to representing the solution.
The incompressible Navier-Stokes equations read:
\begin{subequations}\label{intro:eq:navier_stokes}
\begin{align}
  \frac{\partial\mathbf{V}}{\partial t} + (\mathbf{V}\cdot\nabla)\mathbf{V} &= -\frac{1}{\rho}\nabla p + \nu\Delta\mathbf{V} + \mathbf{f}, \label{intro:eq:ns_momentum}\\[4pt]
  \nabla\cdot\mathbf{V} &= 0, \label{eq:ns_incompressibility}
\end{align}
\end{subequations}
for a 2D velocity field $\mathbf{V}(\mathbf{x},t)\in \mathbb{R}^2$. For the density we choose $\rho = 1$ and for the kinematic viscosity $\nu = \frac{1}{1000}$. For the forcing we take $\mathbf{f}(\mathbf{V},\mathbf{x},t) = \begin{bmatrix}
     \sin(4y) & 0
\end{bmatrix}^T - 0.1\mathbf{V} $, as to simulate Kolmogorov flow. This setup is often used for evaluating data-driven turbulence modeling strategies \cite{Kochkov_2021,shankar2024differentiableturbulenceclosurepartial,Agdestein2025MLLES}. The \gls{PDE} is solved numerically using the second-order accurate scheme presented in \cite{harlow1965numerical} on a periodic domain $\Omega = [-\pi,\pi] \times [-\pi,\pi]$. The computational grid consists of $2048 \times 2048$ grid cells. The simulation is initialized with energy content in the large wave-numbers, and a warm-up period of 25 time units is employed before collecting the training data. After the warm-up period, training data is collected on the interval $t = [0,10]$ and validation data on the interval $t = [100,110]$, to remove any temporal correlation. Both the validation and training data consist of 50 snapshots of the vorticity $\boldsymbol{\omega} = \nabla \times \mathbf{V} \in \mathbb{R}^3$. For 2D, only the third component of this vector is nonzero, which we shall refer to simply as $\omega$. We choose to represent the vorticity as for periodic domains in 2D, all the information of the flow is contained within this scalar field \cite{lequeurre2018vorticitystreamfunctionformulations}.

\begin{figure}[ht]
    \centering
    \includegraphics[width = 0.48\textwidth]{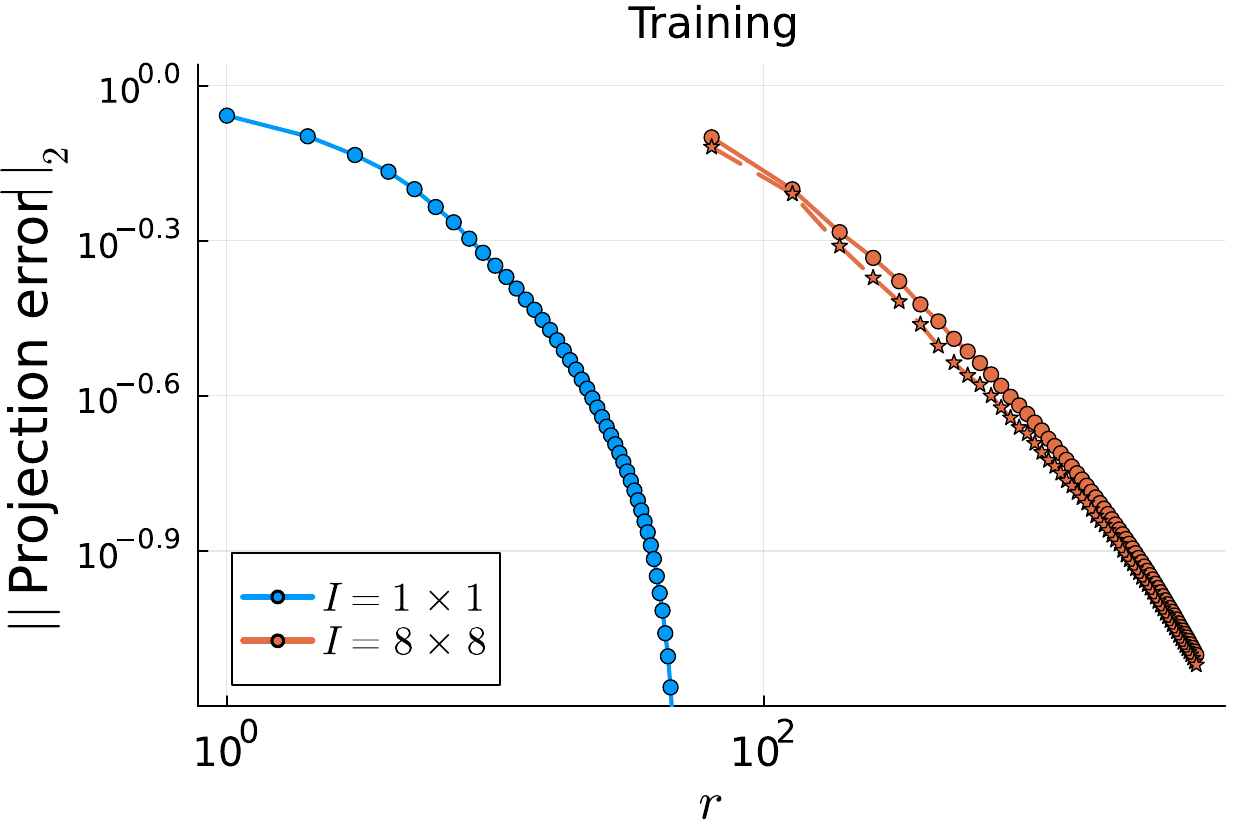}
    \includegraphics[width = 0.48\textwidth]{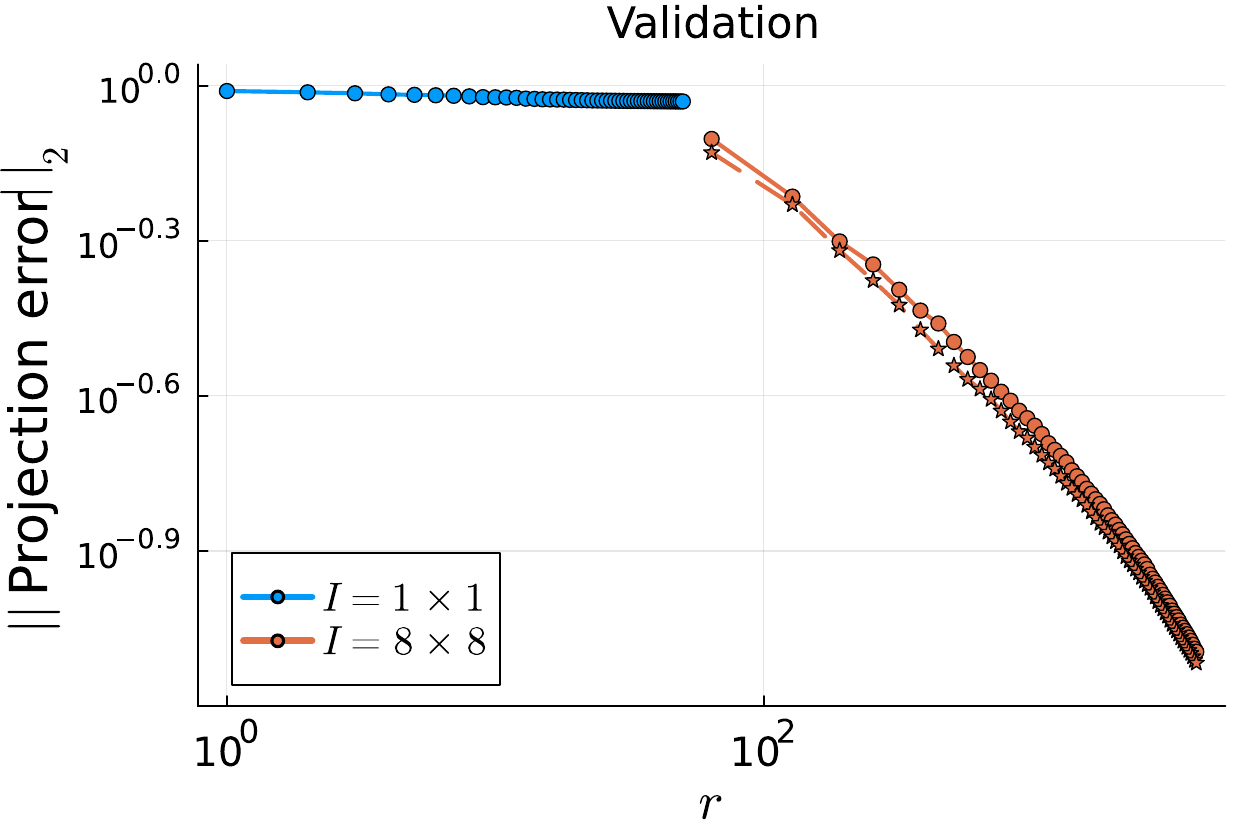}
    \caption{Projection error evaluated over the training (left) and validation (right) data set for the 2D Navier-Stokes test case. $I=1\times 1$ corresponds to \gls{G-POD}. Solid dots correspond to \gls{L-POD}, stars to \acrshort{LO-POD}.}
    \label{fig:2D_NS_convergence}
\end{figure}

To begin our analysis, we examine the projection error for both the training and validation datasets using an $8 \times 8$ decomposition of subdomains, which is identical to the 2D advection-diffusion case, see Figure \ref{fig:2D_NS_convergence}.
Similar to the previous test cases, we find that the projection error of the \gls{G-POD} converges fast on the training set. However, upon examining the validation set, we find that the space-local approaches generalize significantly better. 

\begin{figure}[ht]
    \centering
    \includegraphics[width = 0.32\textwidth]{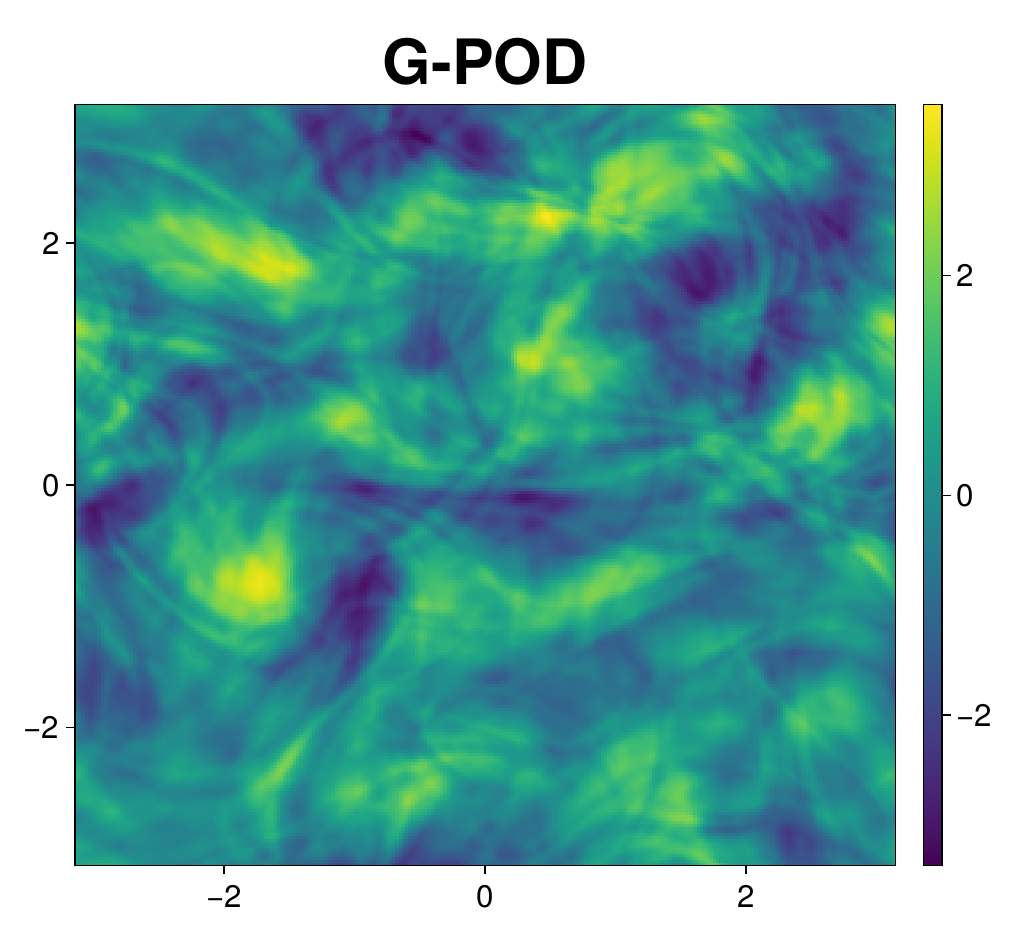}
    \includegraphics[width = 0.32\textwidth]{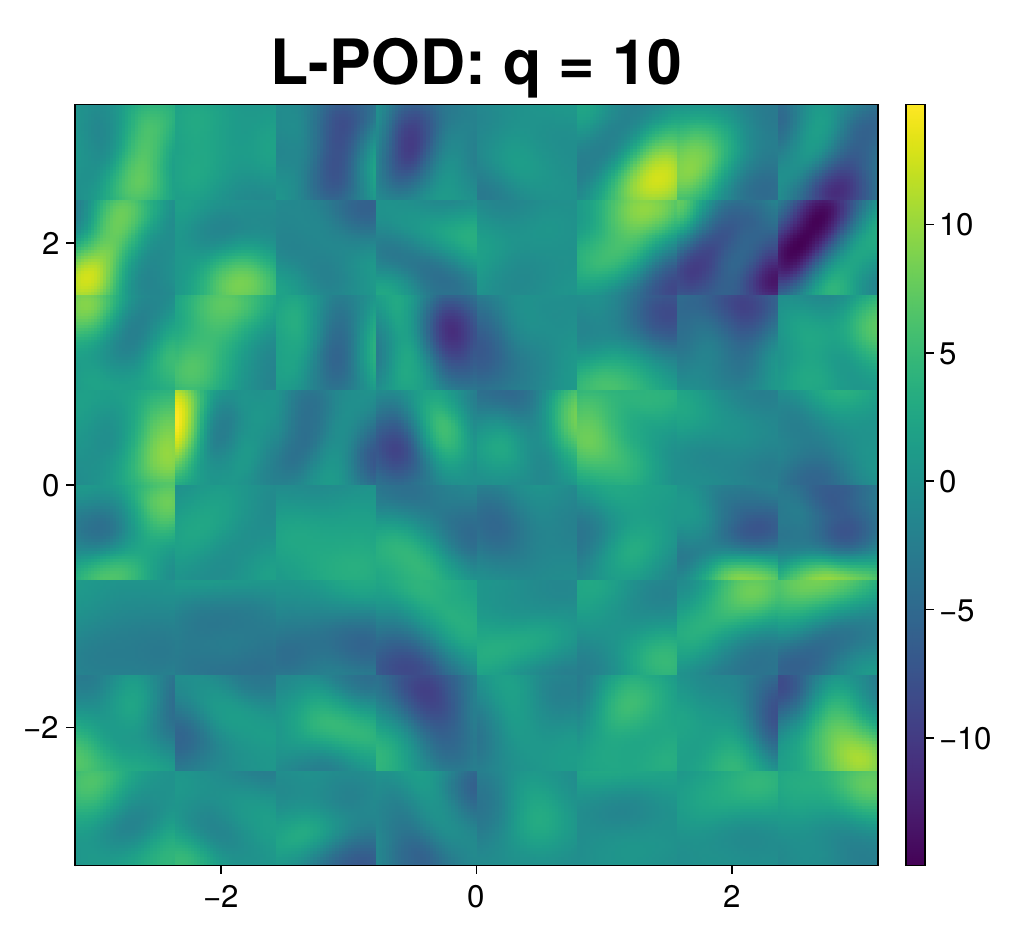}
    \includegraphics[width = 0.32\textwidth]{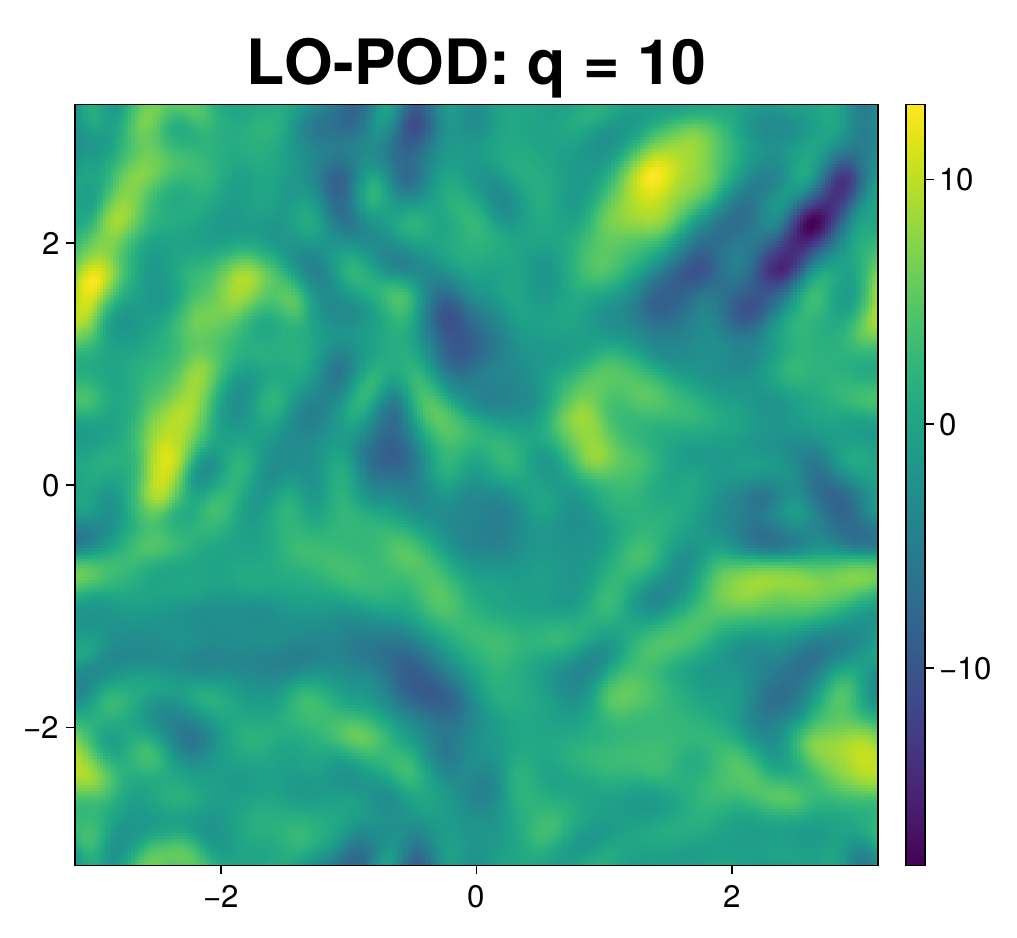}
    \includegraphics[width = 0.32\textwidth]{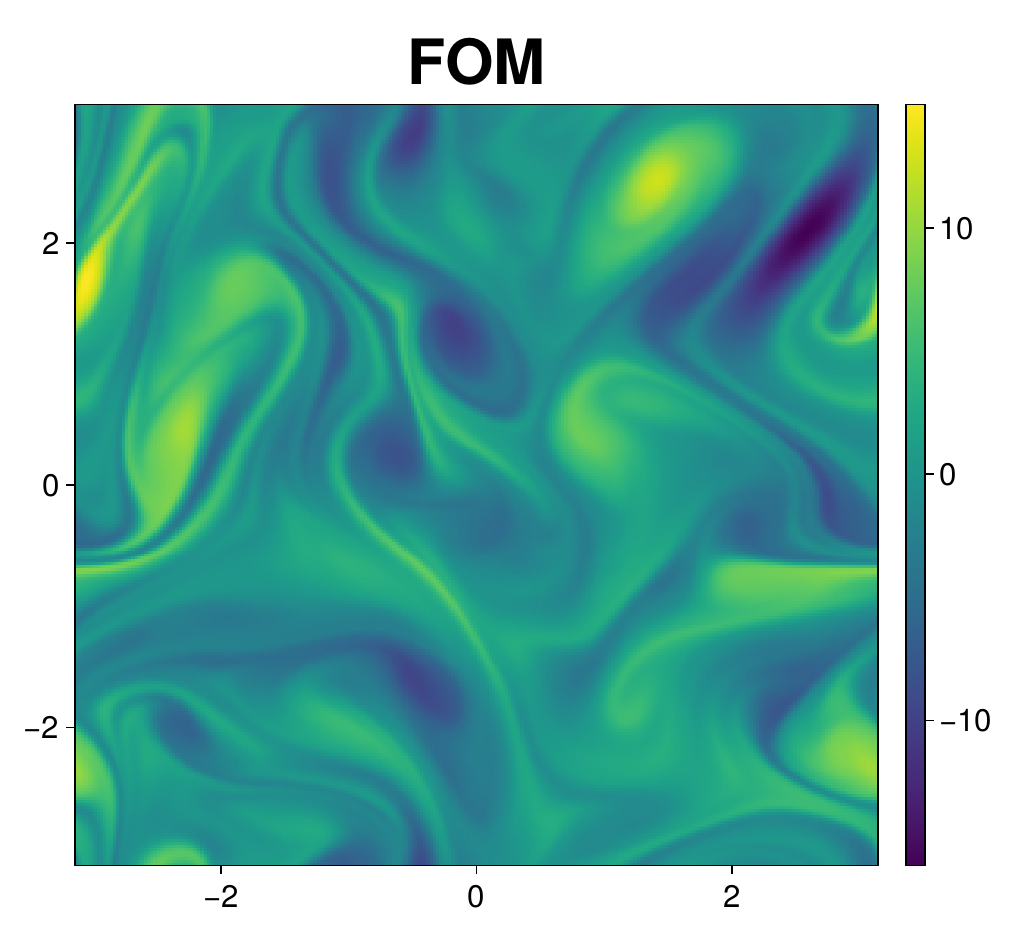}
    \includegraphics[width = 0.32\textwidth]{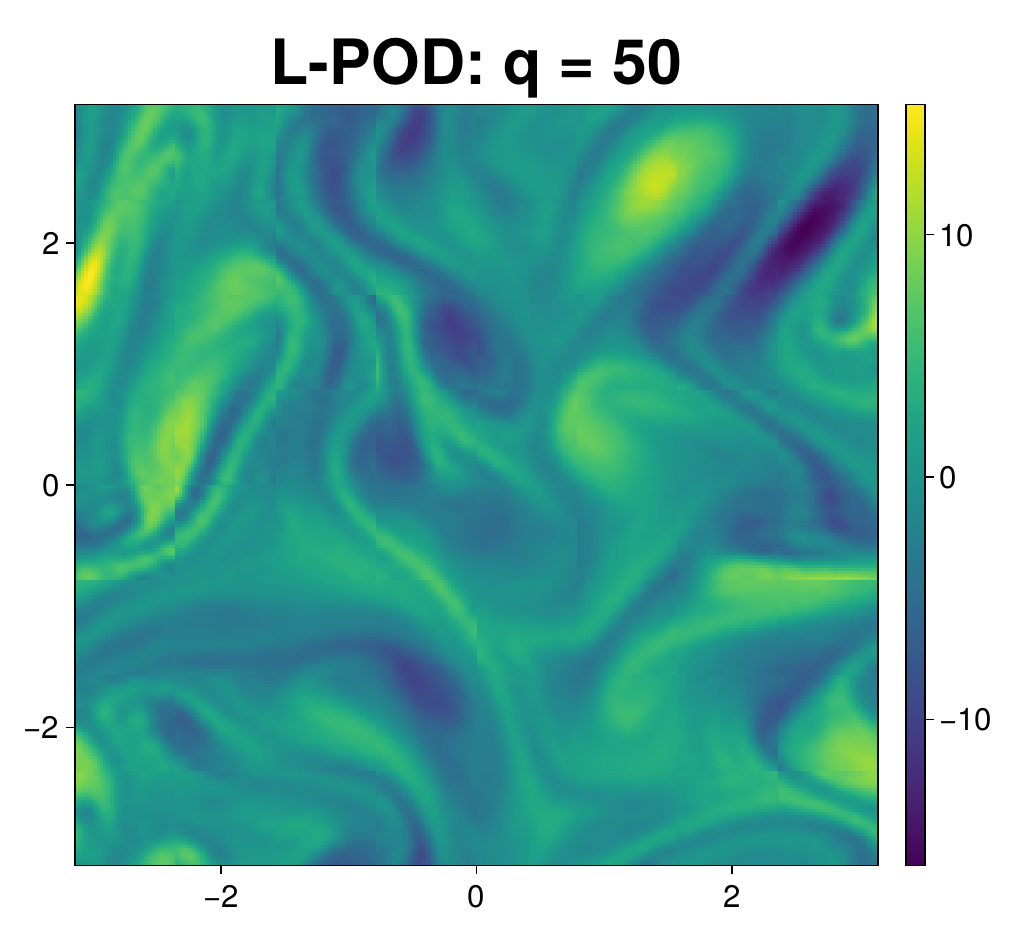}
    \includegraphics[width = 0.32\textwidth]{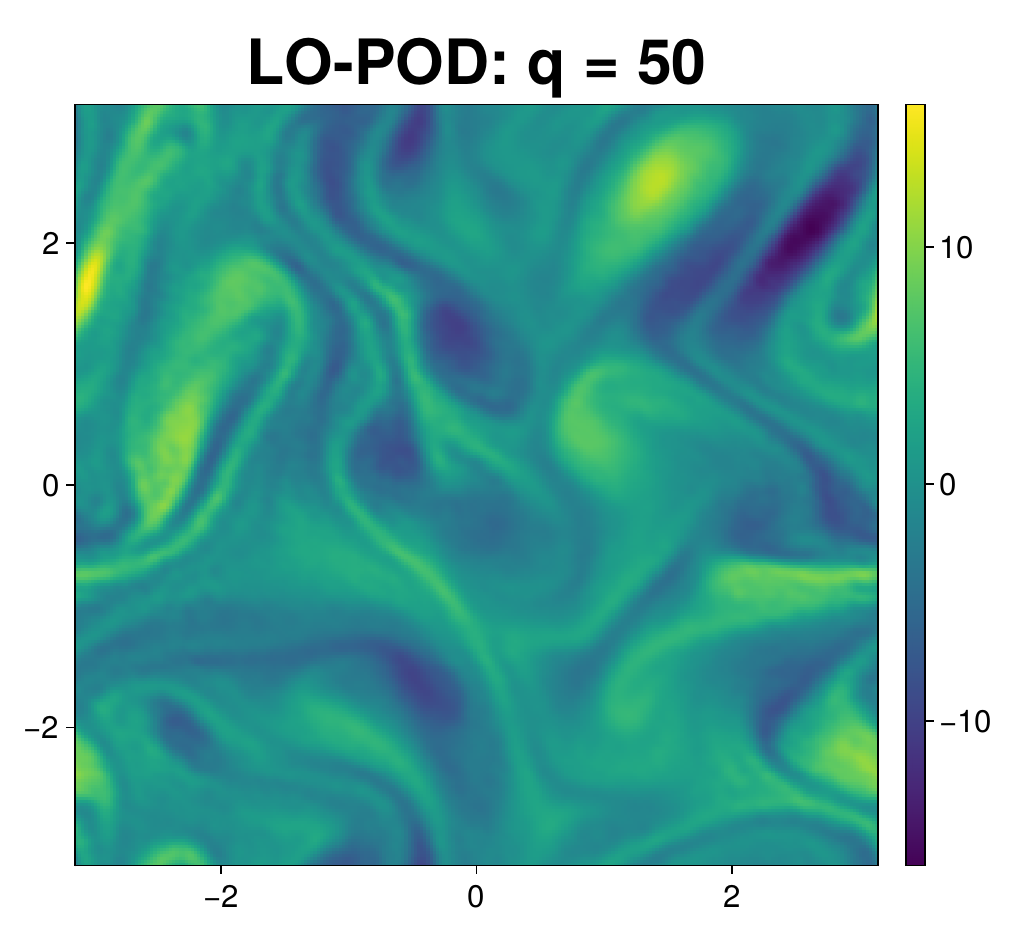}
    \caption{First vorticity snapshot of the validation data projected on the different \gls{POD} bases with varying number of modes per subdomain $q$ for the space-local approaches.}
    \label{fig:2D_NS_heatmaps}
\end{figure}

In Figure \ref{fig:2D_NS_heatmaps}, we display the reconstruction of the first snapshot of the validation set produced by the different bases. 
We depict the reconstructions for both $q = 10$ and $q= 50$ modes per subdomain for the space-local approaches, and all the available \gls{POD} modes, i.e. $r = 50$, for \gls{G-POD}. We find that the space-local approaches capture the flow characteristics quite well, with $8 \times 8 \times 50 = 3200$ \gls{DOF}, whereas \gls{G-POD} does not capture any of the characteristics. Regarding \gls{L-POD}, we find that using a lower number of modes per subdomain $q$ worsens the discontinuities at the subdomain boundaries, whereas \acrshort{LO-POD} suffers from noisier reconstructions for low $q$.

\section{Conclusions}\label{sec:conclusions}

In this work, we presented a novel way of constructing \glspl{ROM} for \glspl{PDE} describing advection-dominated problems. The key properties of our proposed \gls{ROM} are: \textit{sparsity} and \textit{generalizability}  through a space-local \gls{ROM} with overlapping subdomains; \textit{stability} by embedding energy conservation in the \gls{ROM}. The space-local \glspl{ROM} are achieved by building a \gls{POD} basis within each subdomain, as compared to a single common basis for the standard space-global approach. In \cite{farhat_local_POD}, it was shown that such a basis results in a sparse and computationally efficient \glspl{ROM}. \R{comment_1c}\revone{As suggested in \cite{CHUNG2024DDROM}, we modified this approach to generate a common local basis for the entire domain, similar to a finite element basis. This work differs from \cite{CHUNG2024DDROM} in that we applied this to time-dependent \glspl{PDE}, as opposed to steady-state problems. In addition, we introduced overlapping subdomains to ensure a smooth representation of the solution within the space-local \gls{POD} basis.} By generating a common basis, we achieve generalizability in time, as in this way, a feature observed in one subdomain can now be represented in any of the subdomains. In addition, we introduced overlapping subdomains to prevent discontinuities at the subdomain boundaries. To test our methodologies, we made use of the linear advection equation in both 1D and 2D. One of the properties of this system is that the energy is conserved. Our space-local \glspl{ROM} satisfy this property exactly, even when the basis is non-orthogonal (as is the case for the overlapping approach).

We observed that the resulting space-local \glspl{ROM} generalize much better for advection-dominated problems than the standard space-global approach. We also demonstrated that such a space-local approach yields sparser and, consequently, more computationally efficient \glspl{ROM}. In addition, we showed that the space-local \glspl{ROM} also satisfy energy conservation and allow for larger time steps. The latter further decreases the computational cost of the \glspl{ROM}. 
Regarding computational time, we observed an improvement of one to two orders of magnitude with respect to the \gls{FOM} in the 2D test case.
By introducing overlapping subdomains, we demonstrated that we obtain smoother approximations of the solution and require fewer basis functions to represent it compared to the non-overlapping approach. However, this comes at the cost of solving a linear system at each time step. Depending on the application, either overlapping or non-overlapping subdomains might be most suitable. If smoothness of the solution is required, overlapping subdomains might be more suitable, while if computational efficiency is of higher priority, non-overlapping subdomains are likely preferred.

\R{editor_b}\revtwo{To further evaluate the generalizability of the space-local \gls{POD} approaches, we applied it to the 2D incompressible Navier–Stokes equations in a Kolmogorov flow setup. In this case we solely considered the ability of \gls{POD} basis to represent the solution, without building the \glspl{ROM}. \R{editor_c}\revtwo{We consider constructing the actual \gls{ROM} for this non-linear 2D test case outside the scope of this research, as it requires dealing with the divergence freeness condition, and possibly introducing hyperreduction methods for the non-linear term \cite{hyperreduction_1, hyperreduction_2}.} Regarding the representation of the solution space-local methods (\gls{L-POD} and \acrshort{LO-POD}) clearly outperformed \gls{G-POD}, which showed strong overfitting to the training data. In contrast, the local variants maintained low projection errors and accurately reproduced flow features on unseen, temporally decorrelated data. These results indicate that the space-local framework generalizes well for advection-dominated systems and holds strong potential for reduced-order modeling of turbulent flows.}

For future research, we consider constructing a space-local \gls{ROM} for an actual turbulence test case described by the Navier-Stokes equations. To achieve energy conservation, one can make use of the observation presented in \cite{podbenjamin}, where it was shown that carrying out a Galerkin projection on an energy-conserving \gls{FOM} (with quadratic non-linearity) resulted in an energy-conserving \gls{ROM}. However, one important requirement is that the \gls{POD} basis needs to be divergence-free, meaning that the space-local \gls{POD} basis also needs to be divergence-free. Currently, this is still being actively researched in our group. \R{comment_12}\revone{To resolve this, we could take inspiration from \cite{CHUNG2024DDROM} and construct a \gls{POD} basis, not only for the velocity, but also for the pressure \cite{rosenberger2022pressureenergyconsistentromsincompressible}.} Another possible research direction would be to circumvent or speed up the solution of the linear system required in \acrshort{LO-POD}. For this purpose, one could possibly take inspiration from the finite element community \cite{FEM_book}. \R{comment_13_b}\revone{Finally, the application of the space-local \gls{POD} methods to steady-state problems could be investigated. The work in \cite{CHUNG2024DDROM} offers an interesting starting point for this, as they investigated the generalization of the space-local \gls{POD} basis when increasing the size of the domain for steady-state problems. An interesting research avenue would be to compare the space-local \glspl{ROM} to space-global ones, when the spatial domain is kept the same and one of the model parameters is varied.}

\printnoidxglossaries

\section*{CRediT authorship contribution}

\textbf{T. van Gastelen:} Conceptualization, Methodology, Software, Writing - original draft. \textbf{W. Edeling:} Writing - review \& editing. \textbf{B. Sanderse:} Conceptualization, Methodology, Writing - review \& editing, Funding acquisition.

\section*{Data availability}

The code used to generate the training data, construct the \glspl{ROM}, and replicate the presented experiments can be found at  \url{https://github.com/tobyvg/local_POD_overlap.jl}.

\section*{Acknowledgements}

This publication is part of the project ``Unraveling Neural Networks with Structure-Preserving Computing” (with project number OCENW.GROOT.2019.044
of the research programme NWO XL, which is financed by the Dutch Research
Council (NWO)). Part of this publication is funded by Eindhoven University of Technology. Finally, we thank the reviewer and editor for their feedback, enhancing the quality of the article.

\section*{Declaration of generative AI and AI-assisted technologies in the writing process}

During the preparation of this work, the author(s) used ChatGPT to improve language and grammar. After using this tool/service, the author(s) reviewed and edited the content as needed and take(s) full responsibility for the content of the published article.

\bibliographystyle{elsarticle-num}
\bibliography{references}

\appendix

\section{2D Bump kernel}\label{app:2D_kernel}

\begin{figure}[ht]
    \centering
    \includegraphics[width = 0.48\textwidth]{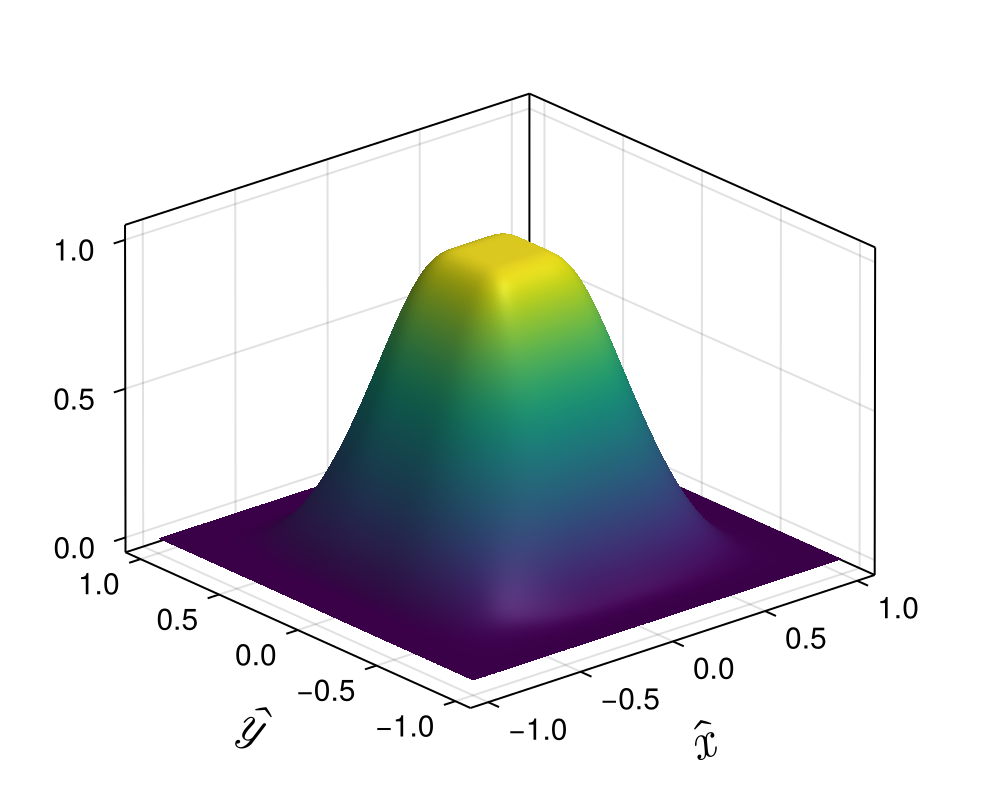}
    \caption{Bump kernel used to construct the \acrshort{LO-POD} basis for the 2D advection test case.}
    \label{fig:2D_kernel}
\end{figure}

Here we discuss the kernel used to construct the \acrshort{LO-POD} basis for the 2D advection test case. This kernel is required to obtain basis functions that smoothly decay to zero at the edge of the subdomain. To achieve this, we use a standard bump function: 
\begin{equation}
    \tilde{\psi}(x) = \begin{cases}
         \exp(\frac{-1}{1 - |x|}) \quad &\text{if } |x| < 1, \\
        0 \quad &\text{elsewhere},
    \end{cases}
\end{equation}
to construct the kernel \cite{lee2012introduction_bump}. This function, along with its derivatives, smoothly decays to zero as $|x|$ approaches $1$.
Next, we apply normalization 
\begin{equation}
    \psi(x) = \begin{cases}
         \frac{\tilde{\psi}(x)}{\tilde{\psi}(x) + \tilde{\psi}(1 - |x|)} \quad &\text{if } |x| < 1, \\
        0 \quad &\text{elsewhere},
    \end{cases}
\end{equation}
to ensure the kernels form a partition of unity, see \eqref{eq:kernel_constraint}.
The 2D kernel is built using a product of two normalized bump functions
\begin{equation}
    k(\mathbf{x}) = \psi(x)\psi(y),
\end{equation}
where $\mathbf{x} = (x,y)^T \in \mathbb{R}^2$ are the spatial coordinates.
For the presentation of this kernel, we exploit the fact that we use a uniform grid. This allows us to conveniently subdivide the domain into square subdomains $\Omega_i$ of the same size. Let $\boldsymbol{\alpha}^i,\boldsymbol{\beta}^i \in \mathbb{R}^2$ denote the coordinates of the bottom-left and top-right vertex. Each subdomain overlaps with its four neighboring subdomains, which share a vertex in the middle of the considered domain. The kernel associated with subdomain $\Omega_i$ is 
\begin{equation}
    k_i(\mathbf{x}) = \psi(\hat{x}_i)\psi(\hat{y}_i),
\end{equation}
where $\hat{x}_i = 2\frac{x - \alpha^i_1}{\beta^i_1 - \alpha^i_1}-1$ and $\hat{y}_i = 2\frac{x - \alpha^i_2}{\beta^i_2 - \alpha^i_2}-1$ are the normalized coordinates. The kernel is depicted in Figure \ref{fig:2D_kernel}.

\section{Choice of subdomains and modes for 2D test case}\label{app:2D_optimization}

\begin{figure}[ht]
    \centering
    \includegraphics[width = 0.48\textwidth]{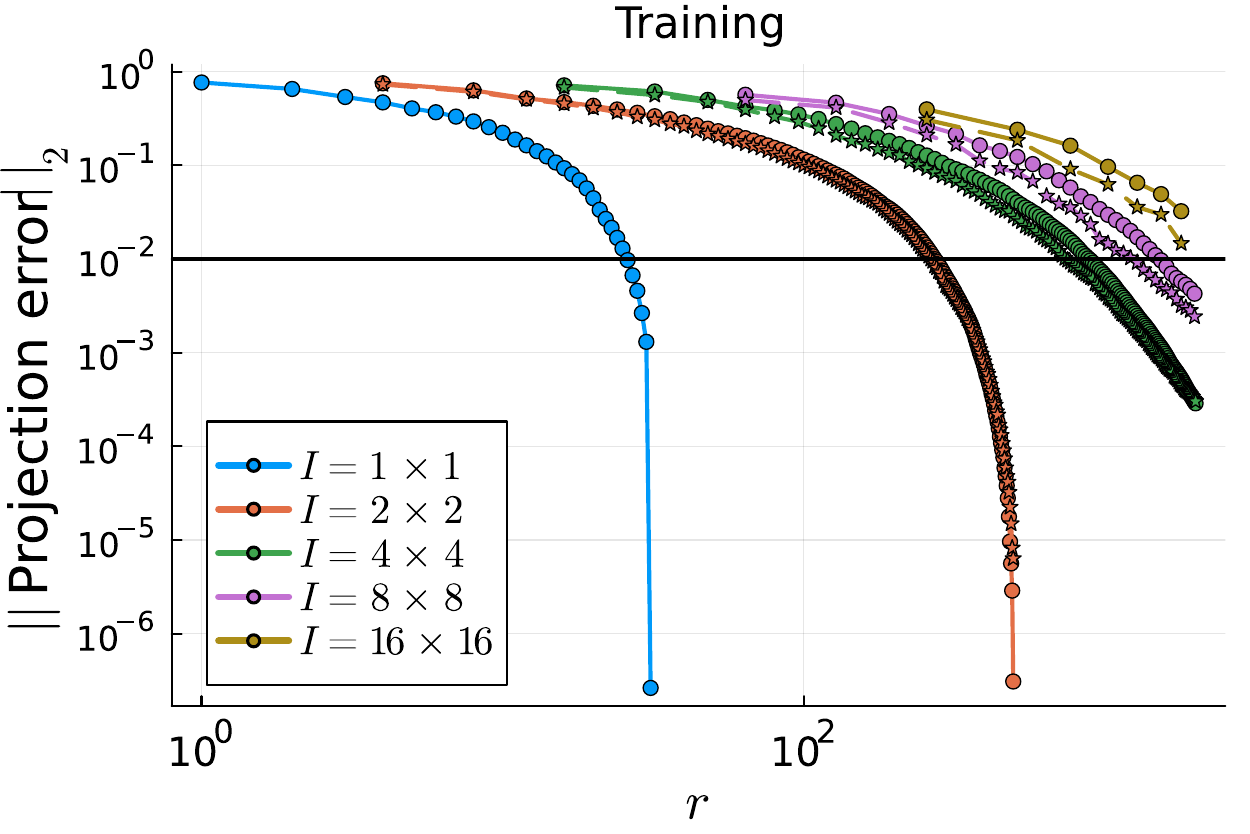}
    \includegraphics[width = 0.48\textwidth]{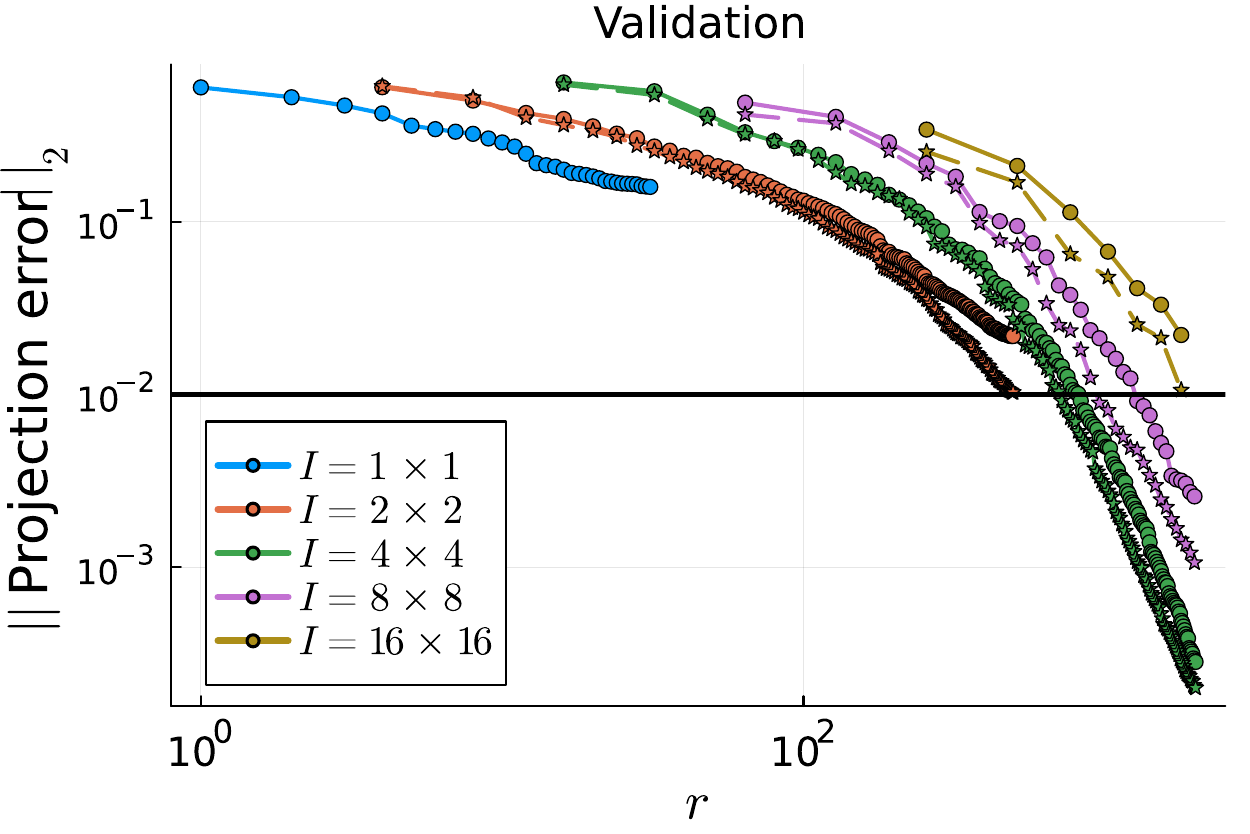}
    \caption{Projection error evaluated over the training (left) and validation (right) data set for the 2D advection test case. $I=1\times 1$ corresponds to \gls{G-POD}. Solid dots correspond to \gls{L-POD} and stars to \acrshort{LO-POD}. The horizontal black line corresponds to a projection error of $10^{-2}$. }
    \label{fig:2D_optimization}
\end{figure}

In this section, we present the projection error on the training and validation sets for the 2D advection test case. This is presented in Figure \ref{fig:2D_optimization}.
Based on these results, we select the number of modes $q$ which first surpasses the projection error of $10^{-2}$ for $I=8\times 8$. This value of $I$ is chosen as it generalizes well, i.e., the performance on the validation set is similar to the training set and the resulting \gls{ROM} is sparser than for $I = 4 \times 4$. For \gls{L-POD} this corresponds to $q = 20$ and for \acrshort{LO-POD} to $q = 15$.

\section{Time step size}\label{app:2D_time_step}

\begin{figure}[ht]
    \centering
    \includegraphics[width = 0.48\textwidth]{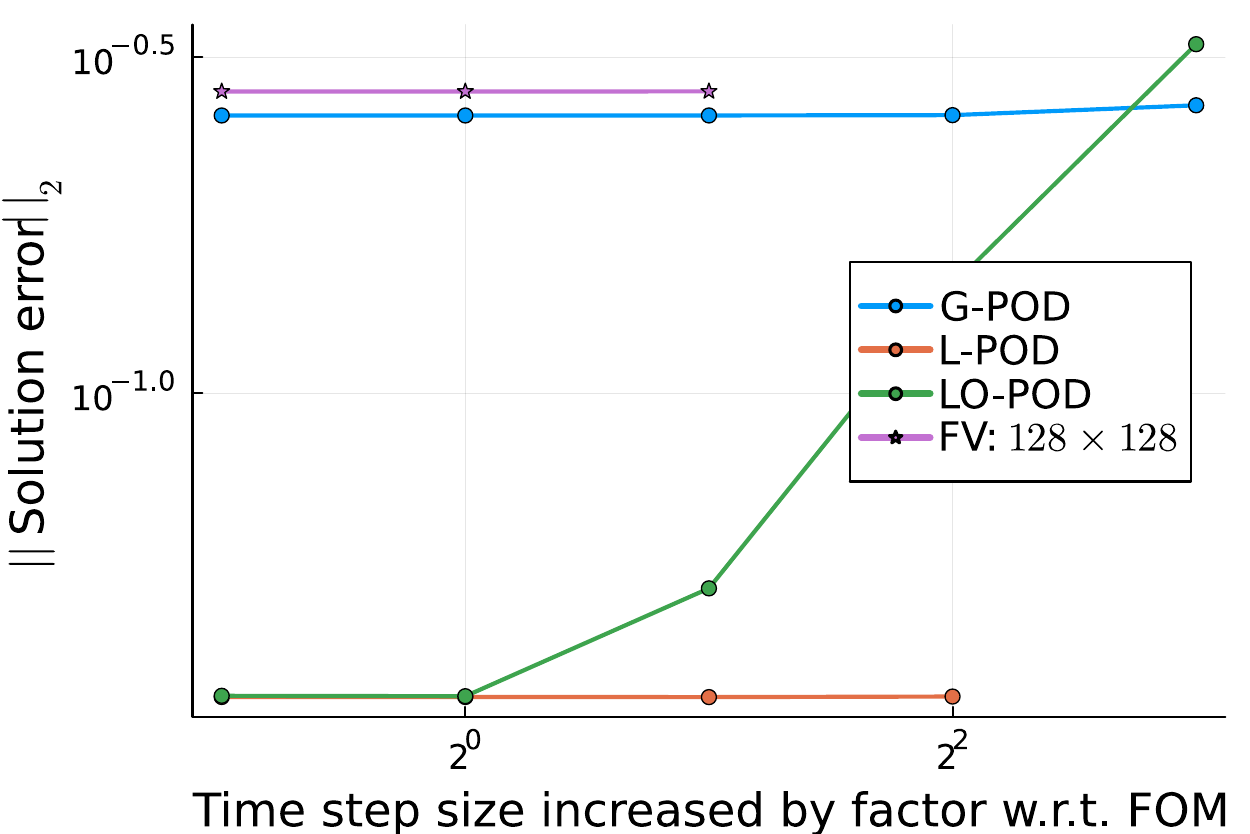}
    \caption{Solution error averaged over the simulation for the 2D advection test case, the \glspl{ROM} for each time step size. Results for a finite volume (FV) discretization on a coarser grid than the \gls{FOM} resolution ($256 \times 256$) are also depicted. Absence of markers indicates an unstable simulation. The \acrshort{LO-POD} \gls{ROM} is integrated using a Crank-Nicolson scheme, as opposed to \gls{RK4} for the others, and is therefore unconditionally stable.}
    \label{fig:2D_time_step}
\end{figure}

In Figure \ref{fig:2D_time_step}, we present the error of \glspl{ROM} integrated using different time step sizes. The maximum time step sizes are selected so that performance degradation or the occurrence of instabilities is avoided. These are presented in Table \ref{tab:ROM_results_overview}.

\end{document}